\documentclass[preprint,review,11pt,authoryear]{elsarticle}

\usepackage{caption}
\usepackage{tikz,graphics,color,fullpage,epsf,caption,subcaption}
\usepackage{float}
\usepackage{algpseudocode}
\usepackage[ruled, commentsnumbered]{algorithm2e}
\usepackage{pdflscape}

\usepackage{amsmath,amssymb}
 \usepackage{amsthm}
 \usepackage{todonotes}
 \usepackage{graphicx}
\usepackage{color}
\definecolor{cornell-red}{RGB}{179,27,27}
\usepackage[hidelinks]{hyperref} 
\usepackage{subcaption}
\usepackage{mathtools}
\usepackage{dsfont}
\usepackage{color}
\usepackage{natbib}
    \usepackage{multirow}
 
\usepackage[left=1in,right=1in,top=1in,bottom=1in]{geometry}

\usepackage{hyperref} 
\usepackage[flushleft]{threeparttable}
\newtheorem{thm}{Theorem}

\newtheorem{prop}[thm]{Proposition}

\theoremstyle{definition}

\newcommand{\black}{\textcolor{black}}  

\theoremstyle{remark}

\newcommand{\blue}{\textcolor{black}}  
\newcommand{\calN}{\mathcal{N}} 
\newcommand{\calT}{\mathcal{T}} 
\newcommand{\calA}{\mathcal{A}} 
\newcommand{\ca}{c_t^{\mbox{\tiny a}}} 
\newcommand{\ch}{c_t^{\mbox{\tiny h}} }
\newcommand{\cp}{c_t^{\mbox{\tiny p}}}
\newcommand{\cu}{c_t^{\mbox{\tiny u}}}
\newcommand{\dU}{\overline{d}}
\newcommand{\dL}{\underline{d}}
\newcommand{\rhoU}{\overline{\rho}}
\newcommand{\rhoL}{\underline{\rho}}

\newcommand{\MU}{\overline{M}}
\newcommand{\ML}{\underline{M}}
\newcommand{\VU}{\overline{V}}
\newcommand{\VL}{\underline{V}}
 \newcommand{\calF}{\mathcal{F}}
 \newcommand{\calR}{\mathcal{R}}

\newcommand{\betaU}{\overline{\beta}}
\newcommand{\betaL}{\underline{\beta}}

\newcommand{\ob}{\pmb{o}}
\newcommand{\zb}{\pmb{z}}

\newcommand{\GammaU}{\overline{\Gamma}}
\newcommand{\GammaL}{\underline{\Gamma}}

\newcommand{\E}{\mathbb{E}}
\newcommand{\Prob}{\mathbb{P}}

\usepackage{lscape}
\usepackage{multirow}
\usepackage{booktabs}
\usepackage{lscape}

\journal{arXiv (under journal review)}

\usepackage{natbib}

\begin{document}
\markboth{Stochastic Optimization Models for Location and Inventory Prepositioning of Disaster Relief Supplies}{\fontsize{10}{10}}

 \begin{frontmatter}

\title{Stochastic Optimization Models for Location and Inventory Prepositioning of Disaster Relief Supplies}
\author[mymainaddress1]{Karmel S. Shehadeh \corref{cor1}}
\cortext[cor1]{Corresponding author.}
 \ead{kshehadeh[at]ehigh.edu, kas720[at]lehigh.edu}
  
 \author[mymainaddress2]{Emily L. Tucker}

 \address[mymainaddress1]{Department of Industrial and Systems Engineering, Lehigh University, Bethlehem, PA,  USA}
\address[mymainaddress2]{Department of Industrial Engineering, Clemson University, Clemson, SC, USA}
\begin{abstract}
\noindent We consider the problem of preparing for a disaster season by determining where to open warehouses and how much relief item inventory to preposition in each. Then, after each disaster, prepositioned items are distributed to demand nodes during the post-disaster phase, and additional items are procured and distributed as needed. There is often uncertainty in the disaster level, affected areas’ locations, the demand for relief items, the usable fraction of prepositioned items post-disaster, procurement quantity, and arc capacity.  To address uncertainty, we propose and analyze two-stage stochastic programming (SP) and distributionally robust optimization (DRO) models, assuming known and unknown (ambiguous) uncertainty distributions. The first and second stages correspond to pre- and post-disaster phases, respectively.  We also propose a model that minimizes the trade-off between considering distributional ambiguity and following distributional belief.  We obtain near-optimal solutions of our SP model using sample average approximation and propose a computationally efficient decomposition algorithm to solve our DRO models. We conduct extensive experiments using a hurricane season and an earthquake as case studies to investigate these approaches computational and operational performance.
\end{abstract}
\end{frontmatter}

\section{Introduction}\label{sec:intro}

\noindent Disasters, such as earthquakes, hurricanes, or tornadoes, can be devastating and unfortunately, frequently occur in disaster-prone areas. Within the United States alone, an average of 13.8 disasters occur annually \citep{noauthor_OCM}. They are hard to predict precisely, and they often strike communities with little warning, leaving devastating impacts on people's lives and infrastructure \citep{rawls2010pre, sabbaghtorkan2020prepositioning}. \black{The ability to quickly meet the urgent need for relief items and provide assistance to disaster-affected populations can be the difference between life and death. To achieve this goal, emergency response experts often seek to preposition emergency relief inventory at strategic locations to distribute when needed.}

\black{A major difficulty in creating an effective and robust prepositioning plan is dealing with uncertainty. We do not know when disasters will hit or what their effects will be, leading to significant variability in post-disaster supply, demand, and road link capacity post-disaster. Such variability is hard to quantify in advance (before a disaster) when prepositioning decisions are made. There is also a need to consider the following trade-offs.  In the immediate aftermath of a disaster, prepositioned items may be the only resources available for distribution, and procuring additional items post-disaster is not easy and is often very expensive. In addition, road networks may be partially functioning or completely damaged, and the items themselves are often subject to higher costs and limited quantities. However, prepositioning too many supplies or putting them in the wrong place can lead to untenably high inventory and transportation costs (see~\ref{Appex:RandomF} for a detailed discussion on these random factors).}

\black{Stochastic optimization approaches can be employed to model uncertainty and help support and optimize prepositioning decisions to better plan and respond to disasters. The existing literature commonly assumes that the exact probability distributions of uncertain factors are known (or there is sufficient and high-quality data to characterize them) and accordingly formulate these problems as two-stage SP models with sample average approximation (SAA). The first and second stages respectively correspond to pre-disaster (e.g., prepositioning of relief inventory) and post-disaster (e.g., distribution of relief items to demand nodes) phases. That is, the SP assumes that the decision-maker has enough data to evaluates the random second-stage costs (e.g., shortage in the aftermath) using the sample average, in which the approximation accuracy improves with the increment of sample size, but the computation often becomes challenging \citep{kim2015guide}.} 

SP remains the state-of-the-art approach to model uncertainty in many application domains.  \textcolor{black}{In the preparedness stage, however, it is unlikely that decision-makers can infer the post-disaster conditions or estimate the actual distributions of random parameters accurately, especially with limited or no information about the situation in the immediate aftermath \citep{altay2014challenges, comes2020coordination, sabbaghtorkan2020prepositioning}.   Even when historical data on past disasters is available, such data is often insufficient to estimate uncertainty distribution accurately, and future disasters often have different characteristics (i.e., distributions) than previous events. In addition, although various types of early warning systems may have been set up, the pre-disasters estimates (e.g., from forecasting and predictive models) of the post-disaster damages and associated demand, for example, are often subject to error and uncertainty. Consequently, the actual costs (benefits) associated with selected deterministic or sample-based decisions (as in the SP approach) tend to be higher (lower) than estimated, causing the decision-maker to experience post-decision disappointment (e.g., significant shortage in the aftermath).}

Alternatively, one can construct an ambiguity set of all distributions that possess some partial information about the random factors before the disaster. This partial information can incorporate available historical data without the assumption of perfect accuracy. Then, using this ambiguity set, one can formulate a two-stage distributionally robust optimization (DRO) problem to minimize the pre-disaster (prepositioning) cost plus the expectation of post-disaster cost over all distributions defined in the ambiguity set. In particular, in the DRO approach, the optimization is based on the worst-case distribution within the ambiguity set. Thus, the distribution of random factors is also a decision variable in DRO.

DRO has received significant attention recently in many application domains due to the following primary benefits. First, as mentioned earlier, while some data on uncertainty may be available  (e.g., from predictive/forecasting models) to support decisions, assuming the decision-maker has complete knowledge of distributions (as in SP) is  unrealistic.  DRO alleviates this assumption by using user-defined ambiguity sets, allowing uncertain variables to follow any distribution defined in this set. Second, intuitive and easy-to-approximate statistics can be used to construct ambiguity sets. For example, decision-makers could estimate the average demand based on their experience in previous disasters or from a forecasting model. Then, one could construct a mean-range ambiguity set of demand distributions, where the range is used to represent the error margin in the estimates. Accordingly, a DRO model could optimize prepositioning decisions against all possible distributions in the ambiguity set that share these mean and range values.

Third, DRO may reflect decision-makers' ambiguity aversion and preference to err on the side of caution  \citep{halevy2007ellsberg, hsu2005neural} as distributionally robust decisions may safeguard performance in adverse scenarios and the direct and indirect costs of operations in the aftermath. Note that even when complete distributional information is available, DRO models often perform well. Forth, various techniques have been developed to derive tractable DRO models \citep{delage2021value, rahimian2019distributionally}, and we have seen many successful applications of DRO to facility location problems (see, e.g.,  \cite{basciftci2019distributionally, saif2020data,  ShehadehSanci, tsang2021distributionally, wang2020distributionally, wang2021two,wu2015approximation}) and disaster management \citep{donmez2021humanitarian, wang2021two}. These benefits suggests that it is worthwhile to consider DRO as an alternative approach for modeling uncertainty. Yet, despite the potential advantages, there are no tractable DRO approaches for the location and inventory prepositioning problem that we study in this paper.

Indeed, there is a trade-off to using different approaches to model uncertainty. For example, SP may provide an excellent basis for the disaster response plan if the model uses the true distribution. However, suppose the SP model uses misspecified distributions or data samples. In that case, the (biased) SP decisions may perform poorly when implemented in practice under the true distribution. In contrast, by focusing on hedging against ambiguity, the optimal DRO decisions may be conservative. A model that minimizes the trade-off between considering distributional ambiguity (DRO pessimism) and following distributional belief (SP optimism) may offer a middle-ground between SP and DRO. Unfortunately, this trade-off has not been studied within the context of the specific location and inventory preposition problem that we study in this paper. 

In this paper, we conduct an analysis to answer the questions of \textit{when} to use each type of modeling approach (SP, DRO, or a trade-off between the two) and \textit{what} is the value of adopting each for a location and inventory prepositioning of disaster relief supplies problem. Specifically, given a set of candidate warehouse locations and different types of relief items, we aim to determine the number of warehouses to open at strategic locations and each relief item's quantity to preposition at each selected location.  Uncertainties considered include the (1) disaster level, (2) locations of affected areas, (3) demand of relief items, (4) usable fraction of prepositioned items post-disaster, (5) procurement quantity, and (6) arc capacity between two different nodes.

We propose a two-stage DRO model for this problem that seeks to find first-stage prepositioning decisions that minimize the sum of the first-stage cost (fixed cost of opening facilities and items acquisitions cost) and the maximum expectation of the second-stage cost (unmet demand, holding, shipping, procurement). We take the expectation over an ambiguity set characterized by the mean and range of the uncertain parameters (1)--(6) (see the motivation and details in Section~\ref{sec3;DRO}). We also consider an SP model that minimizes the expected second-stage cost with respect to assumed known distributions of (1)-(6). In addition, we \textcolor{black}{propose a} model that minimizes the trade-off between considering distributional ambiguity and following distributional belief. 

With the intent of comparing the computational and operational performance of these models, we use a hurricane season and an earthquake as case studies of informational extremes and present extensive computational results. The aim is not to advocate for any of the considered approaches but rather to compare them empirically and theoretically, demonstrating where significant performance improvements may be gained. Finally, we do not claim that our proposed models consider all aspects of disaster preparation and response operations. Instead, we show the pros and cons of each modeling approach and motivate the need to consider multiple approaches when modeling uncertainty for real-world optimization problems.

\subsection{Contributions of the Paper}

We summarize our main contributions as follows. 
\begin{enumerate}\itemsep0em

\item \textbf{Uncertainty Modeling and Optimization Models.}  We propose and analyze a new DRO model for location and inventory prepositioning of disaster relief supplies. We also propose and analyze a new model that minimizes the trade-off between considering distributional ambiguity (DRO pessimism) and following distributional belief (SP optimism). To the best of our knowledge, and per recent surveys of \cite{donmez2021humanitarian, sabbaghtorkan2020prepositioning} and our literature review in Section~\ref{sec2:Lit}, our paper is the first to present an analysis of DRO, SP, and Trade-off modeling approaches for this problem under uncertainty of factors (1)--(6).

\item \textbf{Solution Methods.} We derive equivalent solvable reformulations of the proposed mini-max nonlinear DRO and trade-off models. We propose a computationally efficient decomposition-based algorithm to solve the reformulations. 

\item \textbf{Computational Insights}. We apply the three models to two real-world case studies (Atlantic hurricane season and an earthquake that happened in Yushu County in Qinghai Province, China). We conduct extensive computational analysis to evaluate when it is appropriate to use each approach to modeling uncertainty. Our results demonstrate the (1) superior post-disaster operational performance of the DRO decisions under various distributions compared to SP decisions, (2) trade-off between considering distributional ambiguity and following distributional belief, and (3)  computational efficiency of our approaches. More generally, our results draw attention to the need to model the distributional ambiguity of uncertain problem data in strategic real-world stochastic optimization problems. 

\end{enumerate}

\subsection{Structure of the Paper}
\noindent The reminder of this paper is structured as follows.  In Section~\ref{sec2:Lit}, we review relevant literature. In Section~\ref{sec3.1:statement}, we detail our problem setting. In Section~\ref{sec3:TSM}, we present our proposed SP model. In Sections~\ref{sec3;DRO} and \ref{sec:Trade}, we present our proposed DRO model and trade-off model, respectively. In Section~\ref{sec:SolutionApproches}, we present a decomposition algorithm to solve the DRO model and a Monte Carlo Optimization Procedure to solve the SP model. In Section~\ref{sec:Exp}, we conduct extensive numerical experiments using a hurricane season and an earthquake as case studies. Finally, we draw conclusions in Section~\ref{sec:conclusion}.

\section{Literature Review}\label{sec2:Lit}

\noindent Disaster operations management consists of four phases: mitigation, preparedness, response, and recovery  \citep{altay2006or, sabbaghtorkan2020prepositioning}. Mitigation and preparedness activities occur before a disaster and typically include emergency inventory prepositioning, facility location, and transportation decisions \citep{alkaabneh2020benders, aboolian2013efficient, altay2013capability, tucker2020incentivizing, duran2011pre, sanci2019integrating,lee2009modeling,qi2010effect, shen2011reliable, sheu2007emergency, toregas1971location, yushimito2012voronoi}. Response and recovery are post-disaster relief actions. One essential component is to deliver emergency supplies \citep{barbarosolu2004two, holguin2013appropriate, sheu2010dynamic,tzeng2007multi}. Recent literature has highlighted the need for integrating location, inventory prepositioning, and delivery  decisions for disaster relief as well as assess humanitarian response capacity (see, e.g., \cite{acimovic2016models,arnette2019risk,  balcik2019collaborative,  dufour2018logistics, duran2011pre, charles2016designing, jahre2016integrating, moline2019approaches, mccoy2011efficient, rawls2010pre, sabbaghtorkan2020prepositioning, salmeron2010stochastic,ni2018location, velasquez2020prepositioning, Abazari2021}).

Most of the original, pioneering inventory prepositioning models assume perfect knowledge about the post-disaster conditions (e.g., demand for relief items, road conditions). In reality, it is unlikely that decision-makers be able to predict the exact post-disaster conditions. Thus, a naive, deterministic approach that uses point estimates for random parameters will likely produce suboptimal prepositioning decisions.  To address this, several authors have studied inventory prepositioning under uncertainty. In what follows, we review some papers that are most relevant to this work: papers that use stochastic optimization (specifically SP, Robust Optimization (RO), and DRO) for prepositioning in disaster response operations. For comprehensive surveys, we refer to \citep{altay2006or, anaya2014relief, galindo2013review, gupta2016disaster,Boonmee2017, sabbaghtorkan2020prepositioning}.


Most of these studies employ two-stage SP to model uncertainty, assuming that random parameters' probability distributions are fully known. The first stage of these models consists of the location and amount of prepositioned supplies, and the second-stage consists of transportation decisions. Optimization criteria of these studies include minimizing the expected total cost \citep{chang2007scenario,doyen2012two,  mete2010stochastic, rawls2010pre},  minimizing trade-off of total cost with Conditional Value-at-Risk \citep{noyan2012risk}, minimizing expected response time \citep{duran2011pre}, minimizing expected casualties \citep{salmeron2010stochastic},  maximizing coverage by Social Vulnerability Index \citep{Alem2021}, and maximizing the expected satisfied demand or minimizing expected shortage of relief items \citep{balcik2008facility}.

While SP is a powerful modeling approach for modeling uncertainty, it suffers from the following shortcomings. First, to formulate an SP, one needs access to all possible scenarios and their probabilities, which is not realistic for inventory prepositioning due to the unpredictable nature of disasters and effects \citep{condeixa2017disaster, salmeron2010stochastic, velasquez2020prepositioning}.  If we calibrate an SP to a particular training dataset, the resulting decision policy may have disappointing out-of-sample performance under an unseen data set from the same population. Second, SP approaches suffer from the ``curse of dimensionality,'' and they are often intractable. 

RO assumes a complete ignorance about the probability distribution of uncertain parameters. Instead, RO assumes that uncertain parameters reside in a so-called ``uncertainty set'' of possible outcomes \citep{bertsimas2004price, ben2015deriving, soyster1973convex}, and optimization is based on the worst-case scenario occurring within the uncertainty set. Only a few papers have employed RO for inventory prepositioning.  \cite{zokaee2016robust}'s RO model incorporated the uncertainty associated with demand, supply, and cost parameters, assuming that these uncertain parameters are independent and bounded random variables. They controlled the level of conservativeness using a budget of uncertainty \citep{bertsimas2004price} and reformulated their model as a mixed-integer program.

Two recent papers have also considered RO to preposition inventory \citep{ni2018location,velasquez2020prepositioning}. \cite{ni2018location} proposed a two-stage RO model to minimize prepositioning (first-stage) cost and (second-stage) shortage costs for one disaster.  The uncertainty set included three parameters: demand, usable fraction of prepositioned supplies, and arc capacity. They used an uncertainty budget to control the size of the uncertainty set.   \cite{ni2018location}  used off-the-shelf solvers for small instances and proposed a Benders decomposition approach that can solve larger instances.  The model did not consider the risk of multiple disasters, disaster level, or the possibility of procuring additional relief supplies post-disaster. We incorporate these random factors in our models in addition to the demand and usable fraction of prepositioned supplies.

\cite{velasquez2020prepositioning} proposed a two-stage RO model that hedged against multiple disasters (e.g., hurricanes) and incorporated uncertainty of affected areas. In the first stage, the model determines the location and amount of prepositioned relief supplies before any disaster occurs.  In the second stage, prepositioned relief items are distributed, and additional items are procured as needed. The objective is to minimize the total cost of prepositioning and distributing disaster relief supplies. To solve their model, \cite{velasquez2020prepositioning} proposed a column-and-constraint generation algorithm. The model assumed that if a disaster occurs at a particular location, the demand for a relief item and the usable fraction of prepositioned  supplies at that location are equal to their \textit{mean} values (which depends on the disaster type).  They assumed that procurement quantity and arc capacities post-disaster are deterministic. In our paper, we model uncertainty of procurement quantity and arc capacities post-disaster.

Road and facility vulnerability have been concurrently considered.
For example, \cite{aslan2019pre} proposed a two-stage SP model to design a multi-echelon humanitarian response network. They emphasized that pre-disaster decisions of warehouse location and item prepositioning are subject to uncertainties both in relief item demand and vulnerability of roads and facilities following the disaster. In the second stage (after the disaster), relief transportation is accompanied by simultaneous repair of blocked roads, which gradually increases the network's connectivity at the same time. In addition to considering random arc capacity between two nodes, our paper considers an additional aspect not considered in prior literature. That is, the level of disaster at each node and the associated damage to the link connecting each pair of nodes (see Section~\ref{sec3:TSM} for more details).

An alternative paradigm for modeling uncertainty is DRO. It aims to unify SP and RO while overcoming their drawbacks \citep{rahimian2019distributionally, saif2020data}. In DRO, we assume that the distribution of uncertain parameters resides in a so-called ``\textit{ambiguity set}'' which is a family of all possible distribution of uncertain parameters characterized through some known properties of uncertain parameters \citep{esfahani2018data}. The optimization is based on the worst-case distribution within the ambiguity set. There are three primary benefits to using DRO to model uncertainty. First, DRO alleviates the unrealistic assumption of the decision-maker's complete knowledge of the distribution governing the uncertain problem data. Second, DRO models are often more computationally tractable than their SP and RO counterparts. Finally, DRO avoids the well-known over-conservatism and the poor expected performance of RO and allows for better utilization of the available data. To construct the ambiguity sets, one can use easy-to-compute information such as mean values and ranges of random parameters and derive tractable DRO models that better mimic reality.  We refer to \cite{rahimian2019distributionally} for a comprehensive survey of DRO.

From \cite{sabbaghtorkan2020prepositioning}, we observe that despite the potential advantages, there are no tractable DRO approaches for the specific location and inventory preposition problem that we study in this paper. In fact, \cite{wang2021two} is the only DRO approach for integrated facility location and inventory prepositioning of disaster relief supplies. However, \cite{wang2021two} assume that demand is the only random factor and accordingly define (box and polyhedral) ambiguity sets to model the imprecise probability information of demand. Their model seek prepositioning decisions that minimizes the total cost of warehouse construction cost, resource preposition storage cost, and transportation cost after the disaster. \cite{wang2021two}'s DRO model does not account for (1) the uncertainty of factors (2)-(6) mentioned earlier, (2) post-disaster decisions and costs such as procurement quantity and its cost, holding cost, etc., and (3) trade-offs between SP and DRO approaches.

In this work, we develop a tractable DRO model to address the practical problem of limited distributional information and compare this model to an SP model. We also expand the considered problem by incorporating uncertainty and distributional ambiguity for parameters previously assumed to be deterministic. These include the maximum quantity available to be procured post-disaster, arc capacities, demand, and the usable post-disaster fraction of prepositioned relief items (assumed deterministic in \cite{velasquez2020prepositioning} and not included \cite{wang2021two}). We also model the possibility of procuring an additional amount of relief supplies post-disaster and include the procurement cost in the second-stage objective, in contrast to \cite{ni2018location}'s  recent RO approach for disaster inventory prepositioning. We study the trade-off between considering distributional ambiguity and following distributional belief. This trade-off has not been studied within the context of the specific location and inventory preposition problem that we study in this paper.

\section{\textbf{Problem Description}}\label{sec3.1:statement}

\noindent Let us now introduce a generic description of the problem. We will later (Sections~\ref{sec3:TSM} and \ref{sec3;DRO}) introduce the two-stage stochastic optimization models based on this description. Then, in the computational study, where we use the models to prepare for particular disaster types, we make assumptions about some of the sets, parameters, and variables related to the disaster type.

We consider a directed graph $G (\calN, \calA)$, where nodes $i\in \calN$ are candidate locations for warehouses and arcs $a \in \calA$ represent roads. We assume the set of candidate locations and the demand nodes are the same, without loss of generality. That said, our models can be used with two different sets of candidate locations and demand nodes.

First, let us introduce the parameters and decision variables defining our first stage (pre-disaster). For each $i \in \calN$, we define binary variable $o_i$ that equals 1 if a facility is open at location $i$, and is 0 otherwise. For each $i \in \calN$, we define a  non-negative parameter $f_i$ as the fixed cost of opening a facility at location $i$. We define non-negative parameter $S_i$ as the storage capacity of location $i \in \calN$. We consider a set of relief items $\calT$, and  define non-negative parameter $s_t$ as the capacity needed to store each unit of relief item $t\in \calT$. For each $t \in \calT$ and $i \in \calN$, we define non-negative continuous variable $z_{t,i}$  to represent the amount of relief item $t$ prepositioned at location $i$. We let non-negative parameter $\ca$ represent the unit acquisition cost of relief item $t\in \calT$. We can serve the demand for relief item $t \in \calT$ at location $i\in  \calN$ with prepositioned supplies or with supplies procured post-disaster. We assume that there is no limit on the relief supplies available for prepositioning in the first stage, as in prior studies \citep{sabbaghtorkan2020prepositioning, velasquez2020prepositioning}. Indeed, one can easily add an additional constraint to the first stage of the proposed models to enforce a limit on the relief supplies available for prepositioning.

Let us now introduce additional sets, parameters, and variables defining our second stage. We consider  $|L|$ different types, or levels, of disasters (e.g., major, minor). For all $l \in L$ and $i \in I$, we define $q_{i,l}$ as a 0-1 random variable, which equals to one if a disaster of level $l$ occurs at location $i$, and $q_{i,l}=0$ otherwise. The level of disaster affects: the demand for relief items at each location, how much of the prepositioned inventory is available (some may be destroyed), how much of each relief item can be ordered after the disaster, and the capacity of roads connected to the affected location. As in prior literature \citep{ni2018location,rawls2010pre, sabbaghtorkan2020prepositioning, velasquez2020prepositioning}, in the computational study, we restrict ourselves to a one disaster per disaster-prone node, i.e., $\sum \limits_{l \in L} q_{i,l}\leq 1$, $\forall i \in I$.    However, our models consider the possibility that, for example, one hurricane can impact multiple locations.  Moreover, our models and solution methods can handle preparation for multiple disasters, such as an entire hurricane season.

For each $t\in \calT$, $i\in \calN$, $l \in L$, we let $d_{t,i,l}$ represent the random demand for relief item $t$ if a disaster of type $l$ occurs at node $i$. For all $i \in I$ and $l \in L$, we let $\rho_{i,l}$ represent the fraction of relief supplies prepositioned at location $i$ that remains usable after a disaster of type $l$ occurs. As mentioned earlier, it is possible to procure additional relief items after the disaster strikes in some cases. We assume that relief procurement is limited and more expensive given that supplies are procured with short notice, and price hikes are common \citep{velasquez2020prepositioning}. We assume that the maximum order quantity available post-disaster and the capacity of arc $(i,j)$ are random. For all $t\in \calT$, $i\in \calN$, and $l \in L$, we let random parameter $M_{t,i,l}$ represent the maximum order quantity of relief item $t$ available to procure after a disaster of type $l$ at location $i$. For all $(i,j) \in \calA$ and $l \in L$, we let random parameter $V_{i,j,l}$ represent the random capacity of arc $(i,j)$ after a disaster of type $l$ occurs at either $i$ or $j$. If a disaster does not occur at either, then the capacity of arc $(i,j)$ is equal to its nominal value, $\hat{V}_{i,j}$.

For all $t \in \calT$ and $i \in \calN$, we define the non-negative continuous decision variable $y_{t,i}$ as the amount of relief item $t$ procured post-disaster at location $i$. For all $(i,j) \in \calA$, we define the non-negative continuous decision variable $x_{i,j}^t$ as the amount of relief item $t$ shipped through arc ($i,j$), i.e., flow quantity across link ($i,j$). For all $t \in \calT$ and $i \in \calN$, we define non-negative continuous decision variables $e_{t,i}$ and $u_{t,i}$ to respectively represent the quantity of unused inventory of item $t$ at location $i$ and quantity of unsatisfied demand of relief item $t$ at location $i$.  The non-negative parameters $\ch$, $\cp$, and $\cu$ represent the unit holding cost of relief item $t$, unit cost of procuring relief item $t$ post-disaster, and unit penalty cost of unmet demand for relief item for all $t\in \calT$, respectively. Finally, the  non-negative parameter $c_{i,j}^t$ represents the unit transportation cost of relief item $t\in \calT$ using arc ($i,j)\in \calA$.

In the first stage, we decide (1) number of facilities to open and their location ($\pmb o$), and (2) amount of each relief item to store at each open facility ($\pmb z$). Additional relief items are procured and distributed in the second stage. The quality of these \black{prepositioning} decisions is a function of the (1) fixed costs of opening facilities and item acquisition (first-stage), (2) cost of unmet demand (second-stage), (3) holding cost (second-stage), (4) procurement cost (second-stage), and (5) shipment cost (second-stage).  We summarize the notation in Table~\ref{table:notation}.
 

 \noindent \textbf{Additional notation.}   For $a,b \in \mathbb{Z}$, we define $[a]:=\lbrace1,2,\ldots, a \rbrace$ and $[a,b]_\mathbb{Z}:=\lbrace c \in \mathbb{Z}: a \leq c \leq b \rbrace$, i.e., $[a,b]$ represent the set of running integer indices $\{a, a+1, a+2, \ldots, b \}$. 

\begin{table}[t!]  
\small
\center
   \renewcommand{\arraystretch}{0.5}
  \caption{Notation.}
\begin{tabular}{ll}
\hline
\multicolumn{2}{l}{\textbf{Sets and Indices}} \\
$\calN$ & set of nodes \\
$I$ & set of potential  facility (warehouse) sites to store the prepositioned emergency \\
& supplies, $I \subseteq \calN$ \\
$\calT$ & set of relief supply types \\
$\calA$ & set of arcs\\
$L$ & set of disaster levels\\
\multicolumn{2}{l}{\textbf{Parameters}} \\
$q_{i,l}$ & binary random variable that equals 1 if a disaster of type $l\in L$ occurs\\
& at location $i$ and 0 otherwise \\
$\ca$& unit acquisition cost for \black{prepositioning} relief item $t$\\
$\ch$& unit holding cost for relief item $t$\\
$\cp$&  unit cost of  procuring relief item $t$ post-disaster\\
$\cu$& unit penalty cost of unmet demand for relief item $t$ \\
$c_{i,j}^t$ & unit transportation cost of relief item $t$ using arc ($i,j$) \\
$f_i$& fixed cost of opening facility $i$ \\
$S_i$ & storage capacity at location $i$ \\
$s_t$ & unit storage space required for relief item $t$\\
$v_t$ & arc capacity required to transport a unit of relief item $t$\\
$M_{t,i,l}$ & random maximum order quantity available to procure post-disaster for relief item $t$\\
& at location $i$\\
$\rho_{i,l}$ & random fraction of relief items prepositioned at location $i$ that remains usable  \\
& after a disaster of type $l$ \\ 
$ \hat{V}_{i,j}$ & nominal capacity of arc ($i,j$), i.e.,  if no disaster occurs\\
$ V_{i,j,l} $ & random capacity of arc ($i,j$) if a disaster of type $l$ occurs\\ 
$d_{t,i,l}$ & random demand for relief item $t$ if a disaster of type $l$ occurs at location $i$\\
 \multicolumn{2}{l}{\textbf{Decision Variables } } \\
 $o_i$ &  binary variable that equals 1 if a facility is open at location $i$ and 0 otherwise \\ 
$z_{t,i}$ & amount of relief item $t$ prepositioned at location $i$ \\
$x_{i,j}^t$ &amount of relief item $t$ shipped through arc ($i,j$), i.e., flow quantity across link ($i,j$) \\
$y_{t,i}$ & amount of relief item $t$ procured post-disaster at location $i$ \\
$e_{t,i}$ & quantity of unused inventory of item $t$ at location $i$\\
$u_{t,i}$ &  quantity of unsatisfied demand of relief item $t$ at location $i$\\
\hline
\end{tabular}\label{table:notation}
\end{table} 
\section{\textbf{Two-stage SP Model}}\label{sec3:TSM}

\noindent \textcolor{black}{In this section, we present our proposed two-stage SP formulation that assumes that the probability distribution $\mathbb{P}$ of random parameters $\xi:=[\pmb{q, d, \rho, M, V}]^\top$ is known. A complete listing of the parameters and decision variables of the model can be found in Table~\ref{table:notation}. Using this notation, we formulate the following two-stage SP model.}
 \allowdisplaybreaks
\begin{subequations}\label{TSM}
\begin{align}
\upsilon=& \min \limits_{\ob,\zb} \left \lbrace \sum_{i \in I} f_i o_i + \sum \limits_{i \in I} \sum\limits_{t \in \calT} \ca z_{t,i}  +  \E[Q (o,z,  \xi)] \ \right\rbrace \label{ObjTSM}\\
& \text{s.t.} \ \ \sum \limits_{t \in \calT} s_t z_{t,i} \leq S_i o_i, && \forall i \in I  \label{First-stage:C2}\\ 
& \qquad  \ o_i \in \lbrace 0, 1 \rbrace, \ z_{i,t \geq 0} && \forall i \in I , t \in \calT \label{First-stage:C3}
\end{align}
\end{subequations} 

 where for any feasible first-stage decisions ($\ob,\zb$) and a given joint realization of uncertain parameters $\xi:=[\pmb{q, d, \rho, M, V}]^\top$
 \allowdisplaybreaks
\begin{subequations}\label{2ndstage}
\begin{align}
 Q(\ob,\zb, \xi ):= \min_{\textcolor{black}{\pmb{y, u, e, x}}} &     \ \ \Big \{ \sum \limits_{i \in I} \sum \limits_{t \in \calT} (\cp y_{t,i}+\cu u_{t,i} +\ch e_{t,i} ) + \sum \limits_{t \in \calT} \sum \limits_{(i,j) \in \calA} c_{i,j}^{\mbox{\tiny t}} x_{i,j}^t \Big\} \label{Obj:Q}\\  
\text{s.t.}   &  \sum \limits_{j: (j,i) \in \calA} x_{j,i}^t -\sum \limits_{j:(i,j) \in \calA} x_{i,j}^t+y_{t,i} -e_{t,i}+u_{t,i}  \nonumber  \\
&  \ = \sum_{l \in L}q_{i,l}  \big (d_{t,i,l}-\rho_{i,l}z_{t,i} \big)-(1-\sum_{l \in L}q_{i,l}) z_{t,i}, \qquad  \forall t \in \calT , i  \in I,\label{Const1:Q}\\
& \qquad \  \  y_{t,i} \leq  \sum_{l \in L}q_{i,l}M_{t,i,l}, \qquad \forall t \in \calT, i \in I,\label{Const2:Q} \\
& \qquad \  \sum \limits_{t \in \calT} v_t x_{i,j}^t \leq \sum_{l \in L}q_{i,l} V_{i,j,l}+ (1- \sum_{l \in L}q_{i,l})\hat{V}_{i,j}, \qquad  \forall (i,j) \in \calA, \label{Const3:Q}\\
& \qquad \  \sum \limits_{t \in \calT} v_t x_{i,j}^t \leq \sum_{l' \in L}q_{j,l'} V_{i,j,l'}+ (1- \sum_{l' \in L}q_{j,l'})\hat{V}_{i,j}, \qquad  \forall (i,j) \in \calA, \label{Const32:Q}\\
& \qquad \  y_{t,i} \geq 0, \ e_{t,i} \geq 0, \ u_{t,i}\geq 0, \ \ x_{i,j}^t \geq 0, \qquad \qquad \forall t \in \calT, i \in (i,j)\in \calA.\label{Const4:Q}
\end{align} 
\end{subequations}

The objective function \eqref{ObjTSM} minimizes the sum of the fixed cost of opening facilities (first term) and prepositioning cost (second term) plus the total expected cost of response operations (third term). The expected cost of response is the sum of the costs of procurement (first term in \eqref{Obj:Q}), shortage (second term in \eqref{Obj:Q}), holding (third term in \eqref{Obj:Q}), and transportation (fourth term in \eqref{Obj:Q}). 

Constraints~\eqref{First-stage:C2} enforce the storage capacity for prepositioning relief supplies. Constraints~\eqref{First-stage:C3} define the feasible range of variables $\ob$ and $\zb$. In the second-stage, constraints \eqref{Const1:Q} enforce flow conservation. Constraints~\eqref{Const2:Q} limit the relief supplies procured post-disaster. If no disaster occurs at location $i$ (i.e., $\sum_{l \in L}q_{i,l}=0$), then all of the prepositioned items (if any) remain usable. Constraints~\eqref{Const3:Q}--\eqref{Const32:Q} enforce arc capacities. \textcolor{black}{In addition to considering random arc capacity, these constraints consider an additional aspect not considered in prior literature: the level of disaster at each node and the associated damage to the link connecting each pair of nodes.} If there is not a disaster at either end of arc $(i,j)$, then the arc capacity is not affected, and it is equal to its nominal value, $\hat{V}_{i,j}$. If a disaster hits node $i$ or node $j$, the arc capacity is $V_{i,j,l}$ or $V_{i,j,l'}$, respectively. If a disaster hits both nodes, the arc capacity is defined as $\min \{ V_{i,j,l}, V_{i,j,l'} \}.$ \textcolor{black}{Finally,} constraints~\eqref{Const4:Q} specify feasible ranges of the decision variables.

Note that the recourse problem is a feasible and bounded linear program for any feasible first-stage decision  ($\ob,\zb$) and a given joint realization of uncertain parameters $\xi$. Thus, we have relatively complete recourse.

\section{\textbf{Two-stage DRO model}}\label{sec3;DRO}

\noindent  The distributions of post-disaster parameter uncertainty (e.g., demand for relief supplies) may be difficult to estimate before a disaster occurs. Historical data on past disasters may be either unavailable or insufficient to model the correct distributions of uncertainty. Future disasters often have different characteristics (e.g., distributions) than previous events. Therefore, in this section, we propose a DRO model \textcolor{black}{that} does not assume that the probability distribution of $\mathbb{P}$ of random parameters is exactly known.

\textcolor{black}{First, let us introduce some additional sets and notation defining our proposed DRO model.} We assume that we know or can approximate the mean values and support (i.e., $\overline{\text{upper}}$ and \underline{lower} bound) of the random parameters. Mathematically, we consider support $\calR= \calR^q \times \calR^{\mbox{\tiny d}} \times  \calR^{\tiny \rho}   \times \calR^{\mbox{\tiny M}} \times \calR^{\mbox{\tiny V}}$, where $\calR^q , \calR^{\mbox{\tiny d}} , \calR^{\tiny \rho}  , \calR^{\mbox{\tiny M}},$ and $\calR^{\mbox{\tiny V}}$ \textcolor{black}{defined in \eqref{Support}} are the supports of random parameters $\pmb{q,\  d, \ \rho, \ M},$ and $\pmb{V}$, respectively. The assumption of known mean and support is motivated by the fact that decision-makers often estimate the values of disaster-related uncertainty using historical data.  Advancements in this sector have included developing forecasts for the occurrence and effects of future disasters, and drones and satellite technology are increasingly used for damage assessments. In addition, engineering methods could be used to estimate the mean and range of the potential damage levels ($\rho$) as well as the level of damage for arc capacities.  If estimates are available or can be approximated, the range could represent the potential error margin in these estimates.
 \begin{subequations}\label{Support}
\begin{align}
&  \calR^q:=\{0,1\}^{N}, \qquad \calR^{\mbox{\tiny d}} :=\left\{ \pmb{d} \geq 0: \begin{array}{l} \dL_{t,i,l} \leq d_{t,i,l}\leq \dU_{t,i.l}, \forall (t, i, l) \end{array} \right\},\\
& \calR^{\rho} :=\left\{ \pmb{\rho} \geq 0: \begin{array}{l} \rhoL_{i,l} \leq \rho_{i,l} \leq  \rhoU_{i,l}, \textcolor{black}{\forall (i , l)} \end{array} \right\}, \\
& \calR^{\mbox{\tiny M}} :=\left\{\pmb{M} \geq 0: \begin{array}{l} \ML_{t,i,l} \leq M_{t,i,l} \leq \MU_{t,i,l}, \textcolor{black}{\forall (t, i)} \end{array} \right\},\\
& \calR^{\mbox{\tiny V}} :=\left\{ \pmb{V} \geq 0: \begin{array}{l} \VL_{i,j,l} \leq V_{i,j,l} \leq \VU_{i,j,l}, \forall (i,j) \in \calA, \forall l \end{array} \right\}.
\end{align}
\end{subequations}

We let $\mu^q, \mu^{\mbox{\tiny d}}, \mu^{\tiny \rho}, \mu^{\mbox{\tiny M}},$ and $\mu^{\mbox{\tiny V}} $ represent the mean values of $\pmb{q, \ d, \ \rho, \ M}$, and $\pmb{V}$, respectively. \textcolor{black}{In addition, we define $\mathcal{P}=\mathcal{P}(\calR)$ as the set of all probability measures on $(\calR,\mathcal{B})$, where $\mathcal{B}$ is the Borel $\sigma$-field on $\calR$. Elements in $\mathcal{P}(\calR)$ can be viewed as probability measures induced by the random vector~$\xi$. Finally, given  $\Prob \in \mathcal{P} (\cdot)$, we define $\E_\Prob$ as the expectation under $\Prob$ and we let $ \mu:=\E[\xi]_\Prob=[\mu^q, \mu^{\mbox{\tiny d}}, \mu^{\tiny \rho}, \mu^{\mbox{\tiny M}}, \mu^{\mbox{\tiny V}}]^\top$. Using this notation, we construct the following mean-support ambiguity set}.
\begin{align}\label{eq:ambiguity}
\calF(\calR, \mu) := \left\{ \Prob\in \mathcal{P}(\calR)\middle|  \begin{array}{l}  \E_\Prob[\xi] = \mu \end{array} \right\}.
\end{align}

 \textcolor{black}{Note that we do not consider higher moments of random parameters for three primary reasons. First, the mean and range are intuitive statistics that a decision-maker may approximate and change in the model (e.g., the mean may be estimated from limited data or approximated by subject matter experts, and the range may represent the error margin in the estimates). Second, it is not straightforward for decision-makers to approximate or estimate higher moments \citep{comes2020coordination}.  Third, various studies have demonstrated that incorporating higher moments in the ambiguity set often undermines the computational tractability of DRO models and therefore their applicability in practice. In contrast, DRO models for real-life problems based on mean-range ambiguity sets often allow for tractable reformulations and solution methods. Indeed, as we will show later, using $\calF(\calR, \mu) $ allow us to derive a computationally tractable reformulation and solution methodology.}

Using ambiguity set $\calF(\calR, \mu)$, we formulate our \textit{DR location and inventory \black{prepositioning} of disaster response operations} problem as  the following min-max problem:
\begin{align}
& \min_{\textcolor{black}{\ob, \zb} } \  \left \lbrace \sum_{i \in I} f_i o_i + \sum \limits_{i \in I} \sum\limits_{t \in \calT} \ca z_{t,i}  + \sup \limits_{\Prob\in \calF(\calR, \mu)} \E_\Prob[Q(\ob,\zb,\xi)] \right \rbrace. \label{DRModel}
\end{align}

The DRO model in \eqref{DRModel} aims to find the prepositioning decisions ($o,z$) that minimize the sum of the fixed facility costs, prepositioning costs, and the worst-case expected cost of response operations, over all distribution within the ambiguity set  $\calF(\calR, \mu)$.
\subsection{\textbf{Reformulation of the DRO model }}\label{sec3:Reformulation}

\noindent \textcolor{black}{Recall that $Q(\cdot)$ is defined by a minimization problem. Thus, in \eqref{DRModel}, we have an inner max-min problem. As such, it is not straightforward to solve the DRO model in \eqref{DRModel}. In this section, we derive an equivalent reformulation of the DRO model in \eqref{DRModel} that is solvable. First, in Proposition~\ref{Prop:InnerSup}, we derive an equivalent dual formulation of the inner maximization problem $\sup \limits_{\Prob\in \calF(\calR, \mu)} \E_\Prob[Q(\ob,\zb,\xi)] $ in \eqref{DRModel} (see  \ref{Appex:Prop1} for a proof). }

\begin{prop}\label{Prop:InnerSup}
For any feasible $(\ob,\zb)$, problem $\sup \limits_{\Prob\in \calF(\calR, \mu)} \E_\Prob[Q(\ob,\zb,\xi)] $ in \eqref{DRModel} is equivalent to
\begin{align}\label{eq:FinalDualInnerMax-1}
\min_{\pmb{\alpha, \phi, \gamma, \lambda, \tau}} &\Bigg \{\sum_{l \in L}\Big[ \sum \limits_{i \in I} \sum \limits_{t \in \calT} \Big(\mu_{t,i,l}^{\mbox{\tiny d}}  \alpha_{t,i,l} +   \mu_{t,i,l}^{\mbox{\tiny M}}\phi_{t,i,l} \Big)+ \sum \limits_{i \in I}\Big( \mu_{i,l}^{\tiny\rho}\gamma_{i,l}+\mu_{i,l}^{\tiny q} \lambda_{i,l} \Big)+ \sum \limits_{(i,j) \in \calA} \mu_{i,j,l}^{\tiny V}\tau_{i,j,l} \Big] \nonumber\\
&  \ \ \ + \max \limits_{\xi \in \calR}  \Big\{ Q(\ob,\zb, \xi) +\sum \limits_{l \in L} \big[ \sum \limits_{i \in I} \sum \limits_{t \in \calT}-(d_{t,i,l}\alpha_{t,i,l}+M_{t,i,l} \phi_{t,i,l})+\sum \limits_{i \in I}-(\rho_{i,l}\gamma_{i,l}+q_{i,l}\lambda_{i,l}) \nonumber \\
& \qquad \qquad +\sum \limits_{(i,j) \in \calA}  -V_{i,j,l} \tau_{i,j,l}\big] \Big \}\Bigg\}.
\end{align}
\end{prop}

\textcolor{black}{Again, the problem in \eqref{eq:FinalDualInnerMax-1} involves an inner max-min problem that is not straightforward to solve in its presented form. However, we next derive an equivalent reformulation of the inner problem in \eqref{eq:FinalDualInnerMax-1} that is solvable.} 
First, we observe that for any feasible $(\ob,\zb, \xi)$, $Q(\ob,\zb, \xi)$ is a feasible linear program  (LP).  Let variable $\beta_{t,i}$ represent the dual associated with constraint \eqref{Const1:Q} for all $i \in I, \ t \in \calT$, variable $\Gamma_{t,i}$ represents the dual of constraint \eqref{Const2:Q}, variable $\psi_{i,j}$ represents the dual of constraint \eqref{Const3:Q} for all $(i,j) \in \calA$, and variable $\varphi_{i,j}$ represents the dual of constraint \eqref{Const32:Q} for all $(i,j) \in \calA$. We formulate $Q(\ob,\zb, \xi)$ in its dual form as:
\begin{subequations}\label{DualofQ}
\begin{align}
& Q(\ob,\zb, \xi ):=    \max_{\textcolor{black}{\pmb{\beta, \Gamma, \psi, \varphi}}} \Bigg\{ \sum \limits_{i \in I}  \Big[\sum \limits_{t\in \calT}  \sum_{l \in L} \big(d_{t,i,l}-\rho_{i,l}z_{t,i} \big) q_{i,l}\beta_{t,i} -(1-\sum_{l \in L}q_{i,l}) z_{t,i} \beta_{t,i} \Big]+ \sum \limits_{i \in I} \sum \limits_{t \in \calT} \sum_{l \in L} q_{i,l}M_{t,i,l}\Gamma_{t,i} \nonumber \\
 & + \sum_{(i,j) \in \calA} \Big[ \sum_{l \in L} q_{i,l} V_{i,j,l} \psi_{i,j}+ (1-\sum_{l \in L}q_{i,l}\big)\hat{V}_{i,j}  \psi_{i,j}  \Big] \nonumber  +\sum_{(i,j) \in \calA} \Big[ \sum_{l' \in L}q_{j,l'} V_{i,j,l'}\varphi_{i,j}+ (1- \sum_{l' \in L}q_{j,l'})\hat{V}_{i,j} \varphi_{i,j}\Big] \Bigg\}\\
&  \qquad  \qquad  \qquad \text{s.t.}  \ -\beta_{t,i}+ \beta_{t,j}+v_t\psi_{i,j} +v_t \varphi_{i,j}   \leq c_{i,j}^t, \qquad \forall (i,j) \in \calA, \ \forall t \in \calT, \label{Const1:DualofQ}\\
&  \qquad  \qquad  \qquad \qquad \beta_{t,i}+\Gamma_{t,i} \leq c_{t}^{\mbox{\tiny p}}, \qquad  \ \ \ \ \forall i \in I, \  t \in \calT, \label{Const2:DualofQ}\\
&  \qquad  \qquad  \qquad \qquad -\ch \leq \beta_{t,i} \leq \cu, \qquad \ \forall i \in I, \  t \in \calT, \label{Const3:DualofQ}\\
&  \qquad  \qquad  \qquad \qquad  \Gamma_{t,i } \leq 0, \ \psi_{i,j} \leq 0,  \ \varphi_{i,j} \leq 0, \ \qquad \ \  \forall (i,j)\in \calA, \ t\in \calT. \label{Const4:DualofQ}
\end{align}  
\end{subequations}
 Note that region \{\eqref{Const1:DualofQ}--\eqref{Const4:DualofQ}\} of $Q(\ob,\zb, \xi)$ is feasible and bounded, support $\calR$ of $\xi$ is finite, and $\pmb{z} \geq 0$. In view of formulation \eqref{DualofQ}, we derive the following reformulation of the inner maximization problem in \eqref{eq:FinalDualInnerMax-1}: 
\begin{subequations}\label{Final_InnerQ}
\begin{align} 
\max \limits_{\substack{\pmb{\beta, \Gamma,\psi, \varphi} \\ \pmb{q, d, \rho, M, V}} } &\Bigg\{ \sum \limits_{i \in I} \Big[\sum \limits_{t\in \calT}  \sum_{l \in L} \big(d_{t,i,l}-\rho_{i,l}z_{t,i} \big) q_{i,l}\beta_{t,i} -(1-\sum_{l \in L}q_{i,l}) z_{t,i} \beta_{t,i} \Big] +\sum \limits_{i \in I} \sum \limits_{t \in \calT} \sum_{l \in L} q_{i,l}M_{t,i,l}\Gamma_{t,i} \nonumber \\
& + \sum_{(i,j)\in \calA} \Big [ \sum_{l \in L} q_{i,l} V_{i,j,l} \psi_{i,j} + (1-\sum_{l \in L}q_{i,l}\big)\hat{V}_{i,j}  \psi_{i,j}+ \sum_{l' \in L} q_{j,l'} V_{i,j,l'} \varphi_{i,j} + (1-\sum_{l' \in L}q_{j,l'}\big)\hat{V}_{i,j}  \varphi_{i,j}  \Big]  \nonumber \\
& +\sum \limits_{l \in L} \Big[ \sum \limits_{i \in I} \sum \limits_{t \in \calT}-(d_{t,i,l}\alpha_{t,i,l}+M_{t,i,l} \phi_{t,i,l})+\sum \limits_{i \in I}-(\rho_{i,l}\gamma_{i,l}+q_{i,l}\lambda_{i,l})  +\sum \limits_{(i,j) \in \calA}  -V_{i,j,l} \tau_{i,j,l}\Big] \Bigg\}\label{Final_Inner_Obj} \\
 \text{s.t.}  & \  \{ \eqref{Const1:DualofQ}-\eqref{Const4:DualofQ}\}, \  \pmb{d} \in [\pmb{\dL, \dU}], \  \pmb{M} \in [\pmb{\ML, \MU}], \  \pmb{\rho} \in [\pmb{\rhoL, \rhoU}], \ \pmb{V }\in [ \pmb{\VL, \VU}], \ \pmb{q} \in \{0, 1\}  \label{Const5:Final_InnerQ}
\end{align}
\end{subequations}
\noindent To limit the number of disasters of type $l\in L$ to $N_l$, we can add constraint $\sum \limits_{i \in I}q_{i,l}=N_l$ to \eqref{Final_InnerQ}. Note that objective function \eqref{Final_Inner_Obj} contains the interaction terms $q_{i,l}d_{t,i,l}\beta_{t,i}$,   $q_{i,l}\rho_{i,l}\beta_{t,i}$, $q_{i,l} M_{t,i,l}\Gamma_{t,i}$, $q_{i,l}V_{i,j,l}\psi_{i,j}$, and $q_{j,l'} V_{i,j,l'}\varphi_{i,j}$. In Proposition~\ref{Prop:MILP_Inner}, we derive an equivalent mixed-integer linear programming (MILP) reformulation of problem \eqref{Final_InnerQ} (see \ref{Appex:McCIneq} for a proof).
\color{black}
\begin{prop}\label{Prop:MILP_Inner}
For fixed $(\pmb{z, \alpha, \phi, \gamma, \lambda, \tau})$, problem \eqref{Final_InnerQ} is equivalent to the following MILP.\color{black}
\allowdisplaybreaks
\begin{subequations}\label{Final_InnerQ_MILP}
\begin{align} 
& H (\pmb{z,  \alpha, \phi, \gamma, \lambda, \tau})= \nonumber \\
&\max \limits_{\substack{\pmb{\beta, \Gamma,\psi, \varphi, a,\Theta, \kappa, \varrho} \\ \pmb{a', \Theta', \kappa', \varrho', h, g, F, \pi} \\ \pmb{\eta, \Phi, \varpi, b,\Lambda}}}   \ \ \Bigg\{ \sum \limits_{i \in I} \sum \limits_{t\in \calT} \sum_{l \in L}  \Big[\dL_{t,i,l} \big(k_{t,i,l}-\alpha_{t,i,l} \big) + \Delta d_{t,i,l}\big(h_{t,i,l}-\alpha_{t,i,l}a_{t,i,l} \big) \Big] \nonumber\\
&\qquad \qquad \ \ \  +\sum_{i \in I} \sum_{t \in T}  z_{t,i}\big (\sum_{l \in L} k_{t,i,l}- \beta_{t,i}\big) \nonumber \\
&\qquad \qquad \ \ \ - \sum \limits_{i \in I}   \sum_{l \in L} \Big[  \rhoL_{i,l} \big( \sum \limits_{t\in \calT}z_{t,i} k_{t,i,l}+\gamma_{i,l} \big)+ \Delta \rho_{i,l} \big( \sum \limits_{t\in \calT}z_{t,i} g_{t,i,l}+\Theta_{i,l}\gamma_{i,l} \big) \Big] \nonumber \\
&\qquad \qquad \ \ \ + \sum \limits_{i \in I} \sum \limits_{t \in \calT} \sum_{l \in L} \Big[ \ML_{t,i,l}\big(F_{t,i,l}-\phi_{t,i,l} \big)+\Delta M_{t,i,l} \big(\pi_{t,i,l}-\phi_{t,i,l}\kappa_{t,i,l} \big) \Big]  \nonumber \\
&\qquad \qquad \ \ \  +  \sum_{(i,j)\in \calA}  \sum_{l \in L} \Big[ \VL_{i,j,l} \big (\eta_{i,j,l}+\varpi_{i,j,l} -\tau_{i,j,l} \big) + \Delta V_{i,j,l} \big (\Phi_{i,j,l}+\Lambda_{i,j,l} -\varrho_{i,j,l} \tau_{i,j,l} \big)\nonumber \\
&\qquad \qquad \ \ \ + \sum_{(i,j)\in \calA}  \Big[ \big(\psi_{i,j}-\sum_{l \in L}\eta_{i,j,l}\big)\hat{V}_{i,j} + \big(\varphi_{i,j}-\sum_{l \in L}\varpi_{i,j,l}\big)\hat{V}_{i,j}  \Big] -\sum_{i \in I} \sum_{l \in L} q_{i,l} \lambda_{i,l} \Bigg\}  \label{Obj_Inter2} \\
 \text{s.t.}  & \ \  \{ \eqref{Const1:DualofQ}-\eqref{Const4:DualofQ}\}, \pmb{q} \in \{0, 1\},  (\pmb{a,\Theta, \kappa, \varrho}) \in \{0, 1\}, \eqref{Mac1}-\eqref{Mac14},
\end{align}
\end{subequations}
where $\Delta d_{t,i,l}=(\dU_{t,i,l}-\dL_{t,i,l})$, $\Delta \rho_{i,l}=(\rhoU_{i,l}-\rhoL_{i,l})$, $\Delta M_{t,i,l}=(\MU_{t,i,l}-\ML_{t,i,l})$, $\Delta V_{i,j,l}=(\VU_{i,j,l}-\VL_{i,j,l})$.
\end{prop}
\color{black}
Combing the inner problem  $\max \limits_{\xi \in \calR}  \{ \cdot \}$ in the form of \eqref{Final_InnerQ_MILP}  with the outer minimization problems in \eqref{eq:FinalDualInnerMax-1} and \eqref{DRModel}, we derive the following equivalent reformulation of the DRO model in \eqref{DRModel},
 \allowdisplaybreaks
\begin{subequations}\label{Final_DR_MILP}
\begin{align}
\min_{\substack{ \pmb{o,z,\alpha, \phi} \\ \pmb{\gamma, \lambda, \tau, \delta} } }& \ \Bigg \{ \sum_{i \in I} f_i o_i + \sum \limits_{i \in I} \sum\limits_{t \in \calT} \ca z_{t,i} +  \sum \limits_{i \in I} \sum \limits_{t \in \calT} \sum_{l \in L}\Big(\mu_{t,i,l}^{\mbox{\tiny d}}  \alpha_{t,i,l} +   \mu_{t,i,l}^{\mbox{\tiny M}}\phi_{t,i,l} \Big) \nonumber \\
&  \ + \ \  \sum \limits_{i \in I}\sum_{l \in L} \Big( \mu_{i,l}^{\tiny\rho}\gamma_{i,l}+\mu_{i,l}^{\tiny q} \lambda_{i,l} \Big) +  \sum \limits_{(i,j) \in \calA} \sum_{l \in L}\mu_{i,j,l}^{\tiny V}\tau_{i,j,l} + \delta \Bigg\} \label{FinalDR} \\
\text{s.t. } & \eqref{First-stage:C2}-\eqref{First-stage:C3}, \\
& \delta \geq H (\pmb{z,  \alpha, \phi, \gamma, \lambda, \tau}).\label{Const1:FinalDR}
\end{align}
\end{subequations}
\begin{prop}\label{Prop:Convex}
For any fixed values of variables $\pmb{z, \ \alpha, \ \phi, \ \gamma, \ \lambda}$, and $\pmb{\tau}$, $H (\pmb{z,  \alpha, \phi, \gamma, \lambda, \tau}) <\infty$. Furthermore, function $(\pmb{z,  \alpha, \phi, \gamma, \lambda, \tau}) \mapsto H (\pmb{z, \alpha, \phi, \gamma, \lambda, \tau}) $ is a convex piecewise linear function in $\pmb{z,  \ \alpha, \ \phi, \ \gamma, \ \lambda}$, and $\pmb{\tau}$  with a finite number of pieces (see \ref{Appix:Prop_Covex} for a detailed proof).
\end{prop}

\section{\textcolor{black}{Trade-off Model}\label{sec:Trade}}
\noindent As pointed out by \cite{chen2019robust}, \cite{hurwicz1951generalized} is arguably the first to present a decision criterion that model the trade-off between pessimistic and optimistic objectives. In this paper, we are interested studying the trade-off between considering distributional ambiguity and following distributional belief. We formulate a model that minimizes this trade-off as follows.
\begin{align}\label{Trade-off_Model2}
\min_{(o,z)} \Big\{ \text{first-stage objective}+ (1-\theta) \sup_{\mathbb{P}' \in \calF} \E_{\mathbb{P}'}[ Q(o, z,\xi)]+ \theta  \ \E_\mathbb{P}[ Q(\ob,\zb, \xi)] \Big\}.
\end{align}

Problem \eqref{Trade-off_Model2} finds first-stage planning decisions that minimize the first-stage cost and the trade-off between considering distributional ambiguity and following a distributional belief for $\xi$.  Parameter $\theta \in [0,1]$ represents the level of optimism. If $\theta=0$, then  problem \eqref{Trade-off_Model2} recovers the DRO criterion that solves the DRO model. On the other hand, if $\theta=1$, then problem \eqref{Trade-off_Model2} recovers the optimistic criterion which solves the SP model under prefect distributional belief. $0 < \theta <1$ represent a trade-off between the optimistic and pessimistic perception of the objective value. We apply the same techniques in Section~\ref{sec3:Reformulation} to derive an equivalent reformulation of \eqref{Trade-off_Model2}.

\section{Solution Approaches}\label{sec:SolutionApproches}

\noindent In this section, we propose a decomposition algorithm to solve the DRO model and a Monte Carlo Optimization procedure to obtain near-optimal solutions to the SP model.

\subsection{DRO--Decomposition Algorithm}\label{sec:DRO_decomp}

\noindent Proposition~\ref{Prop:Convex} suggests that constraint \eqref{Const1:FinalDR} describes the epigraph of a convex and piecewise linear function of decision variables in formulation \eqref{Final_DR_MILP}. Therefore, given the two stage characteristic of our problem, it is natural to attempt to solve \eqref{Final_DR_MILP} (equivalently, the DRO model in \eqref{DRModel}) via a decomposition algorithm. Algorithm~\ref{Alg1:CAG} presents our decomposition algorithm. Algorithm~\ref{Alg1:CAG}  is finite because we identify a new piece of the function \textcolor{black}{$H (\pmb{z,  \alpha, \phi, \gamma, \lambda, \tau})$}  each time the set $\{\mathcal{L} (\pmb{z}, \pmb{\alpha, \phi, \gamma, \lambda, \tau}, \delta) \} \geq 0$  is augmented in step 4, and the function has a finite number of pieces according to Proposition~\ref{Prop:Convex}. \textcolor{black}{Thus, the algorithm will terminate after a finite number of iterations}. 

\begin{algorithm}[t!]
\small
\caption{DRO--Decomposition algorithm.}
\label{Alg1:CAG}
  \renewcommand{\arraystretch}{0.3}
\noindent \textbf{1. Input.} Feasible region $\{\eqref{First-stage:C2}-\eqref{First-stage:C3}  \}$; set of cuts $ \lbrace \mathcal{L} (\pmb{z}, \pmb{\alpha, \phi, \gamma, \lambda, \tau}, \delta) \geq 0 \rbrace=\emptyset $; $LB=-\infty$ and $UB=\infty.$

\noindent \textbf{2. Master Problem.} Solve the following master problem
\begin{subequations}\label{Master}
\begin{align}
Z=\min_{\substack{ \pmb{o,z,\alpha, \phi} \\ \pmb{\gamma, \lambda, \tau, \delta} } }& \ \Bigg \{ \sum_{i \in I} f_i o_i + \sum \limits_{i \in I} \sum\limits_{t \in \calT} \ca z_{t,i} +  \sum \limits_{i \in I} \sum \limits_{t \in \calT} \sum_{l \in L}\Big(\mu_{t,i,l}^{\mbox{\tiny d}}  \alpha_{t,i,l} +   \mu_{t,i,l}^{\mbox{\tiny M}}\phi_{t,i,l} \Big) \nonumber \\
&  \ + \ \  \sum \limits_{i \in I}\sum_{l \in L} \Big( \mu_{i,l}^{\tiny\rho}\gamma_{i,l}+\mu_{i,l}^{\tiny q} \lambda_{i,l} \Big) +  \sum \limits_{(i,j) \in \calA} \sum_{l \in L}\mu_{i,j,l}^{\tiny V}\tau_{i,j,l} + \delta \Bigg\} \\
\text{s.t.} & \ \eqref{First-stage:C2}-\eqref{First-stage:C3}\\
&  \ \mathcal{L} (\pmb{z}, \pmb{\alpha, \phi, \gamma, \lambda, \tau}, \delta) \geq 0 
 \end{align}
\end{subequations}
and record an optimal solution $(\pmb{o^*, z^*, \alpha^*, \phi^*, \gamma^*, \lambda^*, \tau^*}, \delta^*)$ and set $LB=Z^*$.

\noindent \textbf{3. Sub-problem.} 
\begin{enumerate}
\item[3.1] with $(\pmb{z, \alpha, \phi, \gamma, \lambda, \tau})$  fixed to $(\pmb{ z^*, \alpha^*, \phi^*, \gamma^*, \lambda^*, \tau^*})$, solve problem $H (\pmb{z,  \alpha, \phi, \gamma, \lambda, \tau})$ in \eqref{Final_InnerQ_MILP}. 
\item[3.2] record optimal solution  $\pmb{k^*, a^*, h^*, \Theta^*,  g^*, F^*,  \kappa^*, \pi^*, \eta^*, \varrho^*, \varpi^*, \Phi^*, \Lambda^*}$ and value $H^*$.

 Then, set $UB=\min \{ UB, \ H^*+ (LB-\delta^*) \} $.
\end{enumerate}

\noindent \textbf{4. if} $\delta^* \geq  H^*$ or $UB-LB \leq \epsilon$ \textbf{then} stop and return $\pmb{o}^*$ and $\pmb{z}^*$ as the optimal solution to~\eqref{Master}.

\noindent $\ \ \ \  $\textbf{else}  add the following cut to the set of cuts $ \lbrace \mathcal{L} (\pmb{z}, \pmb{\alpha, \phi, \gamma, \lambda, \tau}, \delta)\geq 0 \rbrace$  and go to step 2.
\vspace{2mm}
\begin{align*}
\delta &\geq  \Bigg\{\sum \limits_{i \in I} \sum \limits_{t\in \calT} \sum_{l \in L}  \Big[\dL_{t,i,l} \big(k_{t,i,l}^*-\alpha_{t,i,l} \big) + \Delta d_{t,i,l}\big(h_{t,i,l}^*-\alpha_{t,i,l}a_{t,i,l}^* \big) \Big] \nonumber\\
& \ \ \ \ \ +\sum_{i \in I} \sum_{t \in T}  z_{t,i}\big (\sum_{l \in L} k_{t,i,l}^*- \beta_{t,i}^*\big)- \sum \limits_{i \in I}   \sum_{l \in L} \Big[  \rhoL_{i,l} \big( \sum \limits_{t\in \calT}z_{t,i} k_{t,i,l}^*+\gamma_{i,l} \big)+ \Delta \rho_{i,l} \big( \sum \limits_{t\in \calT}z_{t,i} g_{t,i,l}^*+\Theta_{i,l}^*\gamma_{i,l} \big) \Big] \nonumber \\
& \ \ \ \ \ + \sum \limits_{i \in I} \sum \limits_{t \in \calT} \sum_{l \in L} \Big[ \ML_{t,i,l}\big(F_{t,i,l}^*-\phi_{t,i,l} \big)+\Delta M_{t,i,l} \big(\pi_{t,i,l}^*-\phi_{t,i,l}\kappa_{t,i,l}^* \big) \Big]  \nonumber \\
& \ \ \ \ \  +  \sum_{(i,j)\in \calA}  \sum_{l \in L} \Big[ \VL_{i,j,l} \big (\eta_{i,j,l}^*+\varpi_{i,j,l}^* -\tau_{i,j,l} \big) + \Delta V_{i,j,l} \big (\Phi_{i,j,l}^*+\Lambda_{i,j,l}^* -\varrho_{i,j,l}^* \tau_{i,j,l} \big)\nonumber \\
& \ \ \ \ \  + \sum_{(i,j)\in \calA}  \Big[ \big(\psi_{i,j}^*-\sum_{l \in L}\eta_{i,j,l}^*\big)\hat{V}_{i,j} + \big(\varphi_{i,j}^*-\sum_{l \in L}\varpi_{i,j,l}^*\big)\hat{V}_{i,j}  \Big] -\sum_{i \in I} \sum_{l \in L} q_{i,l}^* \lambda_{i,l} \Bigg\} 
\end{align*}
\noindent $\ \ $ \textbf{end if}
\end{algorithm}

 \subsection{Monte Carlo Optimization (MCO)}\label{sec:MCO}

 \noindent Note that it is difficult to obtain an exact optimal solution to the two-stage SP in \eqref{TSM}. Indeed, evaluating the values of $\E[Q (o,z,  \xi)]$ involves taking multi-dimensional integrals \citep{birge2011introduction}. Thus, we resort to Monte Carlo  approximation approach to obtain near-optimal solutions to \eqref{TSM} in a reasonable time. In this approach,  we replace the distribution of $\xi$ with a (discrete) distribution based on $N$ samples of $\xi$, and then we solve the following sample average approximation (SAA) formulation of the SP in \eqref{TSM}.
  \begin{subequations}\label{SAA}
   \begin{align}
 & \ \  \upsilon_N= \min \Big \{  \sum_{i \in I} f_i \hat{o}_i + \sum \limits_{i \in I} \sum\limits_{t \in \calT} \ca \hat{z}_{t,i} + \hat{F}_{N} \Big \} \label{SAA_Obj}\\
& \qquad \eqref{First-stage:C2}-\eqref{First-stage:C3}, \  \qquad \eqref{Const1:Q}-\eqref{Const4:Q}, \ \text{for all } n=1, \ldots,N, \label{SAA_Const2}
  \end{align}
   \end{subequations}

where $ \hat{F}_{N}:= \sum_{n=1}^{N} \frac{1}{N} \Big[ \sum \limits_{i \in I} \sum \limits_{t \in \calT} (\cp y^n_{t,i}+\cu u^n_{t,i} +\ch e^n_{t,i} ) + \sum \limits_{t \in \calT} \sum \limits_{(i,j) \in \calA} c_{i,j}^{\mbox{\tiny t}} x^{n,t}_{i,j}\Big]$.  Note that in the SAA formulation \eqref{SAA}, we associate all scenario-dependent parameters, variables, and constraints with a scenario index $n$ for all $n=1,\ldots,N$.  For example, we replace parameters $d_{t,i,l}$ by $d_{t,i,l}^n$ to represent the demand for relief item $t$ if a disaster of type $l$ occurs at location $i$ in scenario $n$.  In addition, constraints \eqref{Const1:Q}-\eqref{Const4:Q} are incorporated in each scenario.   \textcolor{black}{The sample average $ \hat{F}_{N}$ is an unbiased estimator of the expected value $F:=\E[Q(\ob,\zb,\xi)]$ in \eqref{TSM} \citep{shapiro2003monte}.  By the Law of Large Numbers and \cite{shapiro2003monte}, we have $\hat{F}_N \rightarrow F$ with probability one (w.p.1) as $N\rightarrow \infty$  \citep{homem2014monte, kleywegt2002sample}. It follows that $\upsilon_N \rightarrow \upsilon$ w.p.1 as $N \rightarrow \infty$, i.e., the optimal value of the SAA formulation converges to that of the SP as $N \rightarrow \infty$. However, for a fixed $N$, formulation \eqref{SAA} reduces to an MILP. Hence, one would expect solution time of solving the SAA formulations to increase as $N$ increases.} Algorithm \ref{alg:algorithm1} in  \ref{Appex:MCOProcedure} summarizes the MCO algorithm that determines an appropriate sample size $N$ and obtain near-optimal solutions to the SP model based on its SAA  within a reasonable time and high accuracy.

\section{Computational Experiments}\label{sec:Exp}

\noindent In this section, \textcolor{black}{we conduct experiments using the proposed DRO, SP, and trade-off approaches for two very different disaster relief contexts--a hurricane season (Section~\ref{sec5:Hurrican}) and an earthquake (Section~\ref{sec:Earthquack}). These case studies differ in several characteristics, including the length of forewarning, the size of at-risk and affected areas, and the time from onset to recovery. For each case study, we analyze and compare the DRO, SP, and trade-off models' optimal decisions and their in-sample and out-of-sample simulation performance. To further analyze computational performance of the DRO approach,} we apply the DRO decomposition algorithm on larger, randomly generated networks (Section~\ref{sec:performance}).  We implemented the models and algorithms using the AMPL modeling language and CPLEX (version 12.6.2) as the solver with its default settings. We ran the experiments on a laptop with an Intel Core i7 processor, 2.6 GHz CPU, and 16 GB (2667 MHz DDR4) of memory.

\subsection{Hurricane case study}\label{sec5:Hurrican}

\noindent In this section, we consider the planning process for the Atlantic hurricane season in the US. The season runs from June through November every year and affects the North Atlantic Ocean, the Caribbean Sea, and the Gulf of Mexico. \textcolor{black}{Multiple disasters may occur over the course of this time period.} One month before the season, experts release predictions on how many hurricanes and major storms to expect (Atlantic Oceanographic and Meteorological Laboratory). However, they do not release specific landfall locations; these are only available immediately before a storm hits land. \textcolor{black}{This information availability mirrors the two-stage modeling approach, where pre-disaster season decisions are made based on estimates of number of occurrences, and post-disaster decisions are made when affected areas, demand, and other uncertainty is known. The key items planners seek to preposition are water, food, and medical kits. Pre-disaster, they decide where and how many of these items to keep in selected warehouses, and post-disaster, they are distributed, with possible additional procurement.} 

\begin{figure}[t!]
\center
    \includegraphics[scale=0.7]{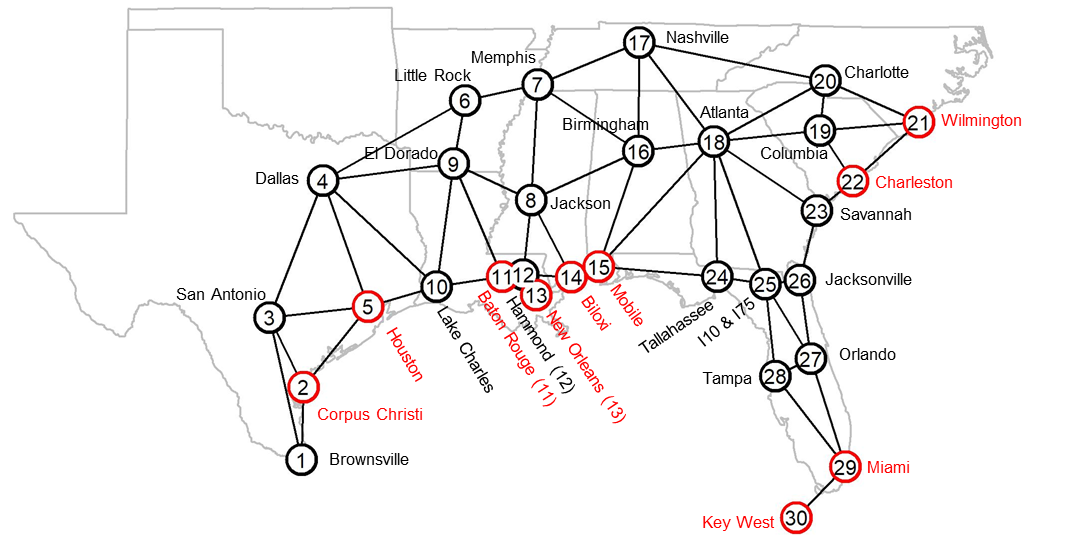}
    \caption{\black{Map of 30 nodes and 112 arcs in Southeast US. Disaster-prone nodes are marked in red.}}  

    \label{Fig_map}
\end{figure}

Our data is based on a case study presented in \cite{rawls2010pre} and \cite{velasquez2020prepositioning}. We present the network of 30 nodes and 112 transportation arcs in Figure~\ref{Fig_map}. Each node represents a potential warehouse or demand location. \black{Ten of the nodes, $i \in  \{2, 5, 11, 13, 14, 15, 21, 22, 29,$ $30\}$, are at risk of being hit by a disaster (designated as red in Figure~\ref{Fig_map}), and twenty nodes are not (black)}. For illustrative purposes, we consider two types of disasters minor ($l=1$) and major ($l=2$), and three relief supplies ($|\calT|=3$): water, medical kits, and food.  Water is measured in 1,000 gallon units, medical kits are in single units, and food is stored in the form of ready-to-eat meals (MREs), and is measured in units of 1000 MREs. We summarize the costs of procurement and transportation as well as storage volumes in Table~\ref{Table:Items_Costs} in  \ref{Appx:Cas1_data}. We set the fixed cost and capacity for potential warehouses to \$188,400 and 408,200ft$^3$, respectively, as in \cite{rawls2010pre} and \cite{velasquez2020prepositioning}.

We make the following assumptions: (1) demand is zero at locations that are not potential landfall nodes and is uncertain at coastal locations susceptible to hurricanes, (2) only one hurricane can impact each node (it is very unlikely that two hurricanes hit the same location in the same hurricane season \citep{velasquez2020prepositioning}), (3) post-disaster procurement is twice as expensive as pre-disaster acquisition (i.e., $c^{\mbox{\tiny p}}=2 c^{\mbox{\tiny a}}$, as in \cite{velasquez2020prepositioning}), (4) the shortage cost is four times the acquisition cost (i.e., $c^{\mbox{\tiny u}}=4 c^{\mbox{\tiny a}}$), and (5) the holding cost is equal to the acquisition cost (i.e., $c^{\mbox{\tiny h}}= c^{\mbox{\tiny a}}$). These are consistent with prior literature \citep{sabbaghtorkan2020prepositioning}.

Our uncertainty is defined as follows. We present the mean demand after a disaster at each landfall node $\mu^{\mbox{\tiny d}}_{t,i,l}$ in Table~\ref{Table:Demand} in  \ref{Appx:Cas1_data}. For notation convenience, we let $t=1$, $t=2$, and $t=3$ respectively represent water, food, and medical kits. We set the standard deviation $\sigma_t^{\mbox{\tiny d}}$ of the demand for each relief item $t$ to $0.5 \mu^{\mbox{\tiny d}}_{t,i.l}$. For the fraction of usable prepositioned items available after the disaster, we let $\mu^{\tiny \rho}_{i,1}=0.5$ and $[\rhoL_{i,1},\rhoU_{i,1}]=[0.4, 0.6]$ (i.e., 10\% below and above $\mu^{\tiny \rho}=50\%$), and $\mu^{\tiny \rho}_{i,2}=0.05$ and $[\rhoL_{i,2}, \rhoU_{i,2}]=[0, 0.10]$. The maximum procurement quantity of each relief item $t$ is  $M_{t,i,l}=2\mu^{\mbox{\tiny d}}_{t,i,1}$.  We let $\mu^{\mbox{\tiny M}}_{t,i,l}= 2\mu^{\mbox{\tiny d}}_{t,i,1}$ and $[\ML_{t,i,l}, \MU_{t,i,l}]=[0.9\mu^{\mbox{\tiny M}}_{t,i,l}, 1.10 \mu^{\mbox{\tiny M}}_{t,i,l}]$. In the text, we designate the number of minor and major disasters to be nMinor and nMajor, respectively. For illustrative purposes, we consider two combinations, where (nMinor, nMajor) is (2, 1) and (4, 2), respectively. To limit the number of minor and major disasters to nMinor and nMajor, respectively, we add constraints $\sum \limits_{i \in I}q_{i,1}=$nMinor and $\sum \limits_{i \in I}q_{i,2}=$ nMajor to the sub-problem $H (\pmb{z, \alpha, \phi, \gamma, \lambda, \tau})$ in \eqref{Final_InnerQ_MILP}. 

To approximate the lower and upper bound of ($\pmb{d}, \pmb{V}$), we follow the same  procedure as  in prior applied DRO studies (see, e.g., \cite{jian2017integer,shehadehDMFRS, wang2020distributionally}) as follows. We first generate $N=1,000$ in-sample data of $\pmb{d}$ and $\pmb{V}$ by following a lognormal (LogN) and a truncated normal distribution (as in \cite{ni2018location}), respectively, with the generated mean values and standard deviations of these random parameters. Second, we respectively use the 20\% and 80\% of the $N$ in-sample data  as the lower and upper bounds.


  \subsubsection{Analysis of optimal prepositioning decisions}\label{sec5:OptimalSols}

\begin{figure}[t!]
     \begin{subfigure}[b]{0.33\textwidth}
 \centering
        \includegraphics[width=\textwidth]{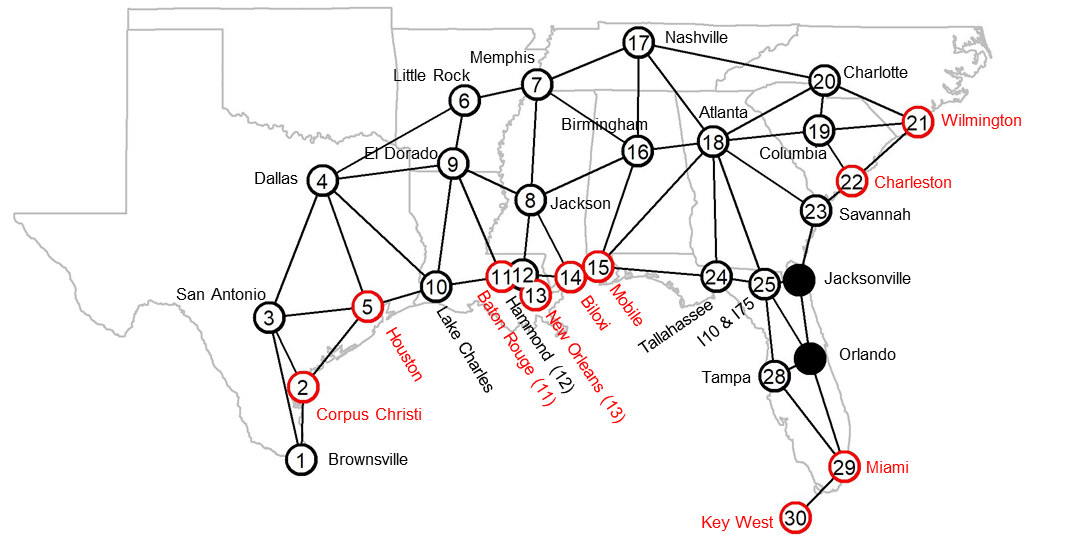}
        \caption{SP, 2 Facilities}
        \label{FigVIa}
    \end{subfigure}%
       \begin{subfigure}[b]{0.33\textwidth}
 \centering
        \includegraphics[width=\textwidth]{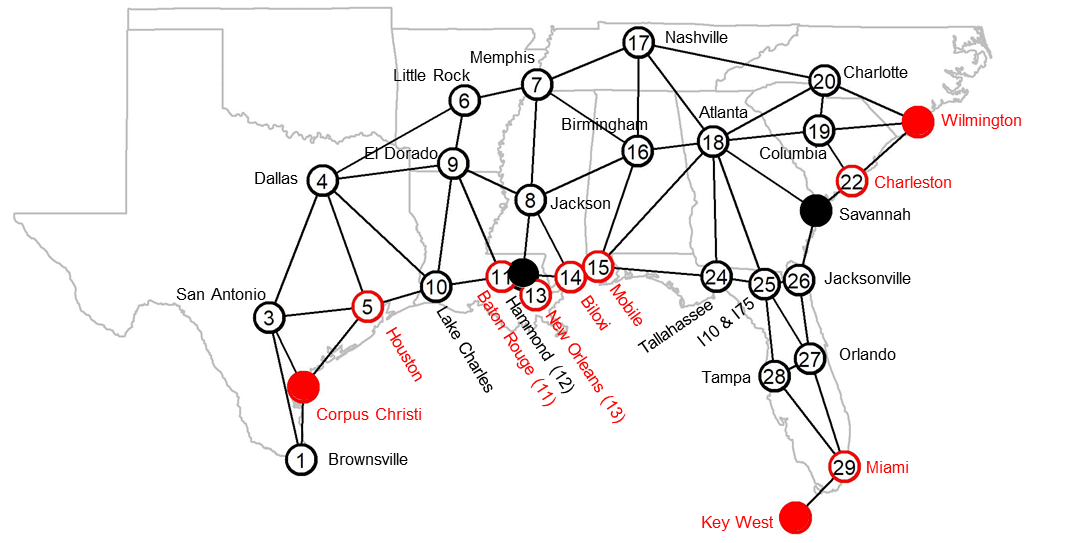}
        \caption{Trade (0.7), 5 Facilities}
        \label{FigVIa}
    \end{subfigure}%
       \begin{subfigure}[b]{0.33\textwidth}
 \centering
        \includegraphics[width=\textwidth]{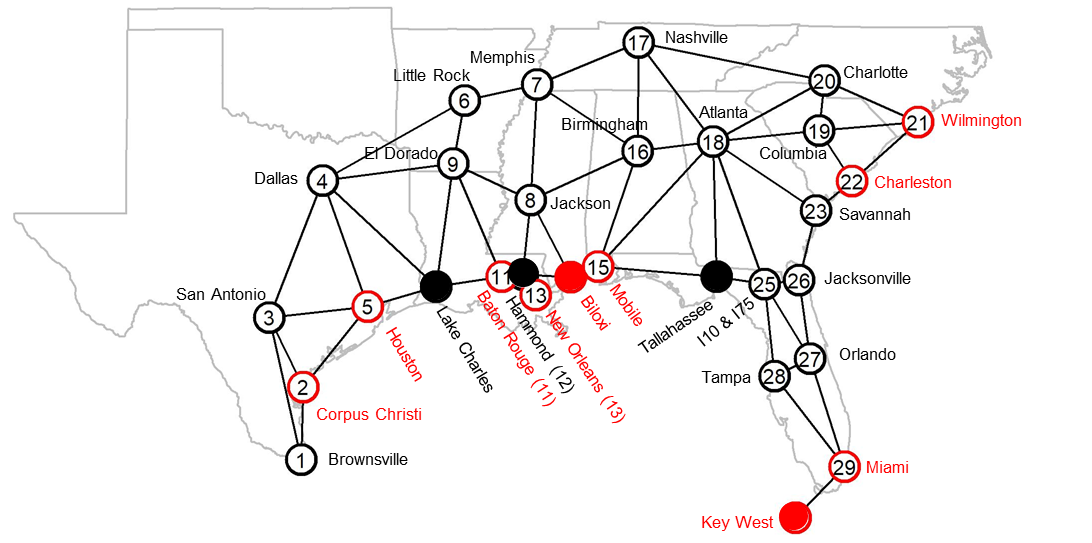}
        \caption{Trade (0.5), 5 Facilities}
        \label{FigVIa}
    \end{subfigure}%
    
            \begin{center}
        \begin{subfigure}[b]{0.33\textwidth}
            \includegraphics[width=\textwidth]{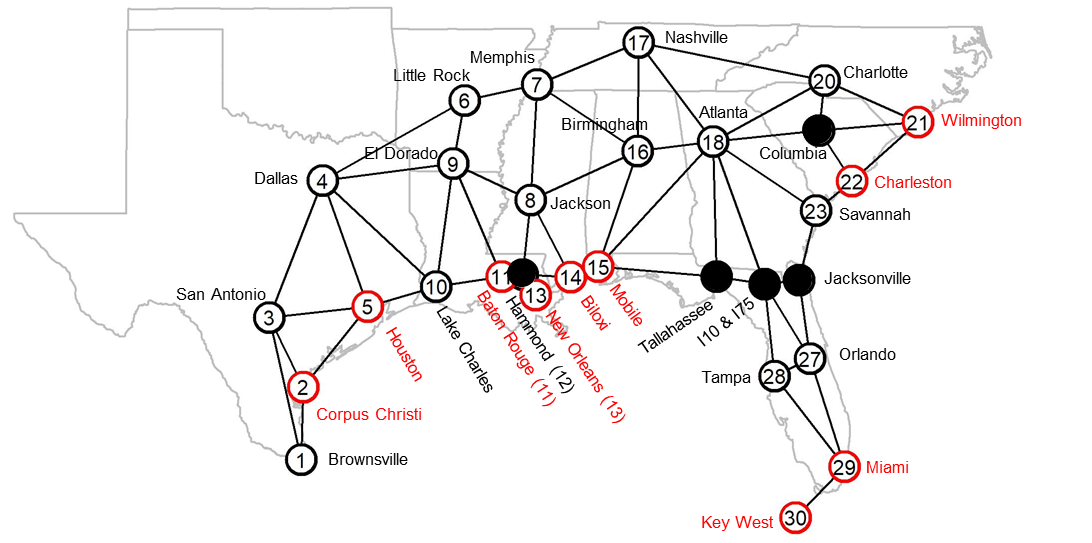}
      \caption{Trade (0.3), 5 Facilities}
      \label{FigVIb}
    \end{subfigure}%
    \begin{subfigure}[b]{0.33\textwidth}
            \includegraphics[width=\textwidth]{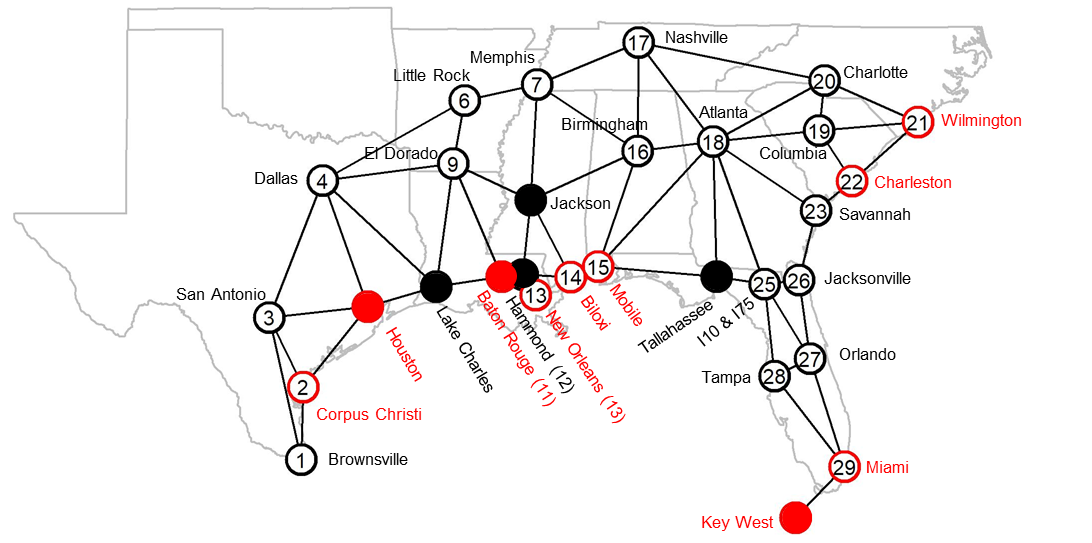}
      \caption{DRO, 7 Facilities}
      \label{FigVIb}
    \end{subfigure}%
        \end{center}
        \caption{Comparison of the optimal medium facility locations for (nMinor,nMajor)=(2,1). Dark and red-filled circles are optimal locations at safe and disaster-prone nodes, respectively.}\label{Fig_Loc_21_medium}
\end{figure}

\begin{figure}[t!]
     \begin{subfigure}[b]{0.33\textwidth}
 \centering
        \includegraphics[width=\textwidth]{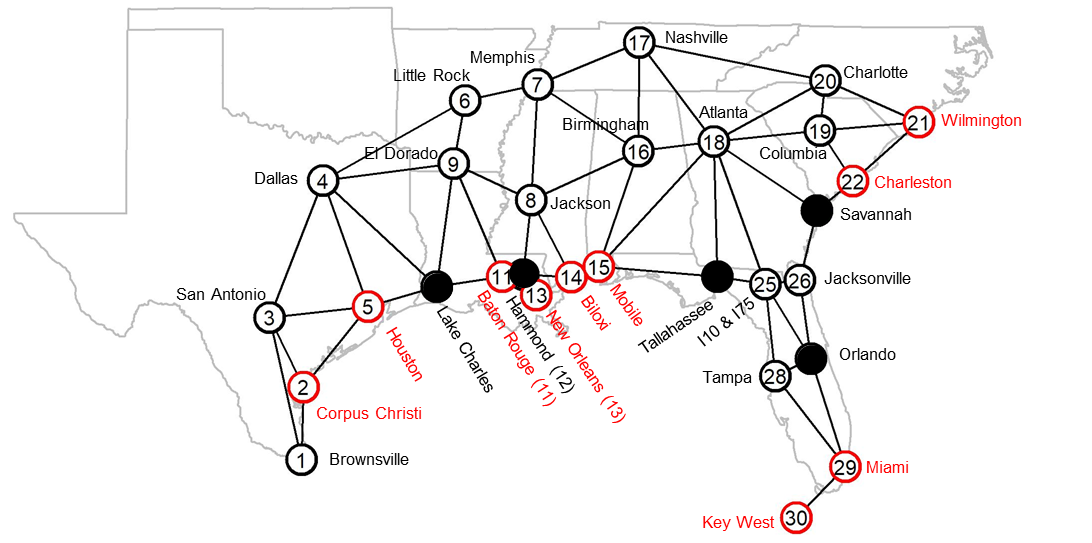}
        \caption{SP, 5 Facilities}
        \label{FigVIaMedium}
    \end{subfigure}%
       \begin{subfigure}[b]{0.33\textwidth}
 \centering
        \includegraphics[width=\textwidth]{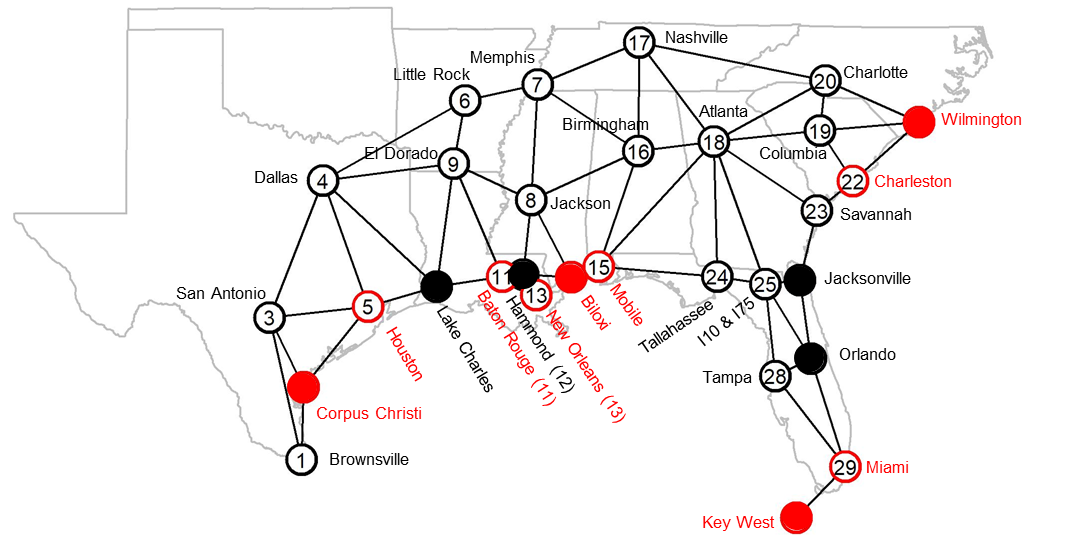}
        \caption{Trade (0.7), 8 Facilities}
        \label{FigVIa}
    \end{subfigure}%
       \begin{subfigure}[b]{0.33\textwidth}
 \centering
        \includegraphics[width=\textwidth]{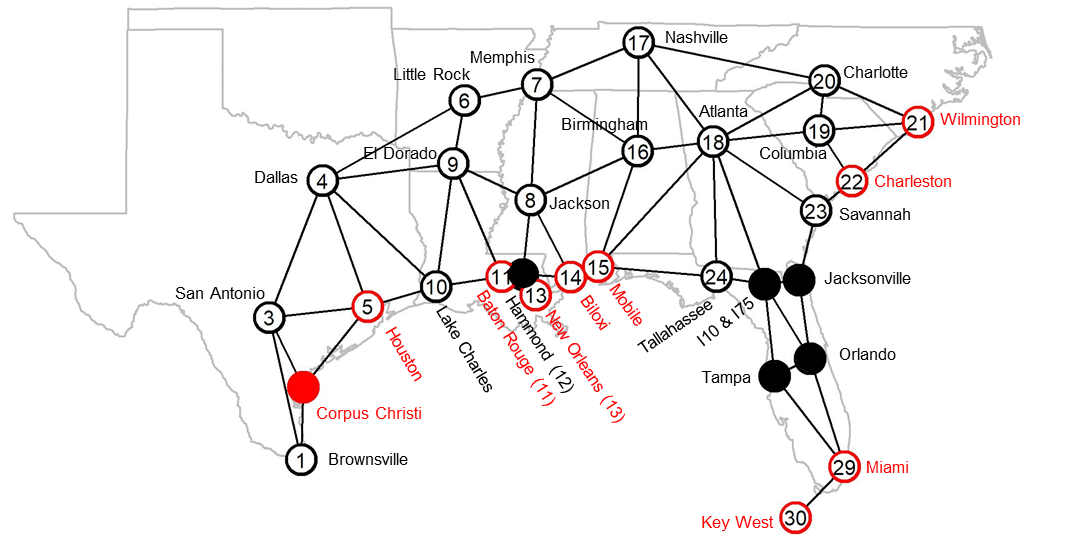}
        \caption{Trade (0.5), 6 Facilities}
        \label{FigVIa}
    \end{subfigure}%
    
    \begin{center}
            \begin{subfigure}[b]{0.33\textwidth}
            \includegraphics[width=\textwidth]{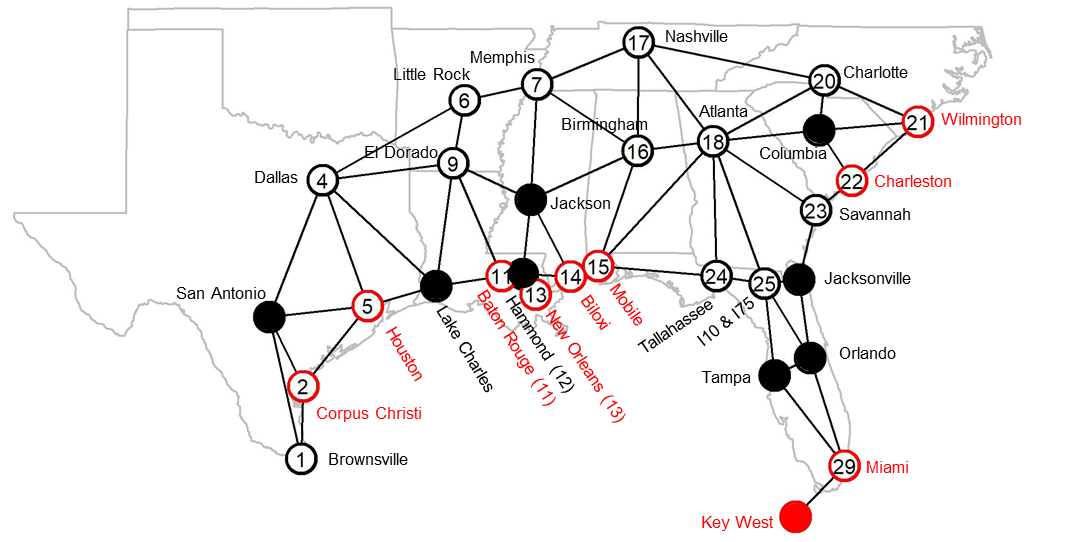}
      \caption{Trade(0.3), 9 Facilities}
      \label{FigVIb}
    \end{subfigure}%
    \begin{subfigure}[b]{0.33\textwidth}
            \includegraphics[width=\textwidth]{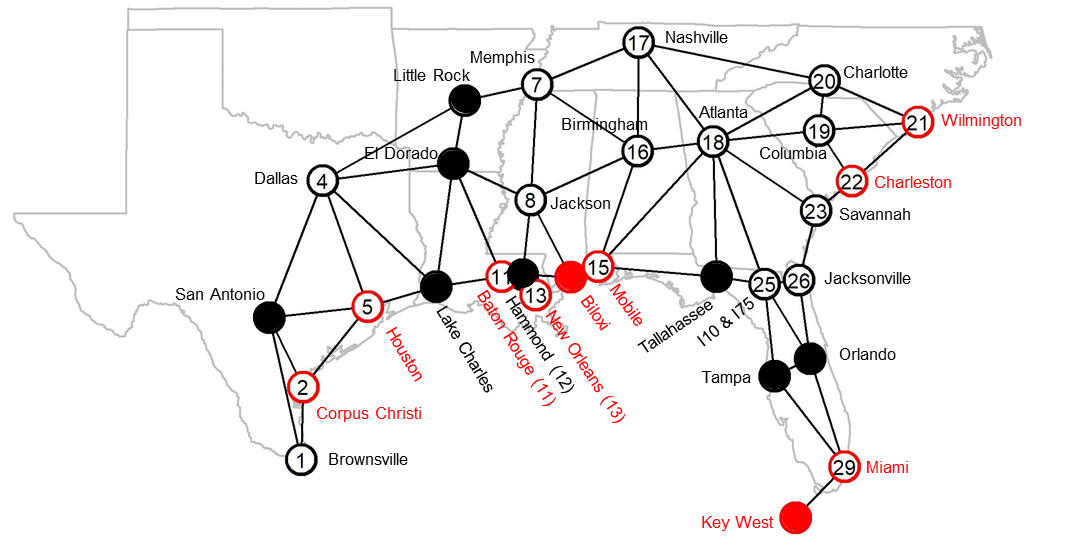}
      \caption{DRO, 10 Facilities}
      \label{FigVIb}
    \end{subfigure}%
        \end{center}
        \caption{Comparison of the optimal medium facility locations for (nMinor,nMajor)=(4,2).}\label{Fig_Loc_42_medium}
\end{figure}

\begin{figure}[t!]
     \begin{subfigure}[b]{0.33\textwidth}
 \centering
        \includegraphics[width=\textwidth]{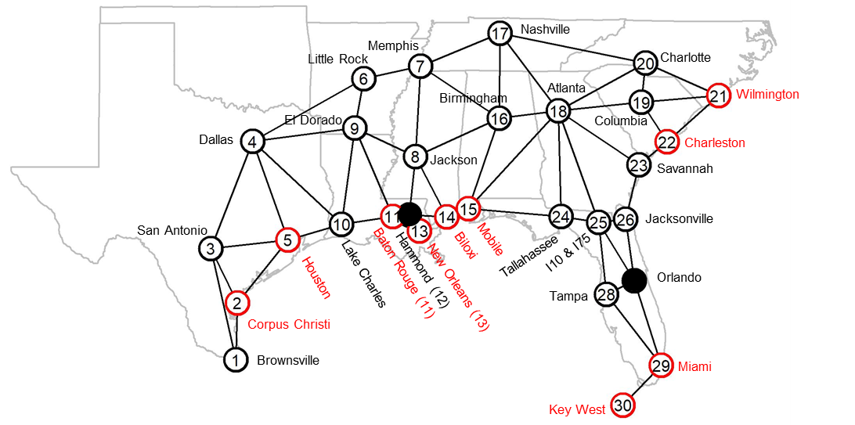}
        \caption{SP, 2 Facilities}
        \label{FigVIa}
    \end{subfigure}%
       \begin{subfigure}[b]{0.33\textwidth}
 \centering
        \includegraphics[width=\textwidth]{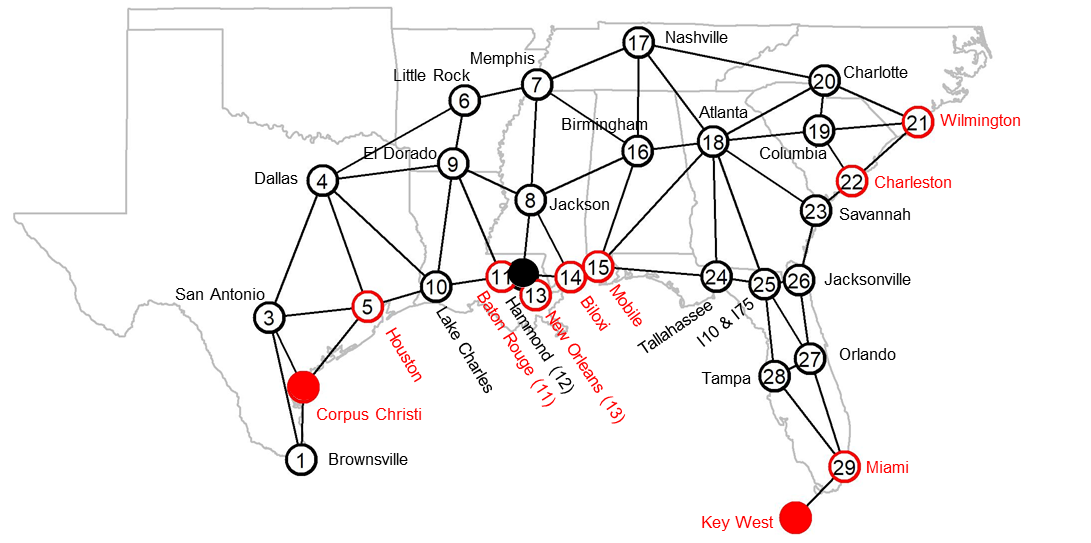}
        \caption{Trade (0.7), 3 Facilities}
        \label{FigVIa}
    \end{subfigure}%
       \begin{subfigure}[b]{0.33\textwidth}
 \centering
        \includegraphics[width=\textwidth]{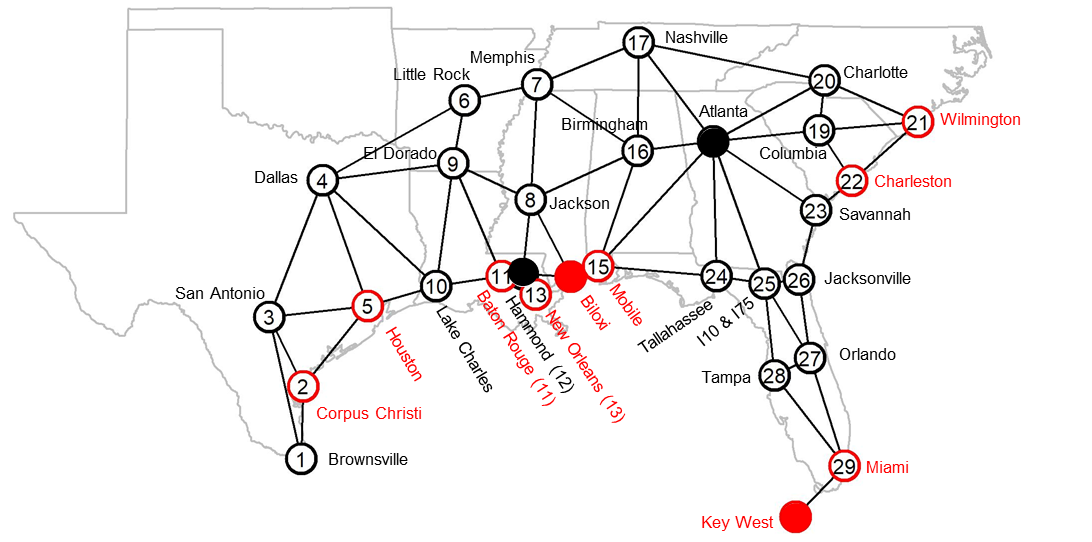}
        \caption{Trade (0.5), 4 Facilities}
        \label{FigVIa}
    \end{subfigure}%
    
            \begin{center}
        \begin{subfigure}[b]{0.33\textwidth}
            \includegraphics[width=\textwidth]{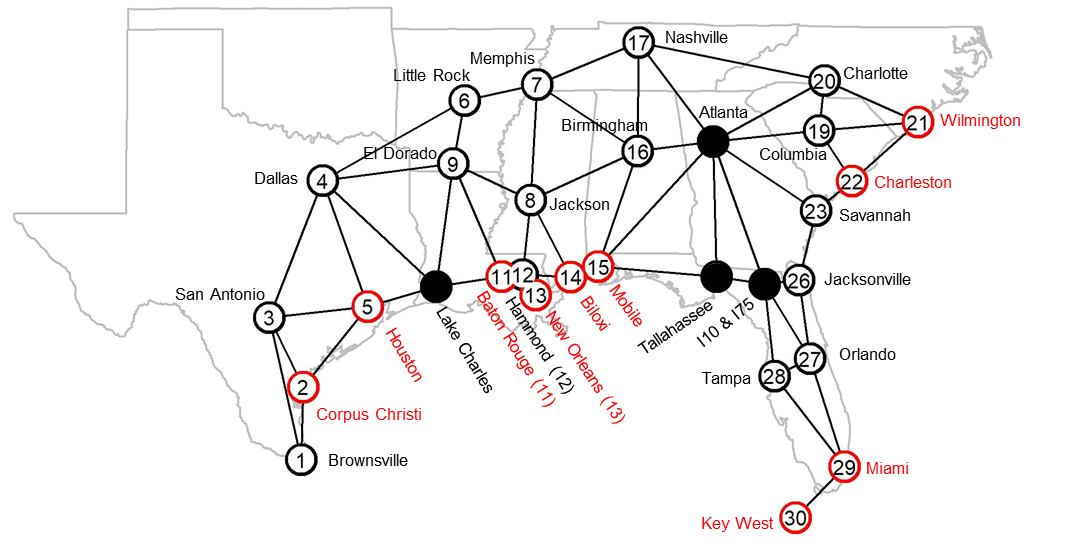}
      \caption{Trade (0.3), 4 Facilities}
      \label{FigVIb}
    \end{subfigure}%
    \begin{subfigure}[b]{0.33\textwidth}
            \includegraphics[width=\textwidth]{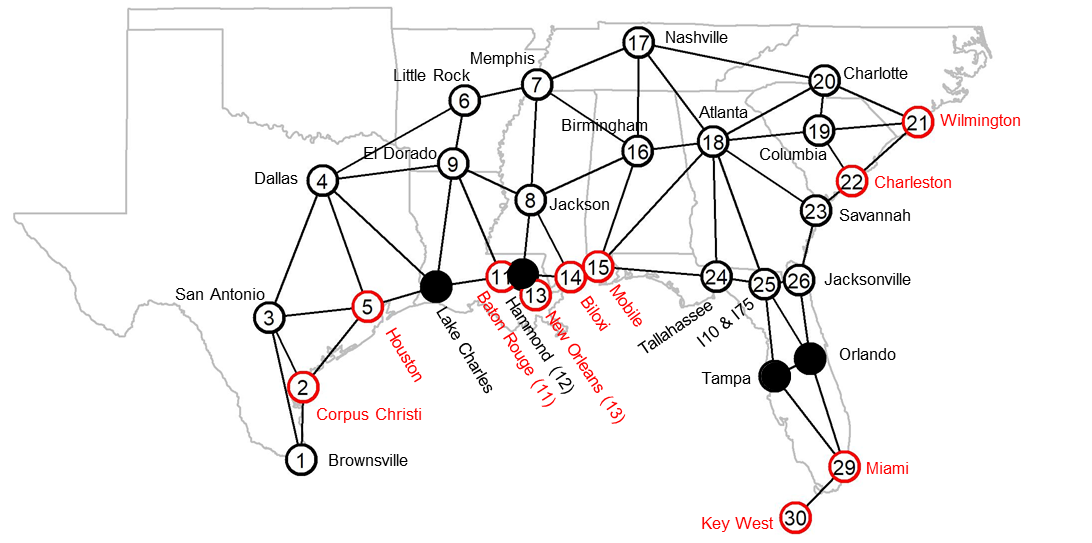}
      \caption{DRO, 4 Facilities}
      \label{FigVIb}
    \end{subfigure}%
        \end{center}
        \caption{Comparison of the optimal large facility locations for (nMinor,nMajor)=(2,1).}\label{Fig_Loc_21_large}
\end{figure}

\begin{figure}[t!]
     \begin{subfigure}[b]{0.33\textwidth}
 \centering
        \includegraphics[width=\textwidth]{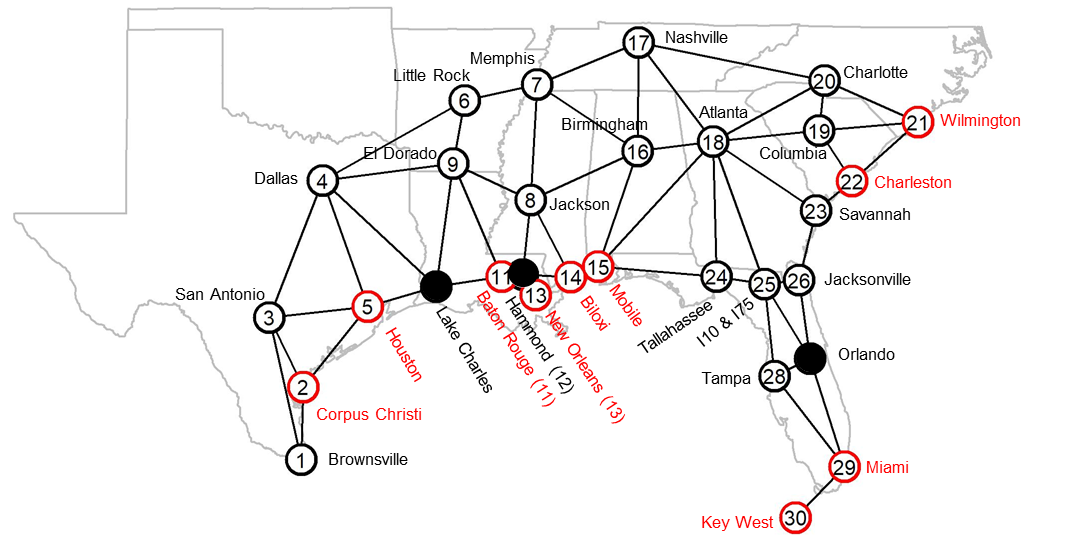}
        \caption{SP, 3 Facilities}
        \label{FigVIa}
    \end{subfigure}%
       \begin{subfigure}[b]{0.33\textwidth}
 \centering
        \includegraphics[width=\textwidth]{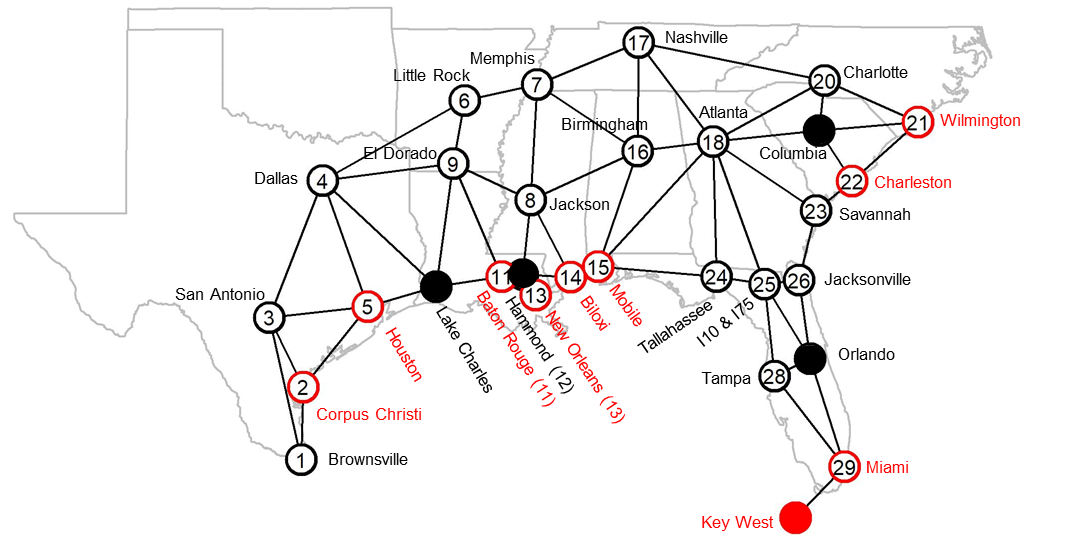}
        \caption{Trade (0.7), 5 Facilities}
        \label{FigVIa}
    \end{subfigure}%
       \begin{subfigure}[b]{0.33\textwidth}
 \centering
        \includegraphics[width=\textwidth]{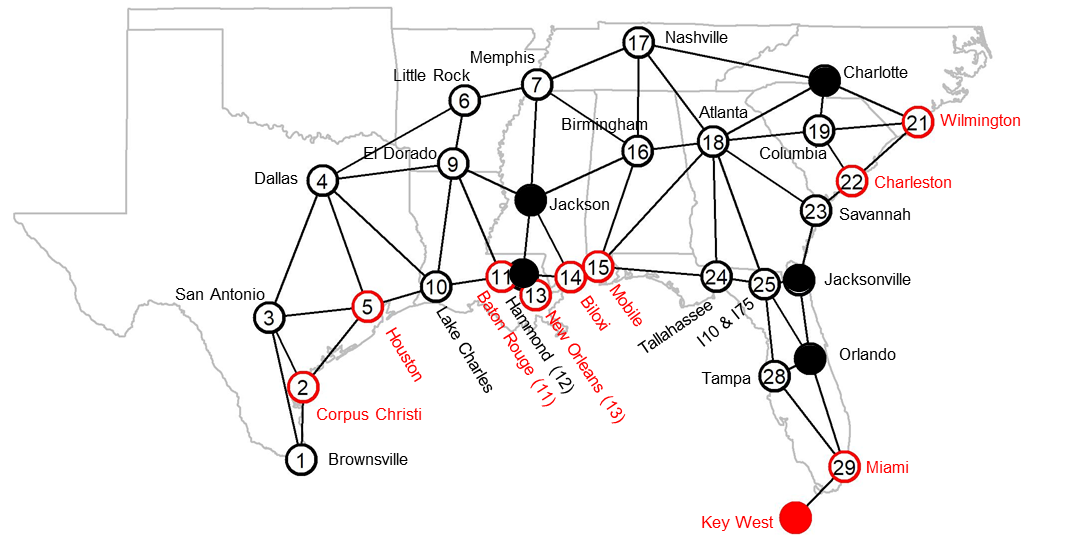}
        \caption{Trade (0.5), 6 Facilities}
        \label{FigVIa}
    \end{subfigure}%
    
            \begin{center}
        \begin{subfigure}[b]{0.33\textwidth}
            \includegraphics[width=\textwidth]{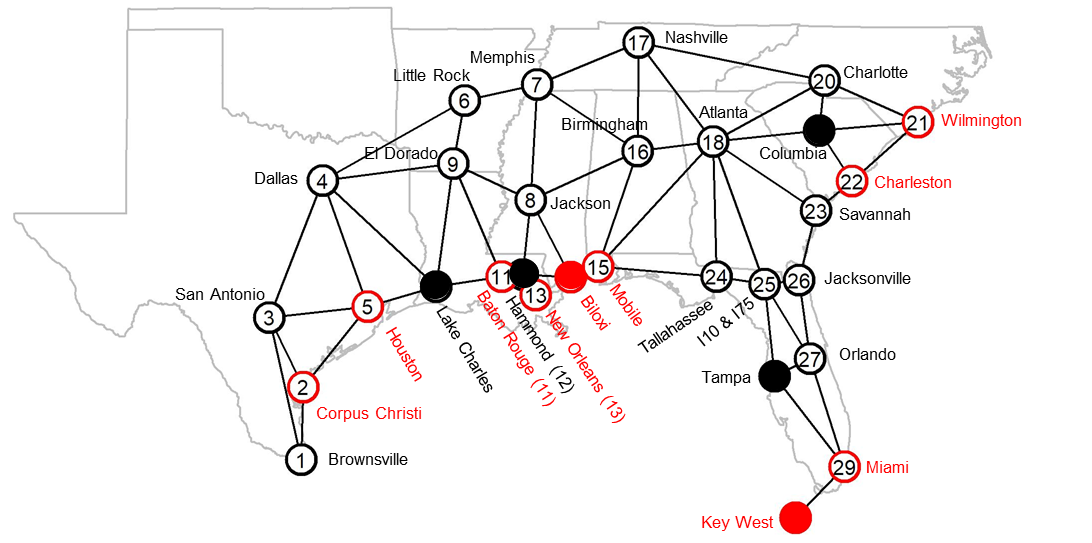}
      \caption{Trade (0.3), 6 Facilities}
      \label{FigVIb}
    \end{subfigure}%
    \begin{subfigure}[b]{0.33\textwidth}
            \includegraphics[width=\textwidth]{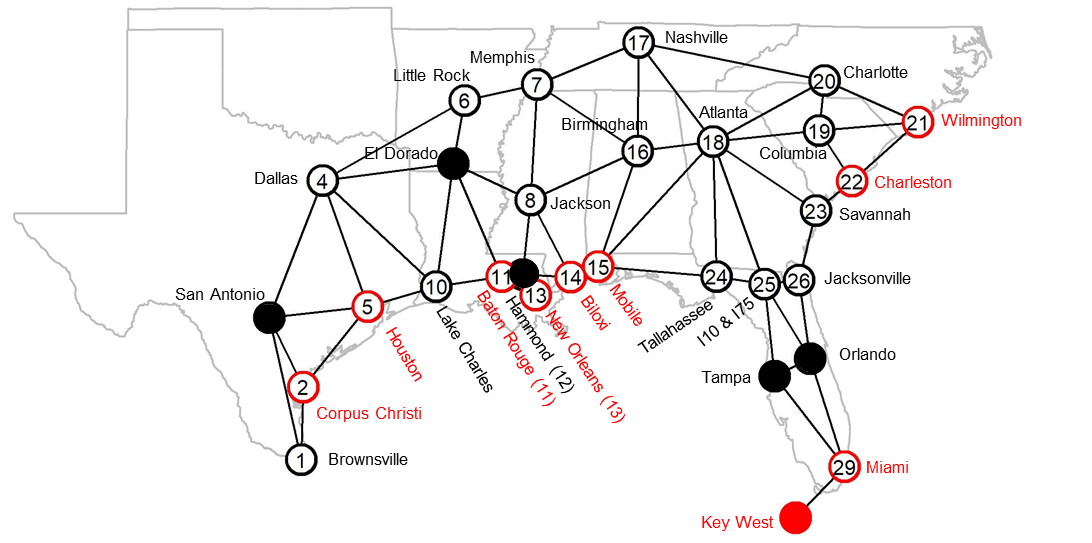}
      \caption{DRO, 6 Facilities}
      \label{FigVIb}
    \end{subfigure}%
        \end{center}
        \caption{Comparison of the optimal large facility locations for (nMinor,nMajor)=(4,2).}\label{Fig_Loc_4_2_large}
\end{figure}

\noindent In this section, we compare the optimal  (first stage) prepositioning  decisions for the hurricane season yielded by the DRO, SP, and trade-off models. The trade-off (denoted as Trade henceforth) considers both the DRO and SP recourse objectives that are weighted by a parameter $\theta$ \black{(see Section~\ref{sec:Trade})}. We present results for $\theta=0.3, 0.5, 0.7$ (low weight on the DRO objective; equally weighted objectives; and high weight on the DRO objective, respectively). We run the experiments twice; first with medium-sized facilities ($S_i$= 408,200, $f_i=\$ 188,400$) and second with large facilities ($S_i$= 780,000, $f_i=\$ 300,000$). For the SP model, we first used the MCO Algorithm in Section~\ref{sec:MCO} to determine an appropriate sample size to use for SAA. Based on the results presented in Tables~\ref{table:SampleSize} in \ref{Appex:MCO_results}, we used a sample size of $N=100$. We present the results of the optimal location of relief facilities/warehouses in Figures \ref{Fig_Loc_21_medium}--\ref{Fig_Loc_4_2_large}.  We present the results for the  amount of prepositioned relief supplies in each location in Tables \ref{table:OptimalLocations_Medium_2_1}--\ref{table:OptimalLocations_Medium_4_2} (medium-sized facilities with (nMinor, nMajor)=(2,1) and (4,2)) and Table \ref{table:OptimalLocations_Large_2_1}--\ref{table:OptimalLocations_Large_4_2} (large-sized facilities with (nMinor, nMajor)=(2,1) and (4,2)).

We make the following observations from these results. First, the SP model opens fewer facilities and prepositions a smaller quantity of each relief item than all models, which, as we show later, leads to a poor post-disaster performance. Moreover, the SP model tends to open more facilities at locations close to the Atlantic coast (e.g., Orlando (27), Tallahassee (24),  Savannah (23) when preparing for 4 minor and 2 major disasters with medium-sized facilities; see Figure~\ref{FigVIaMedium}). The SP model locates few warehouses near the potential landfall nodes in the Gulf of Mexico.


\textcolor{black}{In contrast, the DRO and Trade models better distribute warehouses near potential landfall nodes across the network.} \textcolor{black}{This may better hedge against the uncertainty of landfall locations. We will later observe this in the post-disaster performance; the number of warehouses and their spread tends to improve post-disaster distribution even when landfall is not well-forecast. We also observe that some facility locations are used in most of solutions. These include} \textcolor{black}{node 12 (Hammond, LA), node 10 (Lake Charles, LA), node 27 (Orlando, FL), node 30 (Key West, FL), and node 28 (Tampa, FL). It makes sense to preposition relief items at nodes 10, 12, 27, and 28 as they are not disaster-prone nodes (i.e., hurricane-safe) and are very close to several potential landfall nodes. In particular, node 10 and 12 are respectively close to potential landfall nodes (2,5,11) and  (11, 13-15) in the Gulf of Mexico, and nodes 27-28 are close to potential landfall nodes 29-30.}

\textcolor{black}{Second, to mitigate uncertainty and ambiguity, the DRO  model opens a larger number of facilities and prepositions a larger quantity of each relief item than the SP model (for either medium or large facilities)}. By doing so, the DRO model satisfies a larger amount of demand (reflected by significantly smaller  shortage cost in Table~\ref{table:InSample} and Figures~\ref{Fig1:Out_delta0}-\ref{Fig3:Out_delta50} presented later in section~\ref{sec5:OutSample}) and procures fewer relief items post-disaster (see Section~\ref{sec5:OutSample}).

The Trade model opens fewer (more) facilities and prepositions smaller (larger) amounts of relief items than the DRO (SP) model. This results in satisfying a smaller (larger) demand and procuring a larger (smaller) amount of relief items than the DRO (SP) model (see Section~\ref{sec5:OutSample}). Fourth, we observe that the Trade model with $\theta=0.3$ and $\theta=0.5$ opens the same number or more facilities than when  $\theta=$0.7. For example, when the facility size is medium and (nMinor, nMajor)=(4,2), Trade (0.3), Trade (0.5), Trade (0.7) opens 9, 9, and 8 facilities, respectively. Fifth, we observe that all models open fewer large facilities than medium-sized facilities. This makes sense as facilities have a larger capacity and are more expensive in the former case. 

The results in this section demonstrate how different approaches yield potentially different prepositioning plans. We next analyze the operational performance of these solutions.

\begin{table}[t!]
\color{black}
 \center 
 \footnotesize
   \renewcommand{\arraystretch}{0.2}
  \caption{The amount and locations of prepositioned relief supplies. Medium Facility (nMin,  nMaj)=(2,1)}
\begin{tabular}{lclcllllccclclllll}
\hline
Model & Open & Location & Water  & Food& \textcolor{black}{M-Kits}    \\
\hline 
	 SP 	&	2	&	  \black{Jacksonville (26)}       	&	       \black{2,602}   	&	       \black{383      }	&	0	\\
		 		 	&	 	&	 \black{Orlando (27)}     	&	       \black{948}     	&	       \black{3,124}   	&	       \black{9,283} 	\\  \cline{4-6}
		 		  &	& &	3,550	&	3,507	&	9,283	\\
 		 		 		 		 		 		 		 	\\												
	 DRO 	&	7	&	  \black{Houston (5)}     	&	       -       	&	       \black{210}     	&	       -       	\\
 	&	 	&	\black{Jackson (8)}       	&	       -       	&	       \black{3,754}   	&	       - 	\\
		 	&	 	&	\black{Lake Charles (10)} 	&	       \black{2,823}   	&	       -       	&	       - 	\\
	 	&	 	&	\black{Hammond (12)}      	&	       \black{1,991}   	&	       \black{538}     	&	       \black{65,122}  	\\
	 	&	 	&	\black{Baton Rouge (11)}  	&	7637	&	       \black{9,400}   	&	       \black{65,122}  	\\
		 	&	 	&	\black{Tallahassee (24)}  	&	       -       	&	       \black{4,899}   	&	       -       	\\
		 	&	 	&	\black{Key West (30)}     	&	       \black{2,823}   	&	-      	&	       -       	\\ \cline{4-6}
		 	&	 	&	  	&	15,273	&	18,800	&	130,243	\\
		 		 		 		 		 		 		 	\\
	 Trade (0.3) 	&	5	&	  \black{Hammond (12)}        	&	0	&	       \black{4,202}   	&	       \black{2,881}    	\\
	 	&	 	&	\black{Columbia (19)}     	&	       \black{2,497}   	&	       \black{566}     	&	0	\\
	 	&	 	&	\black{Tallahassee (24)}  	&	0	&	       \black{4,303}   	&	       \black{42,766}  	\\
	 	&	 	&	\black{I10\&I75 (25)}   	&	       \black{2,448}   	&	       \black{608}     	&	       \black{3,054}    	\\
		 	&	 	&	 \black{Jacksonville (26)}        	&	       \black{2,823}   	&	0	&	0	\\ \cline{4-6}
	 	&		&	  	&	7,768	&	9,679	&	48,701	\\ 
		 		 		 		 		 		 		 	\\
 Trade (0.5) 	&	5		&  \black{Lake Charles (10)}     	&	       \black{1,028}   	&	       \black{1,993}   	&	       \black{6,025}    	\\
		 	&	 	&	 \black{Hammond (12)}     	&	       \black{688}     	&	       \black{3,398}   	&	       \black{22,057}   	\\
	 	&	 	&	\black{Biloxi (14)}       	&	       \black{488}     	&	       -       	&	       \black{1,615}  	\\
	 	&	 	&	\black{Tallahassee (24)}  	&	       \black{2,823}   	&	       -       	&	       - 	\\
		 	&	 	&	\black{Key West (30)}     	&	       \black{2,282}   	&	       \black{866}     	&	       \black{5,205}  	\\ \cline{4-6}
	 	&	 		& 	&	       \black{7,309}   	&	       \black{6,257}   	&	       \black{34,902}  	\\
		 		 		 		 		 		 		 	\\													
 Trade (0.7) 	&	5	&	 \black{Corpus Christi (2)}    	&	       -       	&	       -       	&	       \black{9,994}    	\\
	 	&	 	&	\black{Hammond (12)}      	&	       \black{1,528}   	&	       \black{2,235}   	&	       \black{866}  	\\
	 	&	 	&	\black{Wilmington (21)}   	&	       -       	&	       -	&	      \black{3,099}    	\\
	 	&	 	&	\black{Savannah (23)}    	&	       \black{1,130}   	&	       \black{2,825}   	&	       \black{8,085} 	\\
		 	&	 	&	\black{Key West (30)}    	&	       \black{1,914}   	&	71	&	       \black{13,476} 	\\ \cline{4-6}
	 	&	 	&	         	&	4,572	&	5,131	&	35,520	\\
           \hline 
           \end{tabular}
\label{table:OptimalLocations_Medium_2_1}
\end{table} 
\begin{table}[t!]
\color{black}
 \center 
 \footnotesize
   \renewcommand{\arraystretch}{0.2}
  \caption{The amount and locations of prepositioned relief supplies. Medium Facility (nMin,  nMaj)=(4,2)}
\begin{tabular}{lclcllllccclclllll}
\hline
Model & Open & Location & Water  & Food& \textcolor{black}{M-Kits}    \\
\hline
 SP 	&	5	&	  \black{Lake Charles (10)}  	&	1099	&	2961	&	       2179    \\	
 	&	 	&	 \black{Hammond (12)}   	&	2589	&	400	&	       445 \\ 	
 	&	 	&	\black{Savannah (23)}   	&	1796	&	1658	&	       8968    \\	
 	&	 	&	\black{Tallahassee (24)}      	&	1459	&	2189	&	       12822   \\	
 	&	 	&	 \black{Orlando (27)}        	&	1810	&	1641	&	       8328    \\ \cline{4-6}	
 	&	 	&	 	&	8,753	&	8,849	&	       32,742  \\   	
										\\	
 DRO	&	10	&	 \black{San Antonio(3)}        	&	       \black{526}     	&	       \black{1,667}   	&	       \black{35,200}\\	
 	&	 	&	 Little Rock(6)    	&	       -       	&	       \black{2,918}   	&	       -       \\	
 	&	 	&	\black{El Dorado (9)}      	&	       -       	&	       \black{4,899}   	&	       -\\	
 	&	 	&	 \black{Lake Charles (10)}   	&	       \black{2,126}   	&	       \black{1,208}   	&	       0\\	
 	&	 	&	 \black{Hammond (12)} 	&	       \black{1,980}   	&	       \black{1,463}   	&	       0\\	
 	&	 	&	\black{Biloxi (14)}  	&	       -       	&	       \black{694}     	&	       -\\	
 	&	 	&	 \black{Tallahassee (24)}    	&	       \black{1,041}   	&	       \black{2,266}   	&	       \black{59,355}\\	
 	&	 	&	\black{Orlando (27)}      	&	       \black{2,822}   	&	       -       	&	       -\\	
 	&	 	&	\black{Tampa (28)}     	&	       \black{1,374}   	&	2514	&	       0\\	
 	&	 	&	\black{Key West (30)}  	&	       \black{1,949}   	&	       \black{1,072}   	&	       \black{31,941}\\  \cline{4-6}	
 	&	 	&	 	&	11,818	&	18,701	&	126,496\\	
										\\	
 Trade (0.3) 	&	9	&	  \black{San Antonio (3)} 	&	0	&	4,871	&	       1,973   \\	
 	&	 	&	 \black{Jackson (8)} 	&	0	&	4,426	&	       33,916   \\	
 	&	 	&	 \black{Lake Charles (10)}   	&	1,859	&	1,665	&	       570     \\	
 	&	 	&	\black{Hammond (12)}        	&	2,370	&	785	&	       0       \\	
 	&	   	&	  \black{Columbia (19)}   	&	689	&	3,162	&	       38,890   \\	
 	&	   	&	 \black{Jacksonville (26)} 	&	1,551	&	2,063	&	       10,326  \\	
 	&	 	&	 \black{Orlando (27)}  	&	2,823	&	0	&	       0       \\	
 	&	 	&	 \black{Tampa (28)}    	&	2,823	&	0	&	       0       \\	
 	&	 	&	 \black{Key West (30)} 	&	1,471	&	2,273	&	       5,200   \\   \cline{4-6}	
 	&	 	&	 	&	13,587	&	19,245	&	       90,875  \\	
										\\	
  Trade (0.5) 	&	9	&	 \black{Corpus Christi (2)}    	&	1,010	&	       -       	&	       15,438   \\	
 	&	 	&	\black{Hammond (12)}  	&	2823	&	       -       	&	0	\\
 	&	 	&	       \black{Biloxi (14)}      	&	2166	&	0	&	       1529    \\	
 	&	 	&	\black{Wilmington (21)}      	&	614	&	0	&	       0       \\	
 	&	 	&	       \black{I10\&I75 (25)}    	&	869	&	3373	&	       1218    \\ 	
 	&	 	&	 \black{Jacksonville (26)}     	&	0	&	4,620	&	       19995   \\      	
 	&	 	&	 \black{Orlando (27)}  	&	2,823	&	       -       	&	       -       \\	
 	&	 	&	 \black{Tampa (28)}    	&	0	&	4,662	&	       16,984  \\  	
 	&	 	&	 \black{Key West (30)} 	&	2,379	&	361	&	       29,387   \\  \cline{4-6}	
 	&	 	&	      	&	12,684	&	13,017	&	       84,550  \\	
										\\	
  Trade (0.7) 	&	8	&	 \black{Corpus Christi (2)}     	&	       -       	&	       -       	&	       20,629 \\	
 	&	 	&	 \black{Lake Charles (10)}        	&	830	&	3,406	&	       3,779   \\	
 	&	 	&	 \black{Hammond (12)} 	&	2,381	&	738	&	       2,059\\	
 	&	 	&	 \black{Biloxi (14)}     	&	2,606	&	0	&	       9,349   \\	
 	&	 	&	 \black{Wilmington (21)}        	&	423	&	0	&	       3,378   \\	
 	&	 	&	 \black{Jacksonville (26)}      	&	701	&	3,534	&	       10,639  \\	
 	&	 	&	 \black{Orlando (27)}  	&	1,284	&	2,548	&	       8,743 \\	
 	&	 	&	  \black{Key West (30)}        	&	2,645	&	0	&	       22,137\\ \cline{4-6} 	
 	&	 	&	  	&	10,871	&	10,226	&	       80,713 \\ \hline 	
           \end{tabular}
\label{table:OptimalLocations_Medium_4_2}
\end{table} 

\begin{table}[t!]
\color{black}
 \center 
 \footnotesize
   \renewcommand{\arraystretch}{0.5}
  \caption{The amount and locations of prepositioned relief supplies. Large Facility (nMin, nMaj)=(2,1)}
\begin{tabular}{lclcllllccclclllll}
\hline
Model & Open & Location & Water  & Food& \textcolor{black}{M-Kits}    \\
\hline
 SP  	&	2	&	 \black{Hammond (12)} 	&	2,460	&	2,394	&	3,729	\\
 	&	 	&	  \black{Orlando (27)}    	&	1,710	&	1,623	&	9,717	\\ \cline{4-6}
 	&	 	&	 	&	4,170	&	4,017	&	13,446	\\
 		 		 		 		 		 	\\
											
 DRO 	&	4	&	 \black{Lake Charles (10)} 	&	481	&	       -       	&	       - 	\\
 	&	 	&	  \black{Hammond (12)}   	&	5,394	&	       -       	&	       - 	\\
 	&	 	&	 \black{Orlando (27)}     	&	       -       	&	9,360	&	       - 	\\
 	&	 	&	 \black{Tampa (28)}       	&	-      	&	791	&	39,510	\\ \cline{4-6}
 	&	 	&	 	&	5,875	&	10,152	&	39,510	\\
 		 		 		 		 		 	\\

  Trade (0.3) 	&	4	&	 \black{Lake Charles (10)}  	&	205	&	4,026	&	       - 	\\
 	&	 	&	\black{Atlanta (18)}     	&	       -       	&	566	&	1,581	\\
 	&	 	&	 \black{Tallahassee (24)}        	&	2,086	&	5,085	&	47,120	\\
 	&	 	&	 \black{I10\&I75 (25)}  	&	5,394	&	-	&	       -        	\\ \cline{4-6}
 	&	 	&	        	&	7,685	&	9,677	&	48,701	\\
 		 		 		 		 		 	\\

  Trade (0.5) 	&	4	&	  \black{Hammond (12)}      	&	1,116	&	1736	&	-  	\\
 	&	 	&	 \black{Biloxi (14)}    	&	823	&	656	&	16,884	\\
 	&	 	&	 \black{Atlanta (18)}  	&	       19,72   	&	3,789	&	13,082	\\
 	&	 	&	 \black{Key West (30)}        	&	3,500	&	314	&	10,980	\\ \cline{4-6}
 	&	 	&	    	&	7,411	&	6,495	&	40,946	\\
 		 		 		 		 		 	\\

  Trade (0.7) 	&	3	&	 \black{Corpus Christi (2)}    	&	       -       	&	       -       	&	9,852	\\
 	&	 	&	      \black{Hammond (12)}     	&	2,386	&	5,050	&	12,192	\\
 	&	 	&	\black{Key West (30)}       	&	2,283	&	66	&	13,260	\\ \cline{4-6}
 	&	 	&	        	&	4,669	&	5,117	&	35,304	\\
           \hline 
\end{tabular} 
\label{table:OptimalLocations_Large_2_1}
\color{black}
\end{table}

\begin{table}[t!]
\color{black}
 \center 
 \footnotesize
   \renewcommand{\arraystretch}{0.5}
  \caption{The amount and locations of prepositioned relief supplies. Large Facility (nMin, nMaj)=(4,2)}
\begin{tabular}{lclcllllllllll}
\hline
Model & Open & Location & Water  & Food& \textcolor{black}{M-Kits}  &  \\
\hline
 SP 	&	3	&	 \black{Lake Charles (10)}  	&	1,234	&	4,073	&	       11,387&\\	
 	&	 	&	  \black{Hammond (12)}   	&	3810	&	2619	&	       9304    \\  	
 	&	 	&	 \black{Orlando (27)}  	&	4,037	&	2,187	&	       12,075  \\ \cline{4-6} 	
 	&		&	 	&	9,081 &		8,879 &		32,766 \\	
										\\	
 DRO	&	6	&	 \black{San Antonio (3)}      	&	5,394	&	       -       	&	       -       \\	
 	&	 	&	 El Dorado (9)  	&	5,394	&	       -       	&	       - \\	
 	&	 	&	 \black{Hammond (12)}       	&	323	&	7,111	&	       121,251 \\  	
 	&	 	&	 \black{Orlando (27)}  	&	       -       	&	9,360	&	       -       \\ 	
 	&	 	&	 \black{Tampa (28)}  	&	       -       	&	533	&	       - \\ 	
 	&	 	&	  \black{Key West (30)}       	&	2,753	&	1,360	&	       - \\  \cline{4-6} 	
 	&	 	&	 	&	13,864	&	18,364	&	       121,251\\	
										\\	
  Trade(0.3) 	&	7	&	 \black{Lake Charles (10)}     	&	2,305	&	5,361	&	-	\\
 	&	 	&	 \black{Hammond (12)}   	&	5,394	&	       - 	&	-	\\
 	&	 	&	 \black{Biloxi (14)}  	&	       -       	&	1,386	&	       6,035  \\	
 	&	 	&	 \black{Columbia (19)}       	&	2,917	&	1,693	&	       - \\	
 	&	 	&	  \black{Orlando (27)} 	&	2,294	&	4,849	&	       38,152  \\	
 	&	 	&	 \black{Tampa (28)}    	&	       -       	&	4,501	&	       44,674 \\	
 	&	 	&	 \black{Key West (30)}        	&	325	&	1,816	&	       2,947\\   \cline{4-6} 	
 	&	 	&	  	&	13,235	&	19,606	&	       91,808\\	
										\\	
 Trade(0.5) 	&	6	&	 \black{Jackson (8)}        	&	       -       	&	8,201	&	       47,435  \\	
 	&	 	&	      \black{Hammond (12)}     	&	4,468	&	1,583	&	       1,652\\	
 	&	 	&	 \black{Charlotte (20)}        	&	1,843	&	1,517	&	       528\\	
 	&	 	&	 \black{Jacksonville (26)}     	&	151	&	1,444	&	       1,345   \\	
 	&	 	&	 \black{Orlando (27)}  	&	5394	&	       -       	&	       -       \\	
 	&	 	&	 \black{Key West (30)}   	&	1,827	&	271	&	       31,282\\  \cline{4-6}	
 	&	 	&	 	&	13,683	&	13,017	&	       82,241 \\	
										\\	
 Trade(0.7) 	&	5	&	  \black{Lake Charles (10)}  	&	1,241	&	4,262	&	       17,860\\	
 	&	 	&	 \black{Hammond (12)}  	&	4,932	&	754	&	       3,459   \\	
 	&	 	&	 \black{Columbia (19)}        	&	1,086	&	1,780	&	       11,484 \\ 	
 	&	 	&	 \black{Orlando (27)}        	&	3,242	&	3,517	&	       15,602\\	
 	&	 	&	 \black{Key West (30)}      	&	       -       	&	       -       	&	       29,117 \\ \cline{4-6}	
 	&	 	&	 	&	10,501	&	10,313	&	       77,522& \\	
 	\hline
\end{tabular} 
\label{table:OptimalLocations_Large_4_2}
\color{black}
\end{table}

\subsubsection{Analysis of optimal solutions quality}\label{sec5:OutSample}

\noindent In this section, we compare how the optimal solutions to the DRO, SP, and Trade models perform. We consider two cases: \textcolor{black}{when the underlying uncertainty distributions used in the optimization are perfectly specified and misspecified}. For brevity, we present results for (nMinor, nMajor)= (4,2) disasters. We evaluate the operational performance of the models as follows. First, we fix the optimal first-stage decisions ($\pmb{o,z}$) yielded by each model in the SP model. Then, we solve the second-stage recourse problem in \eqref{2ndstage} using  ($\pmb{o,z}$) and the following two sets of $N'=10,000$ out-of-sample data ($\pmb{q^n, d^n, M^n}$), for all $n \in [N']$, to compute the corresponding second-stage procurement, shortage, holding, and transportation costs. To generate the $N'$ out of sample data points:
\begin{enumerate}\itemsep0em
\item[Set 1.] \textit{Perfect distributional information}: we use the same parameter settings and distributions that we use for generating the $N$ in-sample data points to generate the $N^\prime$ data points. 
\item[Set 2.] \textit{Misspecified distributional information}: we vary the distribution type of random parameters to generate the $N^\prime$ data \citep{wang2020distributionally}. That is, we perturb the distribution of the demand by a parameter $\Delta$ and obtain a parameterized uniform distribution [$(1-\Delta) d, (1+\Delta)d$ ], for which a higher value of $\Delta$ corresponds to a higher variation level (we similarly perturb the distribution of the other parameters).  We apply $\Delta \in \{ 0, 0.25, 0.5\}$. A zero value of $\Delta$ indicates that we do not perturb the demand's distribution, but rather we simulate the optimal solutions under a uniform distribution defined on the range of the demand (i.e., we vary the in-sample distribution). We generate 10,000 samples from these uniform distributions to test the performance of the optimal solutions obtained from the DRO, SP, and Trade models. This is to simulate the performance of the models' optimal decisions when the in-sample data is biased (i.e., true distributions are different).
\end{enumerate}
\begin{table}[t!]
 \center 
 \footnotesize
   \renewcommand{\arraystretch}{0.6}
  \caption{The fixed cost and out-of-sample costs  for (nMinor, nMajor)=(4,2) under perfect distributional information (LogN).}
\begin{tabular}{llllllllllllllllll}
\hline
Model 	&	Fixed	& 	Procurement	& 	Shortage 	& 	Holding	& 	Shipping	\\
\hline
SP	&	59,156,500	&	38,601,400	&	2,028,850	&	4,541,360	&	925,687	\\
											\\
Trade (0.7)	&	75,273,100	&	29,351,400	&	547,039	&	9,753,760	&	1003490	\\
											\\
Trade (0.5)	&	97,952,900	&	15,256,300	&	254,479	&	19,033,700	&	1154370	\\
											\\
Trade (0.3)	&	127,527,600	&	3,017,060	&	36,284	&	48,356,700	&	1153970	\\
											\\
DRO	&	128,604,862	&	4,463,470	&	78,920	&	48,577,500	&	1130980	\\
	\hline
	\end{tabular} 
\label{table:InSample}
\end{table}

Let us first analyze simulation results under Set 1 (i.e., perfect distributional information case). Table~\ref{table:InSample} presents the pre-disaster (fixed) planning cost and the average post-disaster (procurement, shortage, holding, and shipping) costs. We make the following observations. The SP model clearly has the lowest fixed cost (i.e., acquisition and facility opening costs) among all models, which makes sense because it opens fewer facilities and allocates less relief items. This, in turn, leads to a substantial shortage and the need to procure significantly larger quantities of relief items in the aftermath (reflected by the higher shortage and procurement costs).  In contrast, the DRO and Trade models satisfy more post-disaster demand by opening more facilities and allocating more relief items, resulting in higher one-time fixed costs (pre-disaster) but lower shortage and procurement costs.  Reducing shortage costs, in particular, translates to better support of vulnerable populations post-disaster.  This suggests there are benefits to using DRO or Trade models even in a perfect information setting, if decision-makers have sufficient budget. Lastly, we note that it is not surprising that the SP model's holding and transportation costs are lower than the other considered models since the SP model opens fewer facilities and prepositions a lower quantity of relief items (there is less available to ship to demand nodes). Moreover, we have assumed perfect information of the exact demand distribution in this simulation (which is highly unlikely in real disaster situation).

\begin{figure}[t!]
     \begin{subfigure}[b]{0.5\textwidth}
 \centering
        \includegraphics[width=\textwidth]{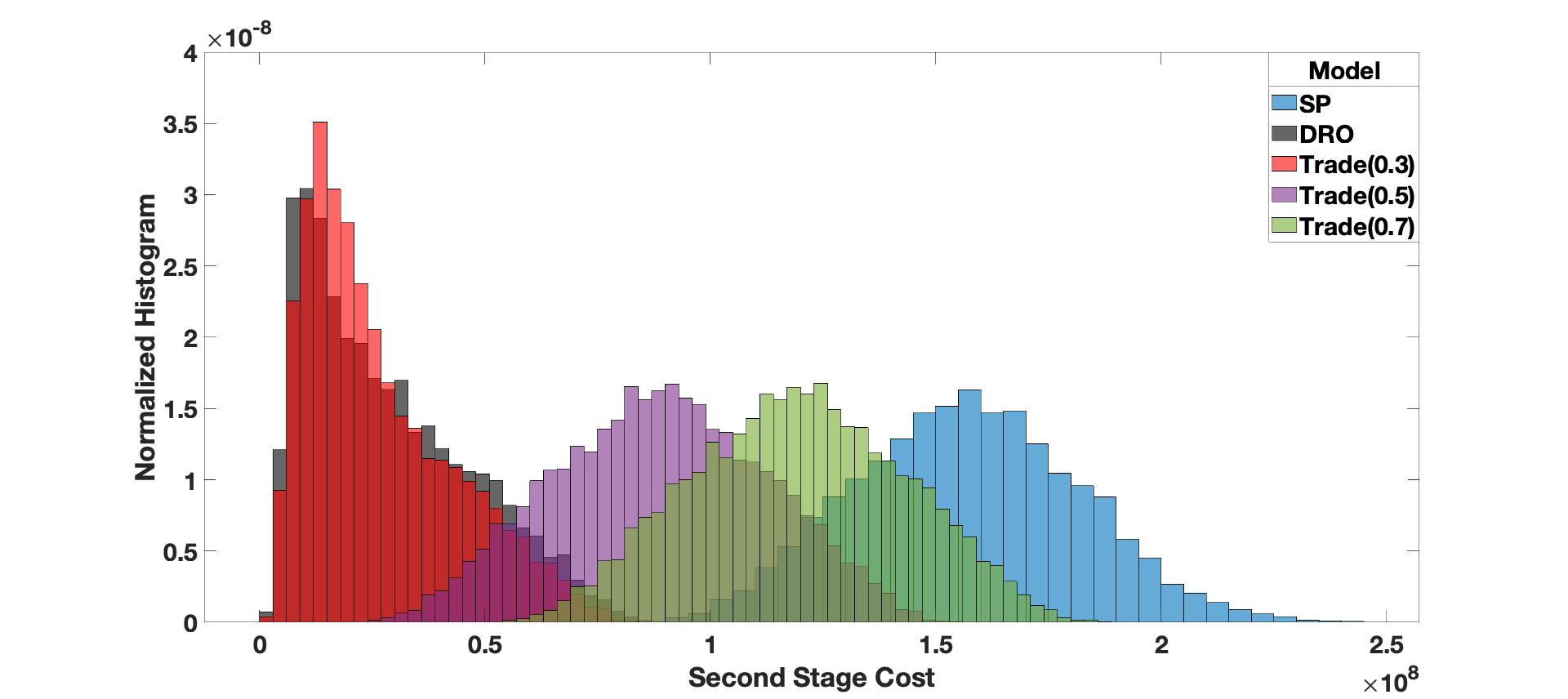}
        \caption{Second stage cost  }
        \label{Fig1a}
    \end{subfigure}%
    \begin{subfigure}[b]{0.5\textwidth}
            \includegraphics[width=\textwidth]{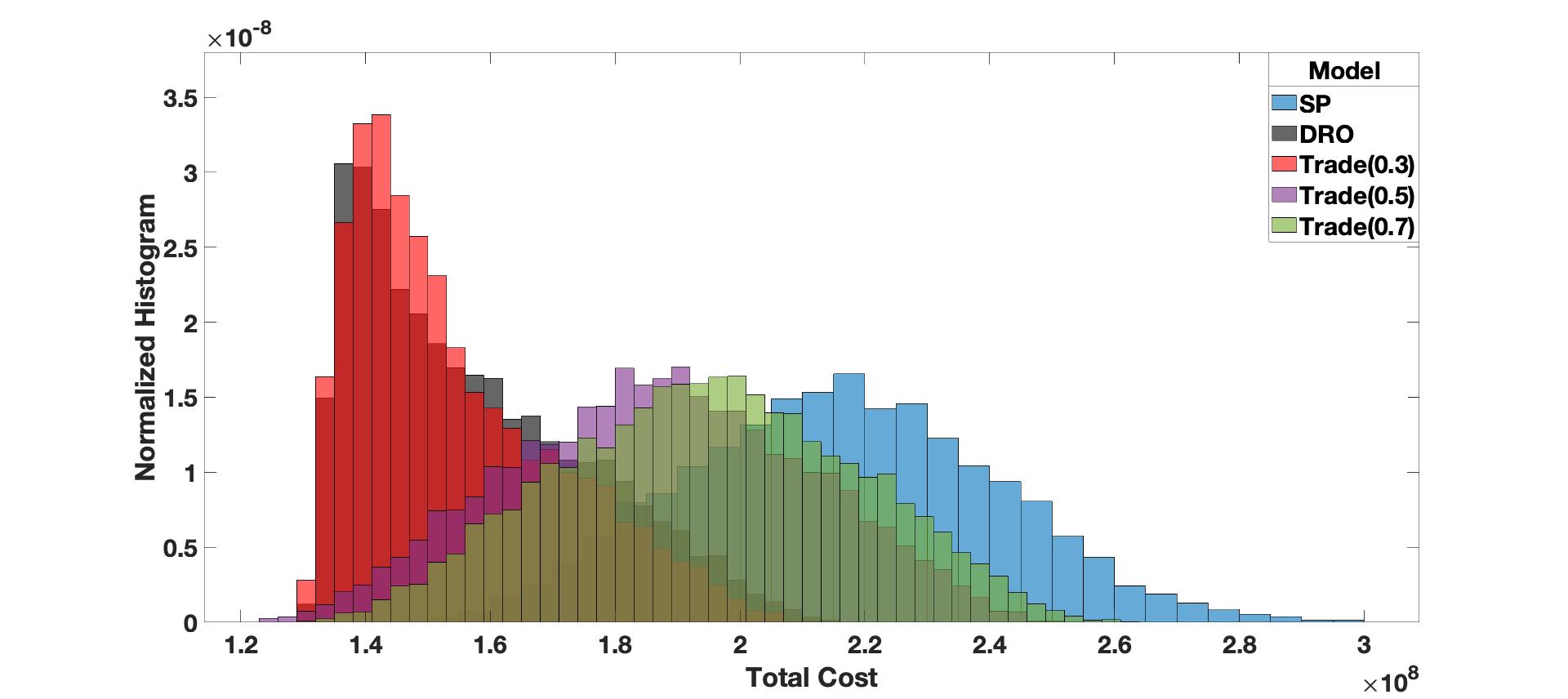}
      \caption{Total cost}
      \label{Fig1b}
    \end{subfigure}%
    \caption{Comparison of simulation results under missspecified distributional information, $\Delta=0$}\label{Fig1:Out_delta0}
\end{figure}

\begin{figure}[t!]
     \begin{subfigure}[b]{0.5\textwidth}
 \centering
        \includegraphics[width=\textwidth]{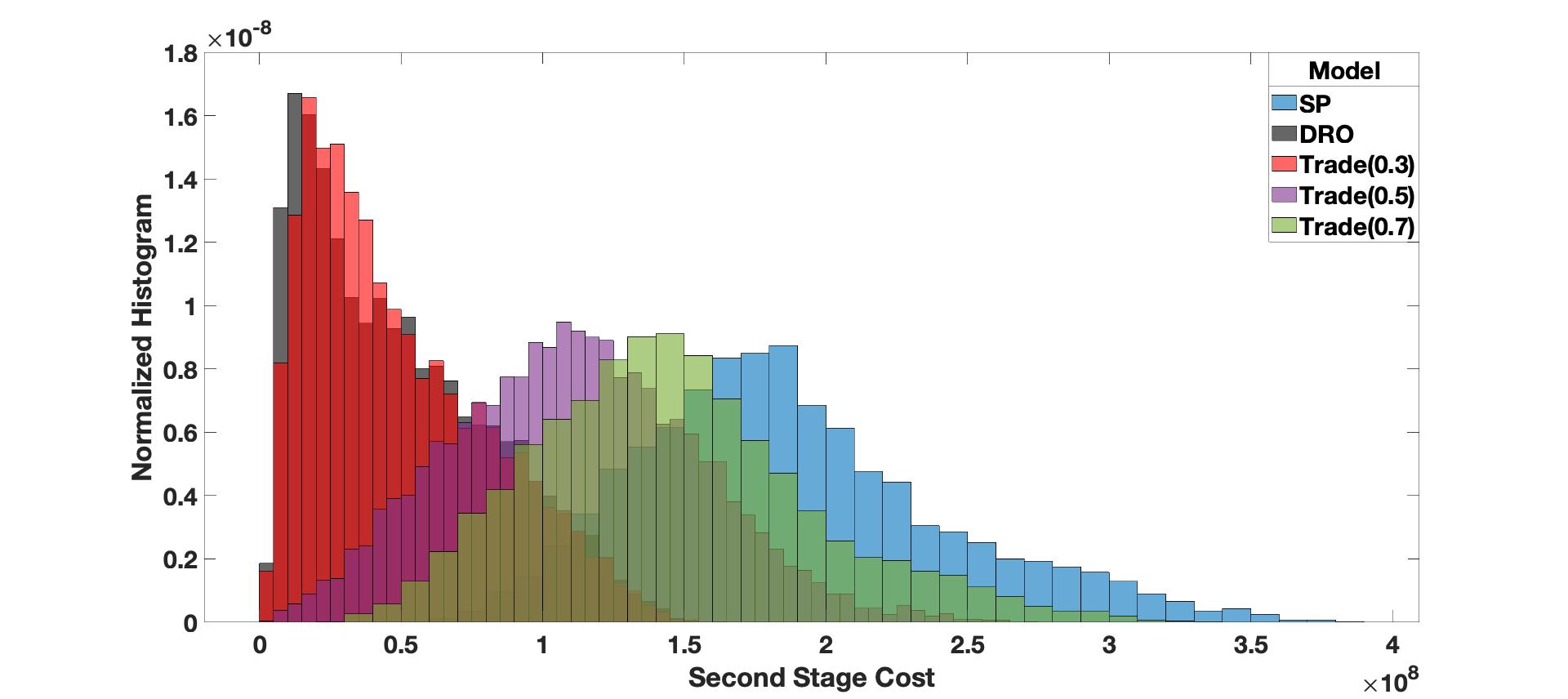}
        \caption{Second stage cost  }
        \label{Fig2a}
    \end{subfigure}%
    \begin{subfigure}[b]{0.5\textwidth}
            \includegraphics[width=\textwidth]{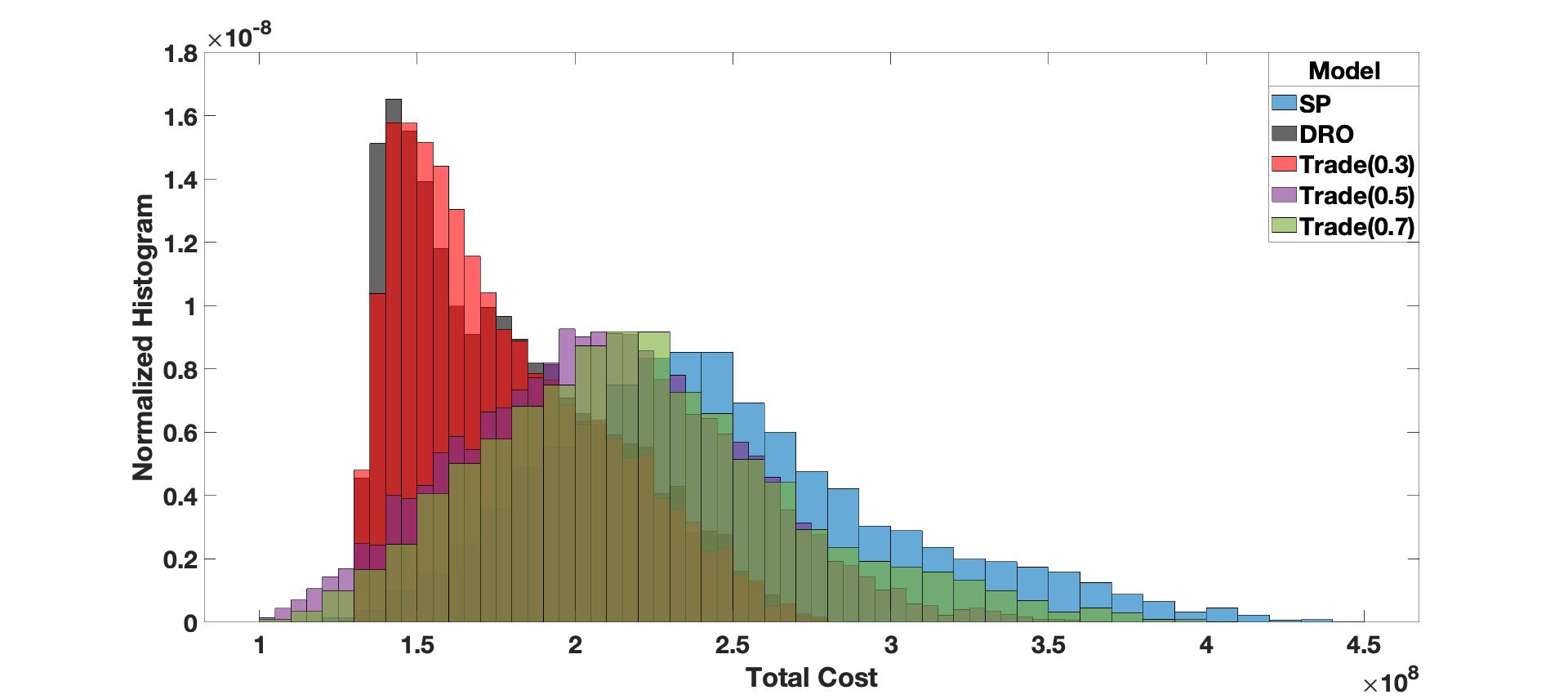}
      \caption{Total cost}
      \label{Fig2b}
    \end{subfigure}%
    
    \caption{Comparison of simulation results under missspecified distributional information, $\Delta=0.25$}\label{Fig2:Out_delta0.25}
\end{figure}

\begin{figure}[t!]
     \begin{subfigure}[b]{0.5\textwidth}
 \centering
        \includegraphics[width=\textwidth]{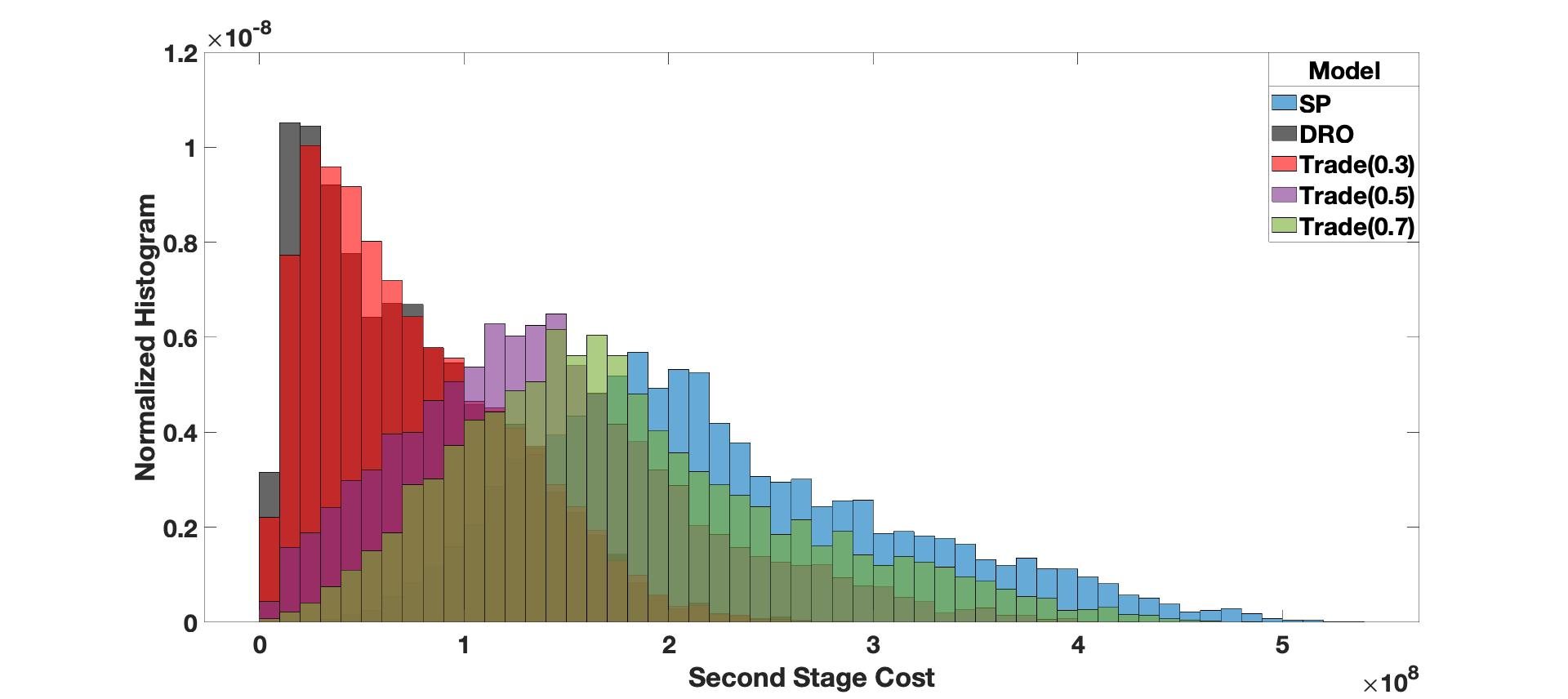}
        \caption{Second stage cost  }
        \label{Fig3a}
    \end{subfigure}%
    \begin{subfigure}[b]{0.5\textwidth}
            \includegraphics[width=\textwidth]{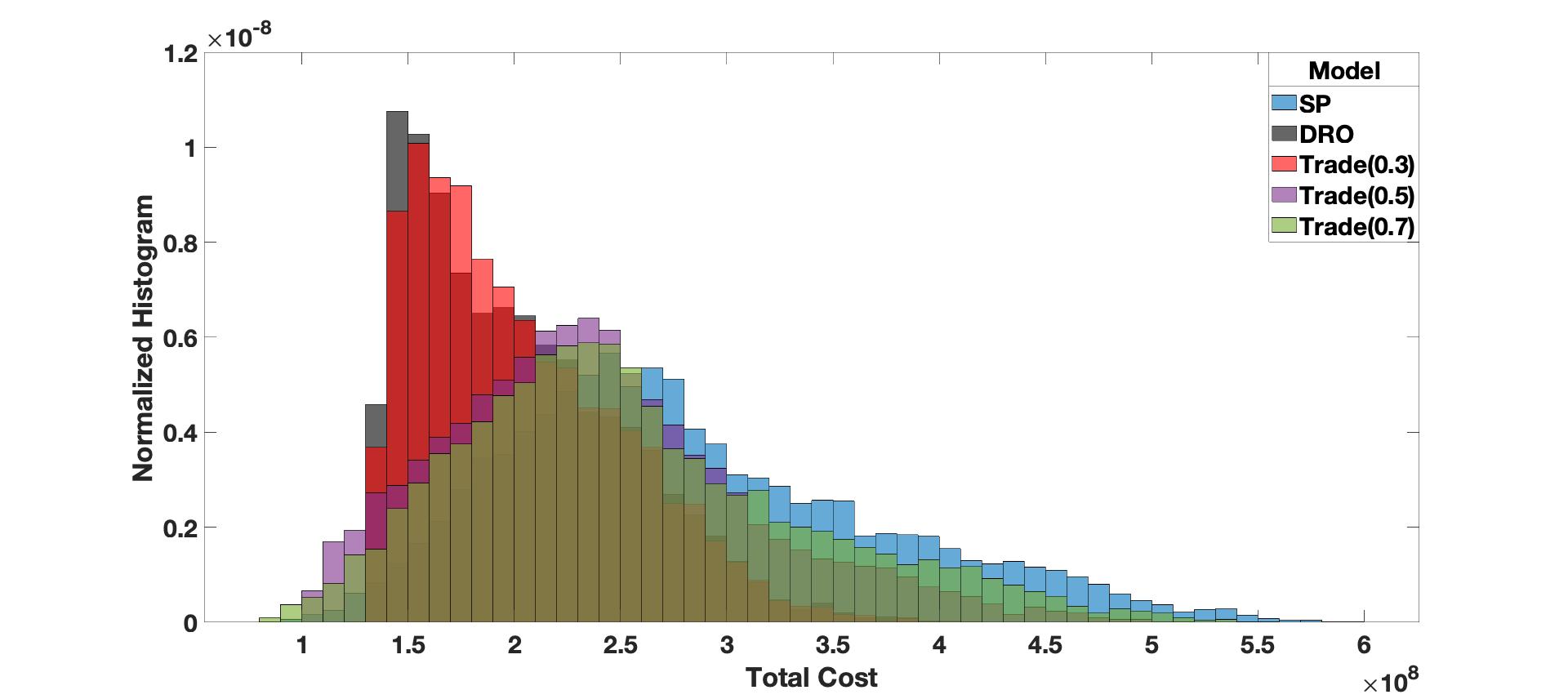}
      \caption{Total cost}
      \label{Fig3b}
    \end{subfigure}%
    \caption{Comparison of simulation results under missspecified distributional information, $\Delta=0.5$}\label{Fig3:Out_delta50}
\end{figure}

Let us now analyze the out-of-sample performance of the models under the case where the true distribution is different than the one used in the optimization (i.e., Set 2. misspecified distributional information case). Figures \ref{Fig1:Out_delta0}--\ref{Fig3:Out_delta50} present histograms of the out-of-sample objective value (i.e., total cost=fixed cost+second-stage cost) and second-stage cost (i.e., operational costs) for each of the three levels of variation, $\Delta$.

It is quite evident in the case of the misspecified distributions that the DRO and Trade (0.3) solutions consistently outperform the solutions from the SP, Trade (0.5), and Trade (0.7) models. This relationship holds for all levels of variation, $\Delta \in \{0, 0.25, 0.5\}$, and across the criteria of mean and quantiles of total and second-stage costs.  In particular, the DRO model has substantially lower shortage costs than the other considered models. For example, when $\Delta=0.25$, the average shortage costs of the DRO is 2,533. In contrast, the average shortage costs for the Trade (0.3), Trade (0.5), Trade (0.7), and SP are 1,021,290, 3,252,750, 10,478,500, and 34,592,400, respectively.  We also observed lower shortage costs of the DRO model as compared to SP model in the perfect information case, but it is particularly pronounced here. In addition, the DRO and Trade (0.3) models seem to be more stable in terms of attaining the lowest standard deviations (i.e., variation) in the total and second-stage costs.

As we mentioned earlier, although various types of early warning systems may have been set up, the pre-disaster estimates of the post-disaster damage and associated demand, for example, are often inaccurate. The superior performance of DRO and Trade (0.3) models, which focus on hedging against distributional ambiguity and uncertainty, reflects the value of modeling such potential distributional ambiguity of random parameters.


\subsection{Earthquake Case Study}\label{sec:Earthquack}

\noindent In this section, we consider a different disaster relief context: earthquakes. There is little-to-no forewarning when they occur, and transportation arcs and prepositioned inventory are often damaged. In addition, post-disaster procurement may be very difficult, though necessary, if inventory is unavailable. Thus, deciding where and how much inventory to locate is challenging but critical to effective response. These planning decisions are generally focused on a single event, rather than for an upcoming season of multiple disasters (as in the hurricane case study). Yet, the same multi-disaster modeling framework can be used to consider how different locations within the affected area have varying levels of damage.

\textcolor{black}{To study the performance of the proposed models in this context,} we use data from the 2010 earthquake that hit Yushu County in Qinghai Province, China presented in \cite{ni2018location}. This 7.1 magnitude earthquake caused large-scale social and economic destruction. Figure~\ref{Fig_map2} in \black{~\ref{Appx:Case2_data}} shows a diagram of the affected area, which consists of 13 nodes and 15 road links. \textcolor{black}{We model the varying post-disaster damage levels by following the outward ripple of damage. That is,} the nodes closest to the epicenter may receive the worst damage, whereas others may have fewer effects. \textcolor{black}{We use} this characteristic \textcolor{black}{to estimate} the post-disaster proportion of usable inventory, $\mu^{\tiny \rho}$ in Table~\ref{table:Case2_input} \black{in ~\ref{Appx:Case2_data}}. We follow the seismic intensity categories presented in \cite{ni2018location}. Table~\ref{table:Case2_input}  also presents data on the fixed, acquisition, shortage, and holding costs. As in \cite{ni2018location}, we set $\mu^{\mbox{\tiny d}}_i=100$, $\sigma^{\mbox{\tiny d}}_i=10$, $\sigma^{\tiny \rho}=0.1$, $\mu^{\mbox{\tiny V}}=300$, and $\sigma^{\mbox{\tiny V}}=30$.  We let $\mu^{\mbox{\tiny M}}=\mu^{\mbox{\tiny d}}$, $\sigma^{\mbox{\tiny M}}=0.5\mu^{\mbox{\tiny d}}$, $[\ML, \MU]=[0.9\mu^{\mbox{\tiny M}}, 1.10\mu^{\mbox{\tiny M}}]$, and $[\rhoL, \rhoU]=[0.9 \mu^{\tiny \rho}, 1.10\mu^{\tiny \rho}]$. To approximate the lower and upper bounds on ($d, V$), we follow the same procedure described earlier in Section~\ref{sec5:Hurrican}. 

\subsubsection{\textcolor{black}{Analysis of optimal prepositioning decisions}}

\noindent In this section, we compare the optimal prepositioning decisions yielded by the DRO, SP, and Trade models presented in Table~\ref{table:OptimalLocations_Case2}. We first observe that each of the models open four facilities. Second, they all select the facility at node 11, though this is the only facility shared between the DRO and SP solutions. It makes sense to preposition items at node 11 because it has the highest usable fraction of prepositioned items after a disaster (i.e., it is less likely that the prepositioned items at this node will be destroyed). Third, the DRO and Trade models always open a facility at node 8, which makes sense because this node is close to an airport and also has a higher $\mu^{\tiny \rho}$ as compared to node 9, which the SP model selects ($\mu^{\tiny \rho}_8=0.7$ vs. $\mu^{\tiny \rho}_9$=0.6). The DRO and Trade (0.3) models select node 6, which has a direct path to the disaster's epicenter at node 1, while the other models select node 7, which does not directly link to node 1. Finally, we observe that the DRO and Trade (0.3) models allocate more relief items than the other considered models. As such, the DRO and Trade (0.3) models result in the higher acquisition costs.

\begin{table}[t]
 \center 
 \footnotesize
   \renewcommand{\arraystretch}{0.3}
  \caption{The amount and locations of prepositioned relief supplies. Earthquake case study }
\begin{tabular}{lccccccccc}
\hline
  Model & \# Open & Location & Amount & Fixed & Acquisition  \\
\hline
SP & 4 &7 &	136 &  659	& 3843\\
&& 9 &	566& \\
&&10&	462& \\
&& 11	&268& \\
\\
DRO 	&	4	&	6	&	800	&	735	&	6009	\\
	&		&	8	&	406	&		&		\\
	&		&	11	&	407	&		&		\\
	&		&	13	&	154	&		&		\\
\\											
Trade (0.3)	&	4	&	6 	&	641	&	735	&	5117	\\
	&		&	8 	&	308	&		&		\\
	&		&	11	&	402	&		&		\\
	&		&	13	&	154	&		&		\\
\\											
Trade (0.5)	&	4	&	7 	&	344	&	691	&	4140	\\
	&		&	8 	&	333	&		&		\\
	&		&	11	&	386	&		&		\\
	&		&	13	&	154	&		&		\\
\\											
Trade (0.7)	&	4	&	7 	&	351	&	691	&	3694	\\
	&		&	8 	&	219	&		&		\\
	&		&	11	&	362	&		&		\\
	&		&	13	&	154	&		&		\\
 \hline 
	\end{tabular} 
\label{table:OptimalLocations_Case2}
\end{table}

\subsubsection{\textcolor{black}{Analysis of optimal solutions quality: perfect vs. imperfect information}}

\begin{figure}[t!]
     \begin{subfigure}[b]{0.5\textwidth}
 \centering
        \includegraphics[width=\textwidth]{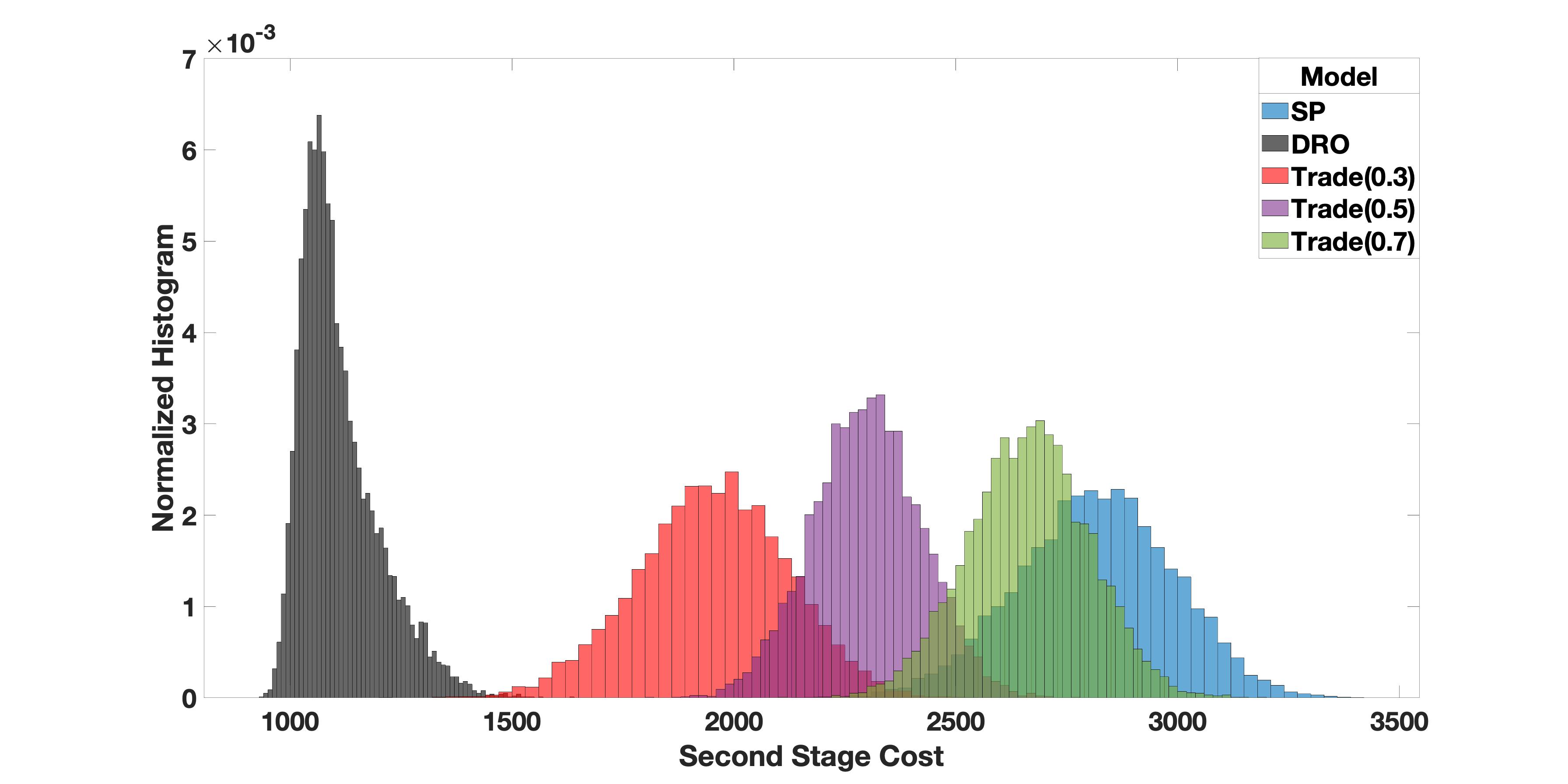}
        \caption{Second stage cost  }
        \label{Fig5a}
    \end{subfigure}%
    \begin{subfigure}[b]{0.5\textwidth}
            \includegraphics[width=\textwidth]{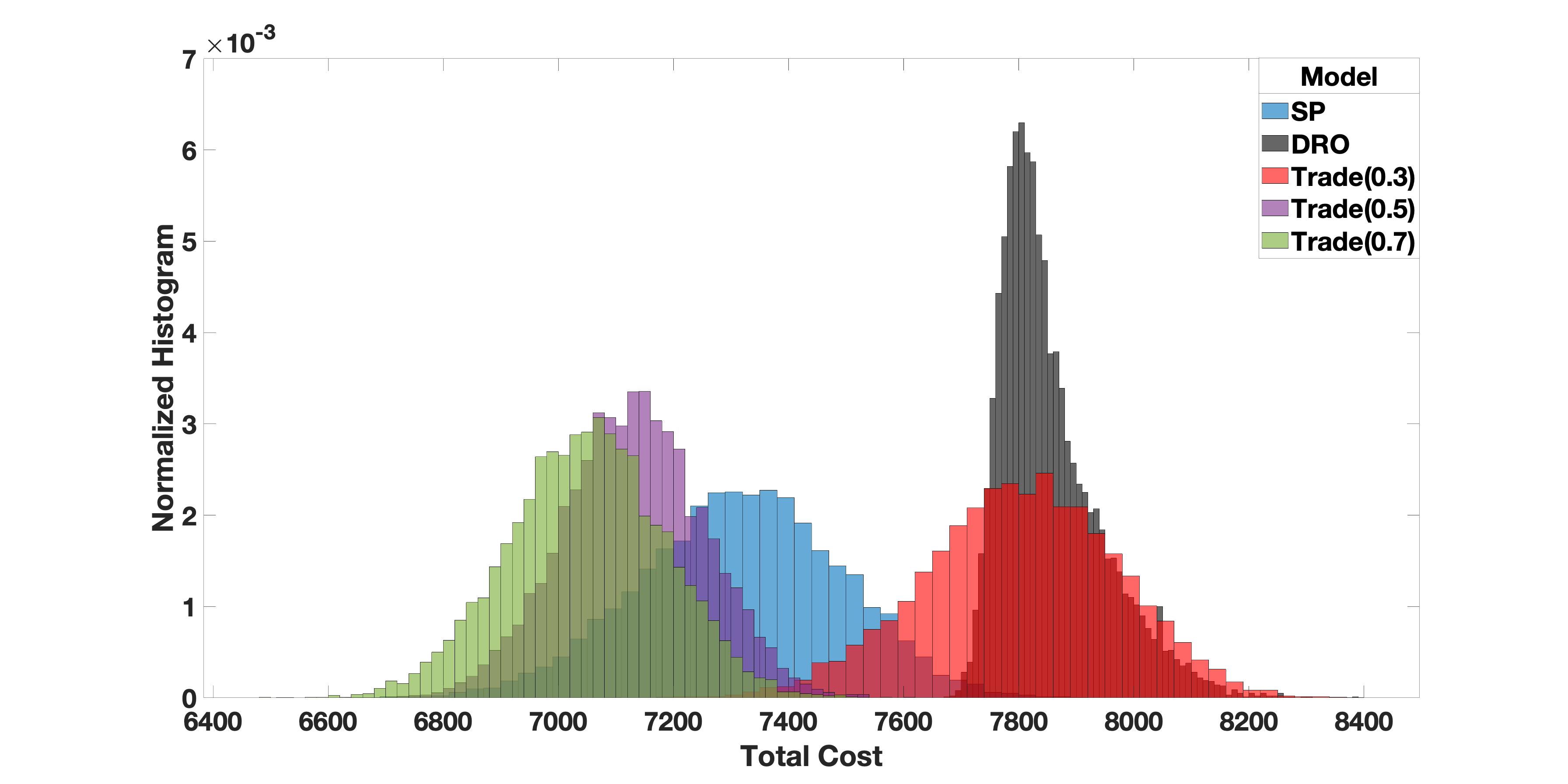}
      \caption{Total cost}
      \label{Fig5b}
    \end{subfigure}%
    \caption{Comparison of in-sample simulation results under perfect distributional information (LogN).}\label{Fig5:InSample}
\end{figure}

\noindent In this section, we use the same out-of-sample simulation procedure in Section~\ref{sec5:OutSample} to compare the operational performance of the optimal solutions to the DRO, SP, and Trade models presented in Table~\ref{table:OptimalLocations_Case2}.   Figure~\ref{Fig5:InSample} presents histograms of the total and second-stage costs under Set 1, i.e., when the decision-maker has perfect distributional information.  We first observe that the DRO and Trade (0.3) models yield larger total costs than the other considered models. This makes sense because these models have the highest prepositioning fixed cost. However, the DRO model has the lowest second-stage cost, followed by the Trade (0.3) model (Figures~\ref{Fig5a}). These models better hedge against uncertainty, which is realized in the second-stage.

\begin{figure}[t!]

 \begin{subfigure}[b]{0.5\textwidth}
 \centering
        \includegraphics[width=\textwidth]{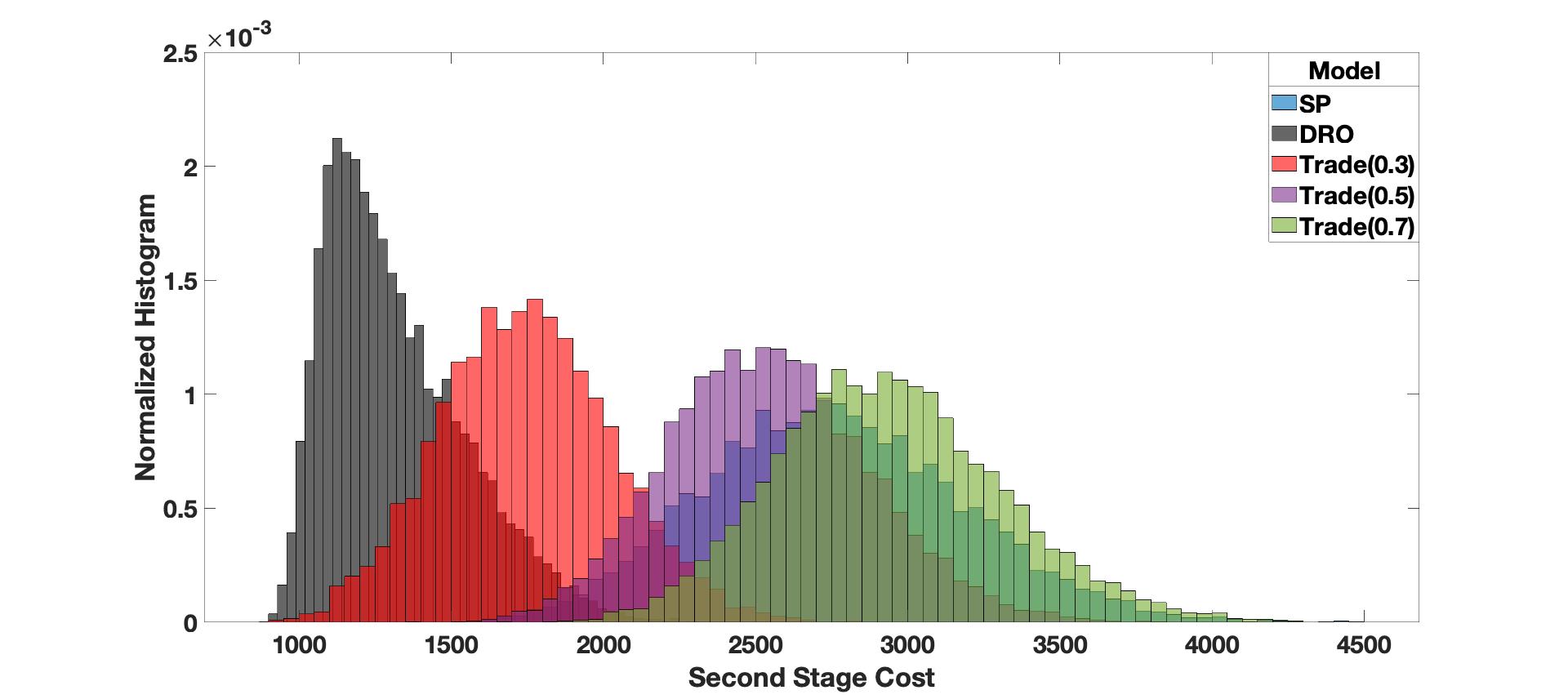}
        \caption{Second stage cost, $\Delta=0.25$}
        \label{Fig7a}
    \end{subfigure}%
    \begin{subfigure}[b]{0.5\textwidth}
            \includegraphics[width=\textwidth]{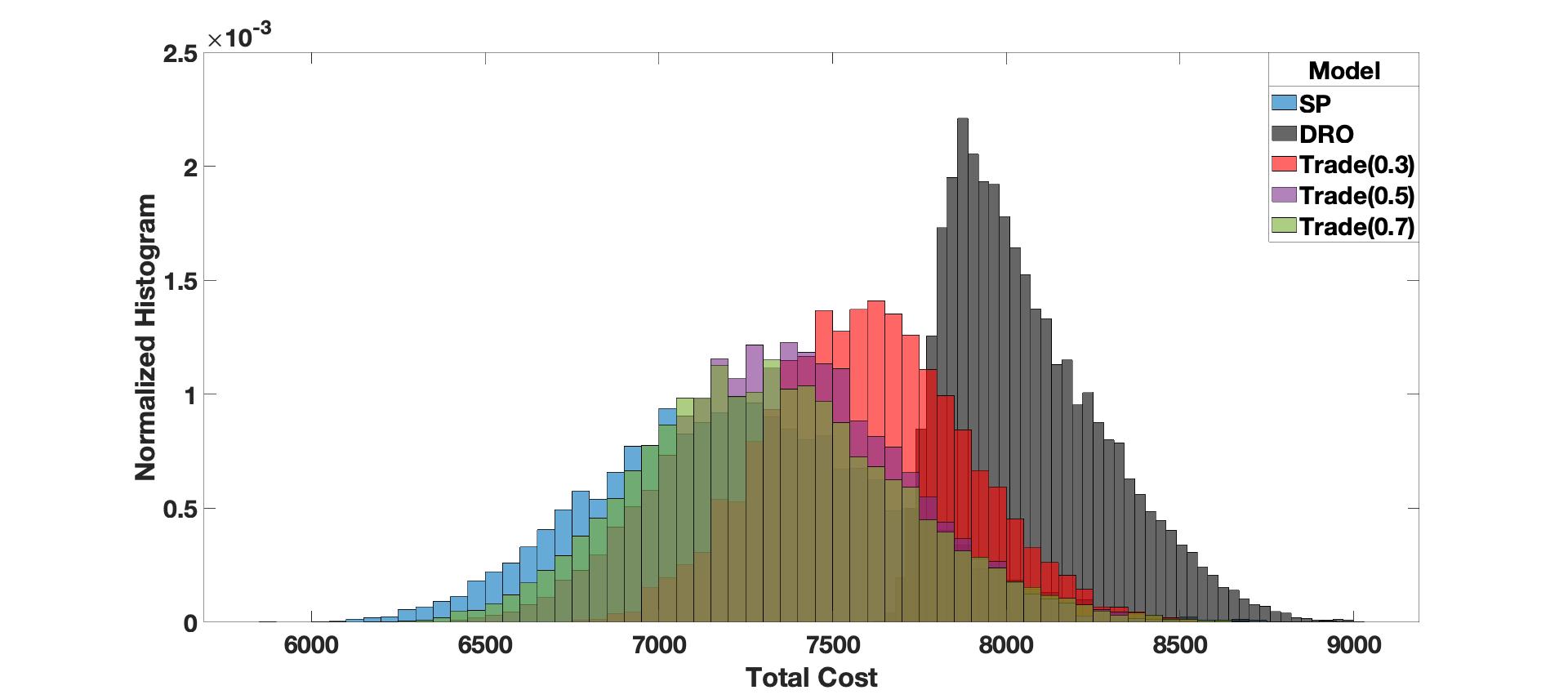}
      \caption{Total cost, $\Delta=0.25$}
      \label{Fig7b}
    \end{subfigure}%

     \begin{subfigure}[b]{0.5\textwidth}
 \centering
        \includegraphics[width=\textwidth]{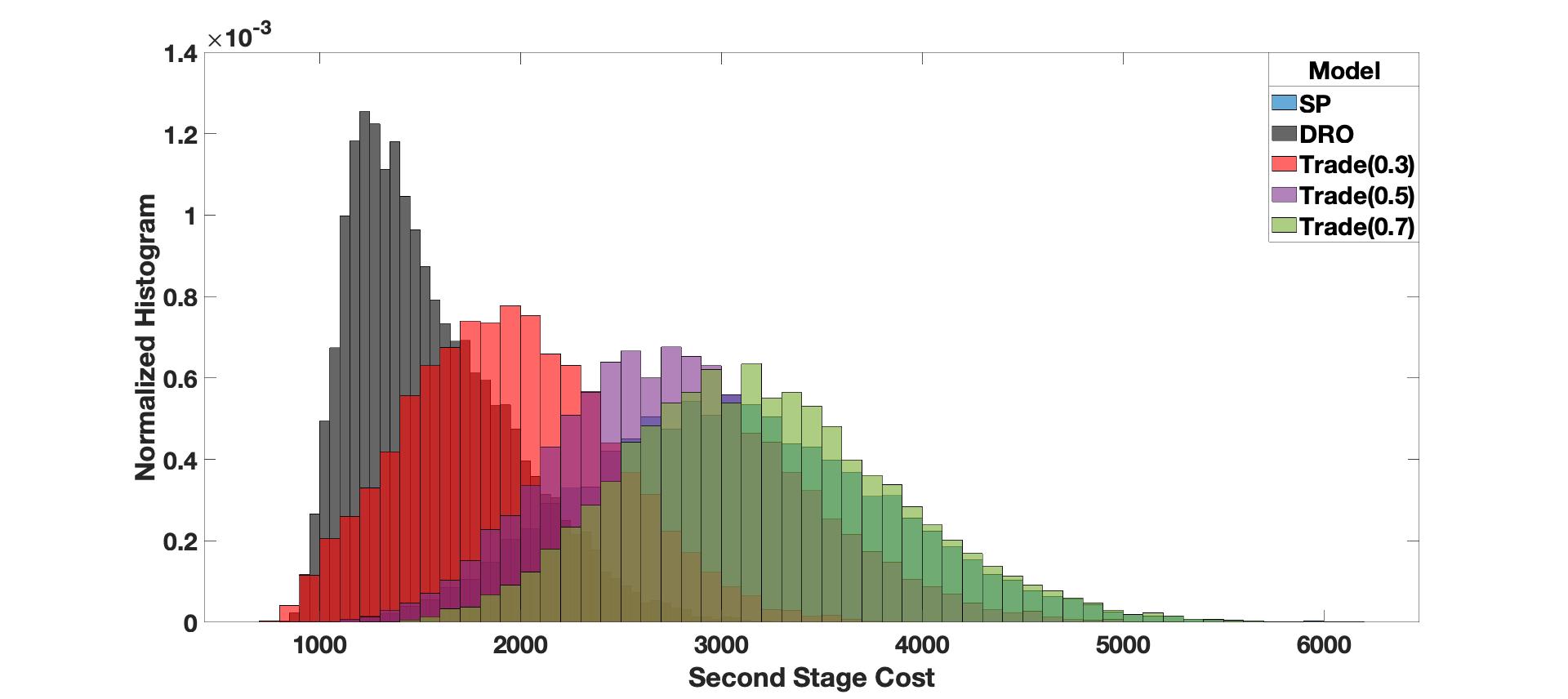}
        \caption{Second stage cost, $\Delta=0.5$  }
        \label{Fig8a}
    \end{subfigure}%
    \begin{subfigure}[b]{0.5\textwidth}
            \includegraphics[width=\textwidth]{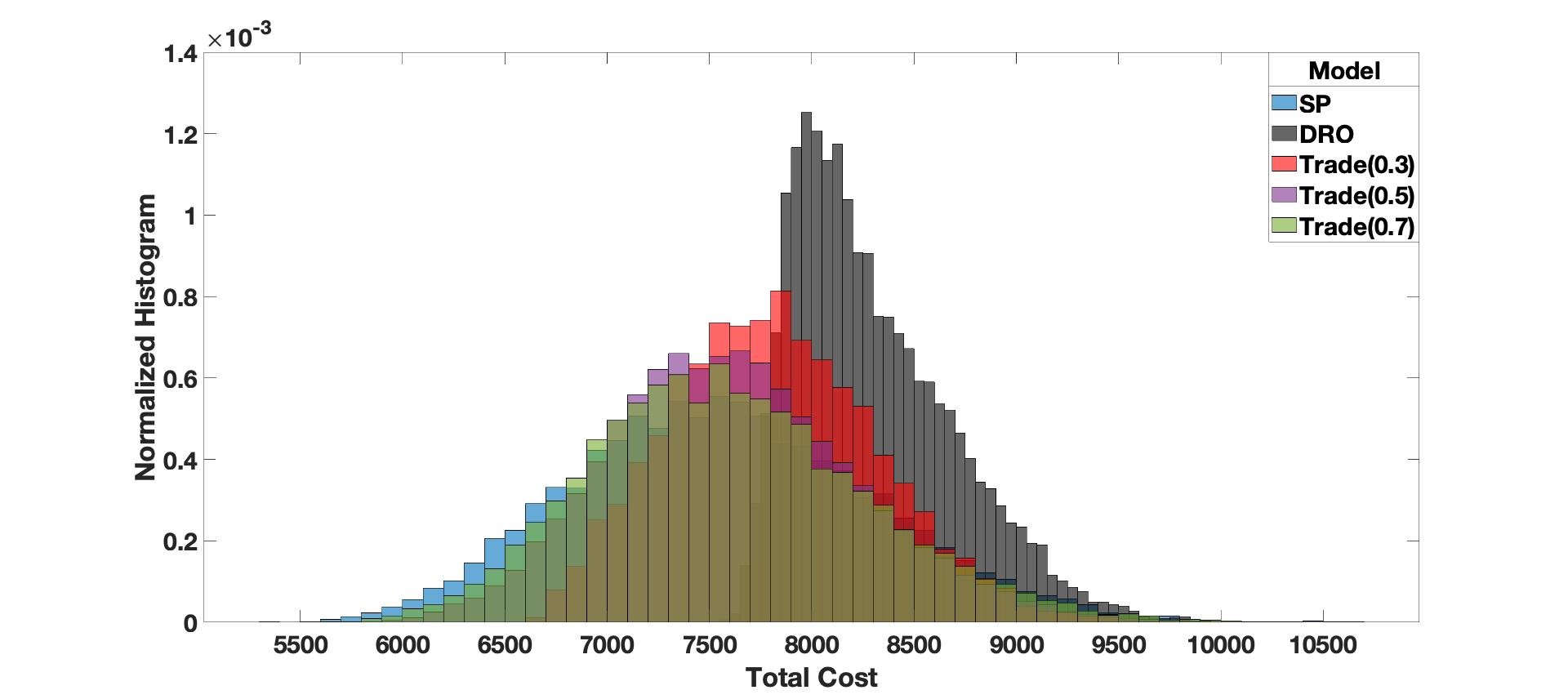}
      \caption{Total cost, $\Delta=0.5$}
      \label{Fig8 db}
    \end{subfigure}%
    \caption{Comparison of out-of-sample simulation results under Set 2 with $\Delta \in \{0.25, 0.5\}$}\label{Fig8:Out_delta50}
\end{figure}

Figure~\ref{Fig8:Out_delta50} presents the  out-of-sample results under Set 2  with $\Delta=0.25$ and $0.5$. These reflect the context where the true distribution is different than the one used in the optimization.  It is clear from this figure that the DRO and Trade (0.3) models maintain a robust performance with substantially lower second-stage costs (and thus better operational performance) than the other considered models. The significantly higher second-stage costs of the SP model indicates that solutions of this model have the worst operational performance in the aftermath.  As observed in Section~\ref{sec5:OutSample}, the Trade (0.5) model outperforms the Trade (0.7) model, and the latter model has approximately the same performance as the SP model.

 Consistent with the results in Section~\ref{sec5:OutSample}, the DRO and Trade (0.3) models show better stability by attaining the lowest variations in the second-stage and total costs across all scenarios. Although the DRO and Trade (0.3) models have higher prepositioning cost, the total cost yielded by all models becomes comparable as $\Delta$ increases (i.e., under high variability).  The DRO and Trade (0.3) models focus on hedging against uncertainty and distributional ambiguity. Their lower second-stage costs indicate that their solutions will have better operational performance in the aftermath, thus underscoring the value of incorporating uncertainty and ambiguity into inventory prepositioning models.

\subsection{Larger networks: computational performance of DRO-decomposition approach}\label{sec:performance}

\noindent In this section, we study the performance of the DRO-decomposition algorithm on larger, randomly generated networks. \textcolor{black}{Recall that the size of the case study problems were: hurricane season $G(\calN,\calA) = (30,112)$ and earthquake $G(\calN,\calA)=(13,15)$. To evaluate whether our DRO approach is tractable in other, larger contexts, we consider networks with $|I| \in \{40, 60, 100 \}$ nodes. We follow the procedure in \cite{ni2018location} to generate these networks.} Specifically, for each instance with $|I|$ nodes, we first generate the nodes on $10 \times 10 $ square, then we label them from 1 to $|I|$.  Second, we construct a spanning tree by connecting nodes $i$ and $j$  for any $i \in I\setminus \{1\}$ and some $j$ randomly selected from the set $\{1,\ldots, i-1\}$. Third, we randomly generate 0.2 $|I|+1$ pairs of nodes and add the corresponding undirected arcs to the network. Each network contains 1.2$|I|$ undirected arcs. We use the Euclidean metric to compute the distance between each pair of nodes $i$ and $j$ of arc $(i,j$).

We use the set of relief items from Case Study 1 (in Table~\ref{Table:Items_Costs}). Similar to Section~\ref{sec5:Hurrican}-\ref{sec:Earthquack}, the uncertain parameters ($\pmb{q, d, \rho, M, V}$) are characterized by their mean values and ranges. For these experiments, we set $\mu^{\mbox{\tiny d}}_{\mbox{\tiny water}}$, $\mu^{\mbox{\tiny d}}_{\mbox{\tiny food}}$, and $\mu^{\mbox{\tiny d}}_{\mbox{\tiny kits}}$ to their average values across all nodes in Table~\ref{Table:Demand}. We uniformly generate $ \mu^{\mbox{\tiny V}}_{i,j}$ from $U(20|I|, 25|I|)$ and $\sigma^{\mbox{\tiny V}}_{i,j}$ from $U(2|I|, 2.5|I|)$, as in \cite{ni2018location}. We generate $\mu^{\tiny \rho}$ and the range of each random parameters as described in Section~\ref{sec5:Hurrican}. For each network $|I| \in \{40, 60, 100 \}$, we conduct experiments with (nMinor, nMajor) $\in \{ (2,1), (4,2), (6,3) \}$ and disaster-prone nodes $\in$\{10, 15, 20\} for a total of 27 instances.

\black{Tables \ref{table:CPU_2_1}--\ref{table:CPU_6_3}  in  \ref{Appix:DROPerform} present the computational details of solving these instances using the DRO-decomposition algorithm}. Specifically, we report the number of iterations of Algorithm 1 before it converges to the optimum (\# of Iter), total CPU seconds taken by the master and the subproblems,  the total number of branching nodes (\# of B\&B), and the total number of MIP simplex iterations (\# of MIPiter). \blue{ From these results, we first observe that solution times and computational effort increase as the number of potential disasters increases.}

\color{black}
Second, we observe that for a fixed network size ($|I|$), the computational effort increases as the number of nodes vulnerable to disaster increases. For example, consider the instance with  (nMinor, nMajor)=$(4,2)$ and $I=60$. From 10 vulnerable nodes to 20, the (solution time; \# of B\&B; \# of MIPiter) increase from (25, 430, 179,054) to (133, 306,238, 2,090,9011). Third, we observe that the solution times of the subproblem are longer than the master problem. This makes sense because our subproblem is a MILP, and the size of this MILP increases as $|I|$, the number of disasters, and the number of nodes prone to disaster increases. As pointed out by \cite{artigues2015mixed, keha2009mixed,  klotz2013practical}, an increase in MILP size suggests an increase in solution time for the linear programming (LP) relaxation of the MILP and thus the solution time via commercial solvers.  These results show that our algorithm can solve large instances of the problem.

\color{black}

\section{Conclusion}\label{sec:conclusion}
\noindent In this paper, we present and analyze three new stochastic optimization models for location and inventory prepositioning of disaster relief supplies. We focus on the choice between different approaches to modeling uncertainty and how each approach yields a different prepositioning plan and performance.  Specifically, given a set of warehouse locations, a set of demand nodes, and a set of relief items,  the proposed models determine the number and locations of warehouses to open and the quantity of each relief item to preposition at each open location. In the aftermath, we consider the distribution of prepositioned relief items to demand locations and procurement and distribution additional supplies as needed. We consider the following random factors  (1) type of disaster, (2) locations of affected areas, (3) demand of relief items, (4) usable fraction of prepositioned items post-disaster, (5) procurement quantity, and (6) arc capacity between two different nodes. To model this uncertainty, we propose and analyze two stochastic optimization models--a two-stage SP and a two-stage DRO, assuming known and unknown distributions of uncertainty, respectively.  We also propose a model that minimizes the trade-off between considering distributional ambiguity and following distributional belief. We propose a decomposition algorithm to solve the DRO model and an MCO procedure to solve the SP model.

We conduct extensive experiments \textcolor{black}{using the three approaches (SP, DRO, and a trade-off between the two)} using a hurricane season and earthquake as case studies. \textcolor{black}{These illustrate very different types of disaster relief efforts. The applicability of the DRO and Trade models in these contexts suggest that considering distributional ambiguity in disaster relief could be worthwhile. We note that} \textcolor{black}{the DRO approach may} outperform the SP approach when there is limited distributional information. \textcolor{black}{We observe this clearly in the hurricane season where both total and second stage costs are lower using DRO and Trade approaches. The results for the earthquake case study are less decisive; the Trade approach leads to the lowest total cost, and SP outperforms DRO. When} the distribution is known with certainty (perfect information) that the SP approach performs better than DRO \textcolor{black}{in terms of pre-disaster fixed costs. However, the DRO solutions always offer a better post-disaster performance (i.e., operational cost)}.

\black{Regarding the effect on people, which we model via the shortage cost, DRO outperforms SP (i.e., DRO always satisfies greater demand in the immediate aftermath) even when there is perfect distributional information for the hurricane season. When the distribution is unknown, or the available estimates on post-disaster conditions are subject to error and uncertainty, the DRO and Trade models perform substantially better. In particular, the optimal DRO prepositioning decisions can satisfy nearly all demand, whereas the optimal SP prepositioning decisions lead to significant shortages. This suggests that the cost of misspecifying distributions \textcolor{black}{may} primarily \textcolor{black}{be} borne by the people the planners aim to serve. \textcolor{black}{Higher shortage costs in an SP model may be needed to mitigate this effect.} In addition, the DRO and Trade models yield lower procurement costs and second-stage costs. These results are consistent with theoretical results that indicate that SP decisions often over-promise (in our context, promise lower total cost and post-disaster operational cost) and under-deliver (in our context, SP solutions yield higher post-disaster operational cost than the estimated optimal costs) while the DRO solutions under-promise (in our context, have higher fixed prepositioning costs) and over-deliver (in our context, have lower  post-disaster operational cost).}

\textcolor{black}{Taken together, our results illustrate the (1)} the applicability of our approaches to multiple types of humanitarian logistics problems; (2) the robustness of DRO prepositioning decisions compared to SP decisions, especially under misspecified distributions and high variability; (3) the trade-off between considering distributional ambiguity (DRO) and following distributional belief (SP); and (4) the computational efficiency of our approaches. More broadly, our results draw attention to the need to model the distributional ambiguity of uncertain problem data in strategic real-world stochastic optimization problems such as planning for disasters. We encourage researchers, especially in uncertain context of humanitarian logistics, to consider the quality of the information they have when parameterizing their models \citep{comes2020coordination}.

We suggest the following areas for future research. First, we assumed that the random parameter distributions are unimodal. We would like to extend our approach by incorporating other random factors in a data-driven DRO approach. Partnerships with disaster relief agencies could provide data for new contexts and \textcolor{black}{further} improve the realism of the models. Second, we aim to include the possibility of restoring arcs functionality provided and the movement on restored arcs or during the restoration process. In particular, we \textcolor{black}{would like} to study how both the availability and restoration of damaged transportation networks impact preparedness and response decisions. Third, incorporating the cost of human suffering incurred by the shortage of relief supplies in the post-disaster relief operations \textcolor{black}{would be} another relevant and \textcolor{black}{useful} extension of our approach. \textcolor{black}{Forth, developing models for disaster response operations that incorporate the dynamic unfolding of the disaster and information would be practically relevant and theoretically interesting.}

\vspace{2mm}

\noindent \textbf{Acknowledgment}

\noindent We want to thank all colleagues who have contributed significantly to the related literature as well as the many disaster relief practitioners who tirelessly plan to reduce the impact of these events. \textcolor{black}{We thank the editor, associate editor, and the anonymous reviewers for their insightful comments and suggestions that allowed improving the paper}. Dr.~Karmel S.~Shehadeh dedicates her effort in this paper to every little dreamer in the whole world who has a dream so big and so exciting. Believe in your dreams and do whatever it takes to achieve them--the best is yet to come for you.

\bibliographystyle{elsarticle-harv}

\bibliography{Manuscript}

\appendix

\newpage 

\begin{center}
    {\Large Stochastic Optimization Models for Location and Inventory Prepositioning of Disaster Relief Supplies   (Appendices)}
\end{center}
\section{Random Factors}\label{Appex:RandomF}

\noindent \textcolor{black}{According to the recent surveys  of \cite{sabbaghtorkan2020prepositioning}, \cite{donmez2021humanitarian}, and our literature review in Section~\ref{sec2:Lit}, our paper is the first to compare the value and performance of SP and DRO approaches to address uncertainty and distributional ambiguity in parameters (1)--(6) in the specific location and inventory preposition problem that we study in this paper. Ignoring the collective variability of these factors may lead to devastating consequences. Ignoring uncertainty of the demand may lead to a significant shortage and thus human suffering in the immediate aftermath. Assuming that the prepositioned items will remain usable in the post-disaster is risky as a hurricane, for example, may hit the warehouses where these items are stored. Thus, ignoring this uncertainty will also lead to shortages and associated costs. Assuming that we can always procure additional relief items in the immediate aftermath is also risky. It may lead to significant shortages, and failure to meet the demand for relief items (i.e., shortage) as procuring additional items post-disaster is not easy and is often very expensive. Ignoring the uncertainty in arc capacity due to road damages, for example, may lead to failure in delivering relief items in the immediate aftermath, and consequently, shortage. Ignoring uncertainty of disaster location may lead to sub-optimal prepositioning decisions. Therefore, addressing these random factors' uncertainty (and potential distributional ambiguity) altogether is important and relevant in practice and theoretically interesting.}

\black{It is not clear when is best to use SP, DRO, or a trade-off model approach to model these factors. This is because different types of disasters have different dimensions of uncertainty and consequences. For example, for a predicted hurricane, we need models that can use some predictions (e.g., winds speed, intensity, and hurricane path, which NHC predicts before its landfall) while recognizing potential errors and variability on the estimates based on these predictions to suggest robust prepositioning decisions. While one may be able to forecast demand in the aftermath, uncertain factors such as the usable fraction of prepositioned items in the aftermath (a hurricane may hit the warehouses where these items are stored) and the ability to procure additional relief items are difficult to predict.  In contrast, there is no such predicted information before it strikes for a disaster like an earthquake, i.e., a lower information context). In addition, it is often difficult to predict which road networks may be partially functioning or completely damaged after an earthquake. We do not have such a situation in hurricanes.}

\newpage

\section{Proof of Proposition~\ref{Prop:InnerSup}}\label{Appex:Prop1}

\noindent \textit{Proof. }  For feasible first-stage decisions $(\ob,\zb)$, we can write the inner maximization problem $\sup \limits_{\Prob\in \calF(\calR, \mu)} \E_\Prob[Q(\ob,\zb,\xi)]$ in \eqref{DRModel} as the following linear functional optimization problem. 
\begin{subequations}
\begin{align}
& \max_{\textcolor{black}{\Prob \in \mathcal{P}(\mathcal{R})}} \ \int_{\textcolor{black}{\mathcal{R}}}Q(\ob,\zb,\xi)  \ d \Prob \\
& \ \text{s.t.} \ \   \int_{\textcolor{black}{\mathcal{R}}} d_{t,i,l}\ d\mathbb{P}= \mu_{t,i,l}^{\tiny d} && \forall  i \in I , \forall t \in \calT, \forall l \in L\label{Const1:InnerMax}\\
& \ \ \ \ \  \ \int_{\textcolor{black}{\mathcal{R}}} M_{t,i,l} \ d\mathbb{P}= \mu_{t,i,l}^{\tiny M} && \forall i \in I, \ \forall t \in \calT, \ \forall l \in L \label{Const6-2:InnerMax}\\
& \ \ \ \ \  \ \int_{\textcolor{black}{\mathcal{R}}} \rho_{i,l} \ d\mathbb{P}= \mu_{i,l}^{\tiny \rho}  && \forall  i \in I, \forall l \in L\label{Const3:InnerMax}\\
& \ \ \ \ \  \ \int_{\textcolor{black}{\mathcal{R}}} V_{i,j,l} \ d\mathbb{P}= \mu_{i,j,l}^{\tiny V} && \forall  (i,j \in \calA, \forall l \in L \label{Const5:InnerMax}\\
& \ \ \ \ \  \ \ \int_{\textcolor{black}{\mathcal{R}}} q_{i,l}\ d\mathbb{P}= \mu_{i,l}^{\tiny q} && \forall i \in I, \forall l \in L \label{Const7:InnerMax}\\
& \ \ \ \ \  \ \ \int_{\textcolor{black}{\mathcal{R}}}  d\mathbb{P}= 1. \label{Const8:InnerMaxDistribution}
\end{align} \label{InnerMax}
\end{subequations}
\noindent \textcolor{black}{Note that the optimization in \eqref{InnerMax} is performed over probability measure $\Prob$.} Letting $\alpha_{t,i,l}$, $\phi_{t,i,l}$, $\gamma_{i,l}$, $\tau_{i,j,l}$, $\lambda_{i,l}$, and $\theta$ be the dual variables associated with constraints \eqref{Const1:InnerMax}--\eqref{Const8:InnerMaxDistribution}, we present problem \eqref{InnerMax} in its dual form:
\begin{subequations}
\begin{align}
 \min_{\pmb{\alpha, \phi, \gamma, \lambda, \tau, \theta}}  &\Bigg \{\sum_{l \in L}\Big[ \sum \limits_{i \in I} \sum \limits_{t \in \calT} \Big(\mu_{t,i,l}^{\mbox{\tiny d}}  \alpha_{t,i,l} +   \mu_{t,i,l}^{\mbox{\tiny M}}\phi_{t,i,l} \Big)+ \sum \limits_{i \in I}\Big( \mu_{i,l}^{\tiny\rho}\gamma_{i,l}+\mu_{i,l}^{\tiny q} \lambda_{i,l} \Big)+ \sum \limits_{(i,j) \in \calA} \mu_{i,j,l}^{\tiny V}\tau_{i,j,l} \Big]+\theta \Bigg\} \label{DualInner:Obj} \\
\text{s.t.} & \sum_{l \in L} \Big[  \sum \limits_{i \in I} \sum \limits_{t \in \calT}d_{t,i,l}\alpha_{t,i,l}+M_{t,i,l} \phi_{t,i,l}+\sum \limits_{i \in I}\rho_{i,l}\gamma_{i,l}+q_{i,l}\lambda_{i,l}+\sum \limits_{(i,j) \in \calA}  V_{i,j,l} \tau_{i,j,l}\Big] \nonumber\\
& + \theta \geq Q(\ob,\zb, \xi), \ \ \forall  \xi \in \calR, \label{DualInner:PrimalVariabl}
\end{align} \label{DualInnerMax}
\end{subequations}
where $\alpha,  \ \phi, \ \gamma, \ \tau,  \ \lambda$, and $\theta$ are  are unrestricted in sign, and constraint \eqref{DualInner:PrimalVariabl} is associated with the primal variable $\mathbb{P}$. Note that strong duality hold between \eqref{InnerMax} and \eqref{DualInnerMax} \citep{bertsimas2005optimal, jian2017integer, shehadehDMFRS}. In addition, we observe that constraint \eqref{DualInner:PrimalVariabl} is equivalent to $\theta \geq \max \limits_{\xi \in \calR}  \big\{ Q(\ob,\zb, \xi) +\sum \limits_{l \in L} \big[ \sum \limits_{i \in I} \sum \limits_{t \in \calT}-(d_{t,i,l}\alpha_{t,i,l}+M_{t,i,l} \phi_{t,i,l})+\sum \limits_{i \in I}-(\rho_{i,l}\gamma_{i,l}+q_{i,l}\lambda_{i,l})+\sum \limits_{(i,j) \in \calA}  -V_{i,j,l} \tau_{i,j,l}\big] \big \}$. Since we are minimizing $\theta$ in \eqref{DualInnerMax}, the dual formulation of \eqref{InnerMax} is equivalent to \eqref{eq:FinalDualInnerMax-1}. This completes the proof. \qed

\newpage

\section{Proof of Proposition~\ref{Prop:MILP_Inner}}\label{Appex:McCIneq}

\color{black}
First, we re-write problem \eqref{Final_InnerQ} as follows
\begin{subequations}\label{Intermediate}
\begin{align} 
\max \limits_{\substack{\pmb{\beta, \Gamma,\psi, \varphi} \\ \pmb{q, d, \rho, M, V}}} &\Bigg\{ \sum \limits_{i \in I} \sum \limits_{t\in \calT}  \sum_{l \in L} d_{t,i,l} \big(q_{i,l}\beta_{t,i}-\alpha_{t,i,l} \big) + \sum_{i \in I} \sum_{t \in T}  z_{t,i}\big (\sum_{l \in L} q_{i,l}\beta_{t,i}- \beta_{t,i}\big) \nonumber \\
& \ \ \ + \sum \limits_{i \in I}   \sum_{l \in L} - \rho_{i,l} \big( \sum \limits_{t\in \calT}z_{t,i} q_{i,l}\beta_{t,i}+\gamma_{i,l} \big) \nonumber \\
&   \ \ \ + \sum \limits_{i \in I} \sum \limits_{t \in \calT} \sum_{l \in L} M_{t,i,l}\big(q_{i,l}\Gamma_{t,i}-\phi_{t,i,l} \big) \nonumber \\
& \ \ \ +  \sum_{(i,j)\in \calA}  \sum_{l \in L} V_{i,j,l} \big (q_{i,l} \psi_{i,j}+q_{j,l} \varphi_{i,j} -\tau_{i,j,l} \big) \nonumber \\
& \ \ \  + \sum_{(i,j)\in \calA}  \Big[ \big(\psi_{i,j}-\sum_{l \in L}q_{i,l} \psi_{i,j}\big)\hat{V}_{i,j} + \big(\varphi_{i,j}-\sum_{l \in L}q_{j,l} \varphi_{i,j}\big)\hat{V}_{i,j}  \Big] -\sum_{i \in I} \sum_{l \in L} q_{i,l} \lambda_{i,l} \Bigg\} \\
 \text{s.t.}  & \ \  \{ \eqref{Const1:DualofQ}-\eqref{Const4:DualofQ}\}, \pmb{q} \in \{0, 1\},  \\
 & \ \    \dL_{t,i,l} \leq d_{t,i,l} \leq \dU_{t,i,l},  \ \ \forall  t \in \calT, l \in I, l \in L ,\label{inter_C1} \\
 & \ \  \ML_{t,i,l} \leq M_{t,i,l} \leq \MU_{t,i,l}, \ \ \forall  t \in \calT, l \in I, l \in L, \label{inter_C2} \\
 & \ \  \rhoL_{i,l} \leq \rho_{i,l} \leq \rhoU_{i,l}, \ \ \forall i \in I, l \in L, \label{inter_C3} \\
 & \ \ \VL_{i,j,l} \leq V_{i,j,l} \leq \VU_{i,j,l}, \ \ \forall (i,j) \in \calA, l \in L.
\end{align}
\end{subequations}
It is easy to verify that in the optimal solution to problem \eqref{Intermediate}, constraint \eqref{inter_C1} is binding at either the lower or upper bound, i.e.,  $d_{t,i,l}=\dL_{t,i,l}$ or $d_{t,i,l}=\dU_{t,i,l}$. Therefore, we define binary variable $a_{t,i,l}$ that equals 1 if $d_{t,i,l}=\dU_{t,i,l}$, and is zero if $d_{t,i,l}=\dL_{t,i,l}$, and replace $d_{t,i,l}$ with $\dL_{t,i,l}+a_{t,i,l} \Delta d_{t,i,l}$, where $\Delta d_{t,i,l}=(\dU_{t,i,l}-\dL_{t,i,l})$. Applying the same logic, we can equivalently replace $\rho_{i,l}$ with $\rhoL_{i,l}+\Delta\rho_{i,l}\Theta_{i,l}$, $M_{t,i,l}$ with $\ML_{t,i,l}+\Delta M_{t,i,l}\kappa_{t,i,l}$, $V_{i,j,l}$ with $\VL_{i,j,l}+\Delta V_{i,j,l}\varrho_{i,j,l}$, where $\Delta \rho_{i,l}=(\rhoU_{i,l}-\rhoL_{i,l})$, $\Delta M_{t,i,l}=(\MU_{t,i,l}-\ML_{t,i,l})$, $\Delta V_{i,j,l}=(\VU_{i,j,l}-\VL_{i,j,l})$, and $(\pmb{\Theta, \kappa,  \varrho}) \in \{0,1\}$. Accordingly, we derive the following equivalent reformulation of problem \eqref{Intermediate}.

\begin{subequations}\label{Intermediate2}
\begin{align} 
\max \limits_{\substack{\pmb{\beta, \Gamma,\psi, \varphi}\\ \pmb{a,\Theta, \kappa, \varrho}}} &\Bigg\{ \sum \limits_{i \in I} \sum \limits_{t\in \calT} \sum_{l \in L}  \Big[\dL_{t,i,l} \big(q_{i,l}\beta_{t,i}-\alpha_{t,i,l} \big) + \Delta d_{t,i,l}\big(q_{i,l}\beta_{t,i}a_{t,i,l}-\alpha_{t,i,l}a_{t,i,l} \big) \Big] \nonumber\\
& \ \ \ +\sum_{i \in I} \sum_{t \in T}  z_{t,i}\big (\sum_{l \in L} q_{i,l}\beta_{t,i}- \beta_{t,i}\big) \nonumber \\
& \ \ \ - \sum \limits_{i \in I}   \sum_{l \in L} \Big[  \rhoL_{i,l} \big( \sum \limits_{t\in \calT}z_{t,i} q_{i,l}\beta_{t,i}+\gamma_{i,l} \big)+ \Delta \rho_{i,l} \big( \sum \limits_{t\in \calT}z_{t,i} q_{i,l} \Theta_{i,l}\beta_{t,i}+\Theta_{i,l}\gamma_{i,l} \big) \Big] \nonumber \\
&   \ \ \ + \sum \limits_{i \in I} \sum \limits_{t \in \calT} \sum_{l \in L} \Big[ \ML_{t,i,l}\big(q_{i,l}\Gamma_{t,i}-\phi_{t,i,l} \big)+\Delta M_{t,i,l} \big(q_{i,l}\Gamma_{t,i}\kappa_{t,i,l}-\phi_{t,i,l}\kappa_{t,i,l} \big) \Big]  \nonumber \\
& \ \ \ +  \sum_{(i,j)\in \calA}  \sum_{l \in L} \Big[ \VL_{i,j,l} \big (q_{i,l} \psi_{i,j}+q_{j,l} \varphi_{i,j} -\tau_{i,j,l} \big) + \Delta V_{i,j,l} \varrho_{i,j,l} \big (q_{i,l} \psi_{i,j}+q_{j,l} \varphi_{i,j} -\tau_{i,j,l} \big)\nonumber \\
& \ \ \  + \sum_{(i,j)\in \calA}  \Big[ \big(\psi_{i,j}-\sum_{l \in L}q_{i,l} \psi_{i,j}\big)\hat{V}_{i,j} + \big(\varphi_{i,j}-\sum_{l \in L}q_{j,l} \varphi_{i,j}\big)\hat{V}_{i,j}  \Big] -\sum_{i \in I} \sum_{l \in L} q_{i,l} \lambda_{i,l} \Bigg\}  \label{Obj_Inter2} \\
 \text{s.t.}  & \ \  \{ \eqref{Const1:DualofQ}-\eqref{Const4:DualofQ}\}, \pmb{q} \in \{0, 1\},  (\pmb{a,\Theta, \kappa, \varrho}) \in \{0, 1\}.
\end{align}
\end{subequations}

 Note that objective function \eqref{Obj_Inter2} contains the interactions terms $q_{i,l}\beta_{t,i}$, $a_{t,i,l}q_{i,l}\beta_{t,i}$, $\Theta_{i,l}q_{i,l}\beta_{t,i}$, $q_{i,l}\Gamma_{t,i}$, $\kappa_{t,i,l}q_{i,l}\Gamma_{t,i}$, $q_{i,l}\psi_{i,j}$, $\varrho_{i,j,l} q_{i,l}\psi_{i,j}$, $q_{j,l} \varphi_{i,j}$, and $\varrho_{i,j,l} q_{j,l'}\varphi_{i,j}$, which consist of binary variables multiplied by continuous variables. To linearize, we define variables $k_{t,i,l}=q_{i,l} \beta_{t,i}$, $a'_{t,i,l}=a_{t,i,l}q_{i,l}$, $h_{t,i,l}=a'_{t,i,l}\beta_{t,i}$, $\Theta'_{i,l}=\Theta_{i,l}q_{i,l}$, $g_{t,i,l}=\Theta'_{i,l}\beta_{t,i}$, $F_{t,i,l}=q_{i,l}\Gamma_{t,i}$, $\kappa'_{t,i,l}=\kappa_{t,i,l}q_{i,l}$, $\pi_{t,i,l}=\kappa'_{t,i,l}\Gamma_{t,i}$, $\eta_{i,j,l}=q_{i,l}\psi_{i,j}$, $\varrho'_{i,j,l}=\varrho_{i,j,l}q_{i,l}$, $ \Phi_{i,j,l}=\varrho'_{i,j,l} \psi_{i,j}$,  $\varpi_{i,j,l}=q_{j,l}\varphi_{i,j}$, $ b_{i,j,l}=\varrho_{i,j,l}q_{j,l}$, and $\Lambda_{i,j,l}=b_{i,j,l} \varphi_{i,j}$.  Also, we introduce the following McCormick inequalities for these variables
\begin{subequations}\label{MacCormick}
\begin{align}
&k_{t,i,l} \geq q_{i,l} \betaL_{t,i}, \ \  k_{t,i,l}\geq \beta_{t,i}+\betaU_{t,i}(q_{i,l}-1), \ k_{t,i,l}\leq q_{i,l}\betaU_{t,i}, \  k_{t,i,l}\leq \beta_{t,i}+ \betaL_{t,i}(q_{i,l}-1) \label{Mac1}\\
& a'_{t,i,l} \geq 0, \ \  a'_{t,i,l} \leq q_{i,l}, \ \  a'_{t,i,l} \leq a_{t,i,l} , \ \ a'_{t,i,l} \geq a_{t,i,l} +q_{i,l}-1 \label{Mac2}\\
& h_{t,i,l} \geq a'_{t,i,l} \betaL_{t,i}, \ \  h_{t,i,l}\geq \beta_{t,i,}+\betaU_{t,i}(a'_{t,i,l}-1), \ h_{t,i,l}\leq a'_{t,i,l}\betaU_{t,i}, \  h_{t,i,l}\leq \beta_{t,i}+ \betaL_{t,i}(a'_{t,i,l}-1) \label{Mac3}\\
& \Theta'_{i,l}\geq 0, \ \  \Theta'_{i,l}\leq q_{i,l}, \ \  \Theta'_{i,l}\leq \Theta_{i,l}, \ \ \Theta'_{i,l}\geq  \Theta_{i,l}+q_{i,l}-1 \label{Mac4}\\
& g_{t,i,l} \geq \Theta'_{i,l} \betaL_{t,i}, \ \  g_{t,i,l}\geq \beta_{t,i,}+\betaU_{t,i}(\Theta'_{i,l}-1), \ \ g_{t,i,l}\leq \Theta'_{i,l}\betaU_{t,i}, \ \  g_{t,i,l}\leq \beta_{t,i}+ \betaL_{t,i}(\Theta'_{i,l}-1), \label{Mac5}\\
& F_{t,i,l}\geq q_{i,l} \GammaL_{t,i}, \ \ F_{t,i,l}\geq \Gamma_{t,i}, \ \ F_{t,i,l}\leq 0, \ \ F_{t,i,l}\leq \Gamma_{t,i}+ \GammaL_{t,i}(q_{i,l}-1), \label{Mac6}\\
& \kappa'_{t,i,l}\geq 0, \ \  \kappa'_{t,i,l}\leq q_{i,l}, \ \  \kappa'_{t,i,l}\leq \kappa_{t,i,l}, \ \ \kappa'_{t,i,l}\geq \kappa_{t,i,l}+q_{i,l}-1 \label{Mac7}\\
& \pi_{t,i,l}\geq \kappa'_{t,i,l}\GammaL_{t,i}, \ \ \pi_{t,i,l}\geq \Gamma_{t,i}, \ \ \pi_{t,i,l}\leq 0, \ \ \pi_{t,i,l}\leq \Gamma_{t,i}+ \GammaL_{t,i}(\kappa'_{t,i,l}-1), \label{Mac8}\\
 & \eta_{i,j,l} \geq q_{i,l} \underline{\psi}_{i,j}, \ \eta_{i,j,l} \geq \psi_{i,j}, \ \eta_{t,i,l} \leq 0, \ \eta_{t,i,l} \leq \psi_{i,j}+\underline{\psi}_{i,j}(q_{i,l}-1)  \label{Mac9}\\
 & \varrho'_{i,j,l}\geq 0, \ \  \varrho'_{i,j,l}\leq q_{i,l}, \ \  \varrho'_{i,j,l}\leq \varrho_{i,j,l} , \ \ \varrho'_{i,j,l}\geq  \varrho_{i,j,l} +q_{i,l}-1 \label{Mac10}\\
& \Phi_{i,j,l} \geq \varrho'_{i,j,l} \underline{\psi}_{i,j},   \ \Phi_{i,j,l} \geq \psi_{i,j}, \ \Phi_{i,j,l}\leq 0, \ \ \Phi_{i,j,l}\leq \psi_{i,j}+ \underline{\psi}_{i,j}(\varrho'_{i,j,l}-1), \label{Mac11}\\
& \varpi_{i,j,l} \geq q_{j,l} \underline{\varphi}_{i,j}, \ \varpi_{i,j,l} \geq \varphi_{i,j}, \  \varpi_{i,j,l}  \leq 0, \  \varpi_{i,j,l}  \leq \varphi_{i,j}+\underline{\varphi}_{i,j}(q_{j,l}-1)  \label{Mac12}\\
& b_{i,j,l} \geq 0, \ \  b_{i,j,l} \leq q_{j,l}, \ \  b_{i,j,l} \leq \varrho_{i,j,l} , \ \ b_{i,j,l} \geq  \varrho_{i,j,l} +q_{j,l}-1 \label{Mac13}\\
& \Lambda_{i,j,l}\geq b_{i,j,l} \underline{\varphi}_{i,j},   \ \Lambda_{i,j,l}\geq \varphi_{i,j}, \Lambda_{i,j,l}\leq 0, \ \Lambda_{i,j,l}\leq \varphi_{i,j}+ \underline{\varphi}_{i,j}(b_{i,j,l}-1), \label{Mac14}
\end{align}
\end{subequations}

Using variables ($\pmb{k,  a', h, \Theta',  g, F,  \kappa', \pi, \eta, \varrho', \varpi, \Phi, b, \Lambda} $) and inequalities \eqref{Mac1}--\eqref{Mac14}, we derive the following equivalent MILP reformulation of problem \eqref{Intermediate2}.

\begin{subequations}\label{Intermediate3}
\begin{align} 
\max \limits_{\substack{\pmb{\beta, \Gamma,\psi, \varphi, a,\Theta, \kappa, \varrho} \\ \pmb{a', \Theta', \kappa', \varrho', h, g, F, \pi} \\ \pmb{\eta, \Phi, \varpi, b,\Lambda}}}  &\Bigg\{ \sum \limits_{i \in I} \sum \limits_{t\in \calT} \sum_{l \in L}  \Big[\dL_{t,i,l} \big(k_{t,i,l}-\alpha_{t,i,l} \big) + \Delta d_{t,i,l}\big(h_{t,i,l}-\alpha_{t,i,l}a_{t,i,l} \big) \Big] \nonumber\\
& \ \ \ +\sum_{i \in I} \sum_{t \in T}  z_{t,i}\big (\sum_{l \in L} k_{t,i,l}- \beta_{t,i}\big) \nonumber \\
& \ \ \ - \sum \limits_{i \in I}   \sum_{l \in L} \Big[  \rhoL_{i,l} \big( \sum \limits_{t\in \calT}z_{t,i} k_{t,i,l}+\gamma_{i,l} \big)+ \Delta \rho_{i,l} \big( \sum \limits_{t\in \calT}z_{t,i} g_{t,i,l}+\Theta_{i,l}\gamma_{i,l} \big) \Big] \nonumber \\
&   \ \ \ + \sum \limits_{i \in I} \sum \limits_{t \in \calT} \sum_{l \in L} \Big[ \ML_{t,i,l}\big(F_{t,i,l}-\phi_{t,i,l} \big)+\Delta M_{t,i,l} \big(\pi_{t,i,l}-\phi_{t,i,l}\kappa_{t,i,l} \big) \Big]  \nonumber \\
& \ \ \ +  \sum_{(i,j)\in \calA}  \sum_{l \in L} \Big[ \VL_{i,j,l} \big (\eta_{i,j,l}+\varpi_{i,j,l} -\tau_{i,j,l} \big) + \Delta V_{i,j,l} \big (\Phi_{i,j,l}+\Lambda_{i,j,l} -\varrho_{i,j,l} \tau_{i,j,l} \big)\nonumber \\
& \ \ \  + \sum_{(i,j)\in \calA}  \Big[ \big(\psi_{i,j}-\sum_{l \in L}\eta_{i,j,l}\big)\hat{V}_{i,j} + \big(\varphi_{i,j}-\sum_{l \in L}\varpi_{i,j,l}\big)\hat{V}_{i,j}  \Big] -\sum_{i \in I} \sum_{l \in L} q_{i,l} \lambda_{i,l} \Bigg\}  \label{Obj_Inter2} \\
 \text{s.t.}  & \ \  \{ \eqref{Const1:DualofQ}-\eqref{Const4:DualofQ}\}, \pmb{q} \in \{0, 1\},  (\pmb{a,\Theta, \kappa, \varrho}) \in \{0, 1\}, \eqref{Mac1}-\eqref{Mac14}.
\end{align}
\end{subequations}

This completes the proof. \qed

Note that the McCormick inequalities often rely on big-M coefficients ($\overline{\cdot}$ and $\underline{\cdot}$) that take large values and thus may undermine computational efficiency. In \ref{Appix:BigM}, we derive tight bounds of these big-M coefficients to strengthen the MILP formulation.
\color{black}

\section{Strengthening the MILP Formulation}\label{Appix:BigM}

\noindent First, observe from constraint \eqref{Const3:DualofQ} that $-\ch \leq \beta_{t,i} \leq \cu$. Thus, w.l.o.o, we can assume that $\betaL_{t,i}=-\ch$ and $\betaU=\cu$. Second, given that $k_{t,i,l}=q_{i,l}\beta_{t,i}$, then w.l.o.o., $k_{t,i,l} \in [ \betaL_{t,i}, \betaU_{t,i}]. $  Third, observe from constraints \eqref{Const2:DualofQ} and \eqref{Const4:DualofQ} that $\Gamma_{t,i}  \leq c_{t}^{\mbox{\tiny p}}-\beta_{t,i}$ and $\Gamma_{t,i} \leq 0$. Note that in the optimal solution $\Gamma_{t,i}  =\min \{ c_{t}^{\mbox{\tiny p}}-\beta_{t,i}, 0 \} $. Therefore, $\GammaU_{t,i}=0$. And $\GammaL_{t,i} =  c_{t}^{\mbox{\tiny p}} -\cu$  if $ c_{t}^{\mbox{\tiny p}} \leq \cu$, and $\GammaL_{t,i} =0$ otherwise. 

Finally, from constraints \eqref{Const1:DualofQ}, we have $v_t\psi_{i,j} +v_t \varphi_{i,j} \leq c_{i,j}^t+\beta_{t,i}-\beta_{t,j}$. Given that $\psi_{i,j} \leq 0$ and $\varphi_{i,i} \leq 0$, then w.l.o.o., $\underline{\psi_{i,j}} \leq  \min_t \{ (1/v_t) (c_{i,j}^t  - \ch- \cu  \}$ if  $ c_{i,j}^t  \leq \ch+ \cu $  and $\underline{\psi_{i,j}}=0$ otherwise. 

\color{black}
\section{Proof of Proposition~\ref{Prop:Convex}}\label{Appix:Prop_Covex}

\color{black}
\noindent \textit{Proof.} First, it is easy to verify that for any fixed values of variables $\pmb{z, \  \alpha, \ \phi, \ \gamma, \ \lambda}$, and $\pmb{\tau}$

\begin{align*}
 & H (\pmb{z, \alpha, \phi, \gamma, \lambda, \tau}):=\\
\max \limits_{\substack{\pmb{\beta, \Gamma,\psi, \varphi, a,\Theta, \kappa, \varrho} \\ \pmb{a', \Theta', \kappa', \varrho', h, g, F, \pi} \\ \pmb{\eta, \Phi, \varpi, b,\Lambda}}} &\Bigg\{ \sum \limits_{i \in I} \sum \limits_{t\in \calT} \sum_{l \in L}  \Big[\dL_{t,i,l} \big(k_{t,i,l}-\alpha_{t,i,l} \big) + \Delta d_{t,i,l}\big(h_{t,i,l}-\alpha_{t,i,l}a_{t,i,l} \big) \Big] \nonumber\\
& \ \ \ +\sum_{i \in I} \sum_{t \in T}  z_{t,i}\big (\sum_{l \in L} k_{t,i,l}- \beta_{t,i}\big) \nonumber \\
& \ \ \ - \sum \limits_{i \in I}   \sum_{l \in L} \Big[  \rhoL_{i,l} \big( \sum \limits_{t\in \calT}z_{t,i} k_{t,i,l}+\gamma_{i,l} \big)+ \Delta \rho_{i,l} \big( \sum \limits_{t\in \calT}z_{t,i} g_{t,i,l}+\Theta_{i,l}\gamma_{i,l} \big) \Big] \nonumber \\
&   \ \ \ + \sum \limits_{i \in I} \sum \limits_{t \in \calT} \sum_{l \in L} \Big[ \ML_{t,i,l}\big(F_{t,i,l}-\phi_{t,i,l} \big)+\Delta M_{t,i,l} \big(\pi_{t,i,l}-\phi_{t,i,l}\kappa_{t,i,l} \big) \Big]  \nonumber \\
& \ \ \ +  \sum_{(i,j)\in \calA}  \sum_{l \in L} \Big[ \VL_{i,j,l} \big (\eta_{i,j,l}+\varpi_{i,j,l} -\tau_{i,j,l} \big) + \Delta V_{i,j,l} \big (\Phi_{i,j,l}+\Lambda_{i,j,l} -\varrho_{i,j,l} \tau_{i,j,l} \big)\nonumber \\
& \ \ \  + \sum_{(i,j)\in \calA}  \Big[ \big(\psi_{i,j}-\sum_{l \in L}\eta_{i,j,l}\big)\hat{V}_{i,j} + \big(\varphi_{i,j}-\sum_{l \in L}\varpi_{i,j,l}\big)\hat{V}_{i,j}  \Big] -\sum_{i \in I} \sum_{l \in L} q_{i,l} \lambda_{i,l} \Bigg\}  <\infty
\end{align*}

\noindent Second, for any fixed and feasible $\pmb{k, a,  h, \Theta,  g, F,  \kappa, \pi, \eta, \varrho, \varpi, \Phi, b, \Lambda}$, function

\begin{align*}
&\Bigg\{ \sum \limits_{i \in I} \sum \limits_{t\in \calT} \sum_{l \in L}  \Big[\dL_{t,i,l} \big(k_{t,i,l}-\alpha_{t,i,l} \big) + \Delta d_{t,i,l}\big(h_{t,i,l}-\alpha_{t,i,l}a_{t,i,l} \big) \Big] \nonumber\\
& \ \ \ +\sum_{i \in I} \sum_{t \in T}  z_{t,i}\big (\sum_{l \in L} k_{t,i,l}- \beta_{t,i}\big) \nonumber \\
& \ \ \ - \sum \limits_{i \in I}   \sum_{l \in L} \Big[  \rhoL_{i,l} \big( \sum \limits_{t\in \calT}z_{t,i} k_{t,i,l}+\gamma_{i,l} \big)+ \Delta \rho_{i,l} \big( \sum \limits_{t\in \calT}z_{t,i} g_{t,i,l}+\Theta_{i,l}\gamma_{i,l} \big) \Big] \nonumber \\
&   \ \ \ + \sum \limits_{i \in I} \sum \limits_{t \in \calT} \sum_{l \in L} \Big[ \ML_{t,i,l}\big(F_{t,i,l}-\phi_{t,i,l} \big)+\Delta M_{t,i,l} \big(\pi_{t,i,l}-\phi_{t,i,l}\kappa_{t,i,l} \big) \Big]  \nonumber \\
& \ \ \ +  \sum_{(i,j)\in \calA}  \sum_{l \in L} \Big[ \VL_{i,j,l} \big (\eta_{i,j,l}+\varpi_{i,j,l} -\tau_{i,j,l} \big) + \Delta V_{i,j,l} \big (\Phi_{i,j,l}+\Lambda_{i,j,l} -\varrho_{i,j,l} \tau_{i,j,l} \big)\nonumber \\
& \ \ \  + \sum_{(i,j)\in \calA}  \Big[ \big(\psi_{i,j}-\sum_{l \in L}\eta_{i,j,l}\big)\hat{V}_{i,j} + \big(\varphi_{i,j}-\sum_{l \in L}\varpi_{i,j,l}\big)\hat{V}_{i,j}  \Big] -\sum_{i \in I} \sum_{l \in L} q_{i,l} \lambda_{i,l} \Bigg\} 
\end{align*}
\noindent is a linear function of $\pmb{z, \ \ \alpha, \ \phi, \ \gamma, \ \lambda}$, and $\pmb{\tau}$. It follows that  $H (\pmb{z, \alpha, \phi, \gamma, \lambda, \tau})$  is the maximum of linear functions of $\pmb{z, \ \ \alpha, \ \phi, \ \gamma, \ \lambda}$, and $\pmb{\tau}$, and hence convex and peicewise linear. Finally, it is easy to verify the number of pieces of this function is finite.  This completes the proof.

\qed

\newpage
\section{The Monte Carlo Optimization (MCO) Procedure}\label{Appex:MCOProcedure}

 \begin{algorithm}[h!] 
 \LinesNumberedHidden
\footnotesize
  \renewcommand{\arraystretch}{0.4}
    \caption{The MCO Procedure}
      \label{alg:algorithm1}
\KwIn{$N_{o}$ is an initial sample size, $K$ is number of replicates, $N^\prime$ is number of scenarios in the Monte Carlo Simulation step,  and $\epsilon$ is a termination tolerance.}
\KwOut{$N$ is sample size, $\bar{v}_{N}$ and $\bar{v}_{N^\prime}$ are  respectively statistical lower and upper bounds on the optimal value of the SP, and $AOI_{N}$ is approximate optimality index.} 
\textbf{Initialization}: $N:=N_{o}$
\vspace{2mm}

\textbf{Step 1.} \textbf{MCO Procedure }

\For{$k=1,\ldots,K,$}{

 $\qquad  $ \textbf{Step 1.1 }\textit{\textbf{Scenario Generation}} 

 $\qquad \qquad \ \ \ $ - Generate $N$ independent and identical distributed (i.i.d.) scenarios of $(\pmb{d, \rho, M, V}$ \

$\qquad$  \textbf{Step 1.2 } \textit{\textbf{Solving the SAA formulation}}

$\qquad \qquad \ \ \ \  $ -  Solve the SAA formulation  in  \eqref{SAA} with the scenarios generated in step 1.1 and record
 
$\qquad \qquad \ \ \ \ \ \  $  the  corresponding optimal objective value $v_{N}^{k}$ and optimal solution $(\pmb{\hat{o}}, \pmb{\hat{z}})_{N}^{k}$.

$\qquad $  \textbf{Step 1.3 } \textit{\textbf{Cost Evaluation using Monte Carlo Simulation}} 

 
 
$\qquad \qquad \ \ \ \  $ - Generate a new $N^\prime$ i.i.d scenarios of $(\pmb{d', \rho', M', V'})$

$\qquad \qquad \ \ \ \  \ $-  Use solution solution $(\pmb{\hat{o}}, \pmb{\hat{z}})_{N}^{k}$  and parameters $(\pmb{d', \rho', M', V'})$  to compute 
$x', u', e', y' $, and 

$\qquad \qquad \ \ \ \ \ \  $evaluate the objective function $\hat{v}_{N^\prime}^{k}$ as follows:

  $$\hat{v}_{N^\prime}^{k}=   \sum_{i \in I} f_i \hat{o}_i + \sum \limits_{i \in I} \sum\limits_{t \in \calT} \ca \hat{z}_{t,i} + \sum_{n=1}^{N'} \frac{1}{N'} \Big[ \sum \limits_{i \in I} \sum \limits_{t \in \calT} (\cp y^n_{t,i}+\cu u^n_{t,i} +\ch e^n_{t,i} ) + \sum \limits_{t \in \calT} \sum \limits_{(i,j) \in \calA} c_{i,j}^{\mbox{\tiny t}} x^{n,t}_{i,j}\Big)  $$
 }
 
 \textbf{Step 2.} Compute the average of $\hat{v}_{N}^{k}$ and $\hat{v}_{N^\prime}^{k}$  among the $K$ replications
$$\overline{v}_{N}= \frac{1}{K} \sum \limits_{k=1}^{K} v_{N}^{k} \quad \quad \quad \overline{v}_{N^\prime}=\frac{1}{K} \sum \limits_{k=1}^{K}\hat{v}_{N^\prime}^{k}$$

\textbf{Step 3.} Compute the Approximate Optimality Index
 $$AOI_{N}=\frac{\overline{v}_{N^\prime}-  \overline{v}_{N}}{\overline{v}_{N^\prime}} $$ 
 
\textbf{Step 4.} If $AOI_{N}$ satisfies a predetermined termination tolerance (i.e., $|AOI_{N}| <\epsilon)$, terminate and output $N$, $ \overline{v}_{N} $, $\overline{v}_{N^\prime}$, and $AOI_{N}$. Otherwise, update $N  \leftarrow 2N$, and go to step 1.

\end{algorithm}

\noindent Starting with an initial candidate value of $N$, the algorithm\ref{alg:algorithm1} proceeds as follows. First, for $k=1,\ldots,K$, we repeat the following steps. In step 1.1, we generate a sample of $N$ i.i.d scenarios of $(d, \rho, M, V)$. In step 1.2, we solve the SAA formulation with the scenarios generated in step 1.1 and record the corresponding optimal objective value $v_{N}^{k}$ and optimal prepositioning decisions $(\pmb{\hat{o}}, \pmb{\hat{z}})_{N}^{k}$. In step 1.3, we evaluate the objective function value  $v_{N^\prime}^{k}$ via Monte Carlo simulation of $(\pmb{\hat{o}}, \pmb{\hat{z}})_{N}^{k}$ with a new sample of $N^\prime>>N$ i.i.d scenarios of $(d, \rho, M, V)$. 

In step 2, we compute the average of $v_{N}^{k}$ and $v_{N^\prime}^{k}$ among the $K$ replications as $\overline{v}_{N}= (1/K) \sum_{k=1}^{\textcolor{black}{K}} v_{N}^{\textcolor{black}{k}}$ and  $\overline{v}_{N^\prime}=(1/K) \sum_{k=1}^{K}\hat{v}_{N^\prime}^{k}$, respectively. The statistical results in \cite{mak1999monte} and \cite{linderoth2006empirical} infer that $\overline{v}_{N} $ and $\overline{v}_{N^\prime}$ are respectively statistical lower and upper bounds of the optimal value of the SP model. In step 3, we compute the approximate optimality index $|AOI_{N}= (\overline{v}_{N^\prime} -\overline{v}_{N})/\overline{v}_{N^\prime}|$ as a point estimate of the relative optimality gap between $\overline{v}_{N} $ and $\overline{v}_{N^\prime}$. Finally, if $AOI_{N}$ satisfies a predetermined termination tolerance, the algorithm terminates and outputs $N$, $ \overline{v}_{N} $, $\overline{v}_{N^\prime}$, and $AOI_{N}$. Otherwise, we increase the sample size (i.e., $N  \leftarrow 2N$), and go to step 1. This algorithm is based on the SAA method in \cite{homem2014monte} and  \cite{kleywegt2002sample} and \cite{shehadeh2020using} with some adaptations to our model.

\section{Data related to Hurricane Case Study}\label{Appx:Cas1_data}

\setcounter{table}{0}

\begin{table}[h!]  
\small
\center
   \renewcommand{\arraystretch}{0.4}
  \caption{Acquisition cost, transportation cost, and storage volume of relief supplies \citep{rawls2010pre, velasquez2020prepositioning}.}
\begin{tabular}{lcccc}
\hline 
 \textbf{Relief item} & \textbf{Acquisition cost} &	\textbf{Transportation cost}  & \textbf{Storage volume }\\
  & $c_{t}^{\mbox{\tiny a}}$ (\$/unit) & $c_{i,j}^t$ (\$/unit-mile) & $s_t$ ($\text{ft}^3/$unit) \\
  \hline 
Water  (1000 gallons) & 647.7 & 0.3 & 144.6 \\
Food (1000 meals) & 5420 & 0.04 & 83.33 \\
Medical kits & 140 & 0.00058 & 1.16\\
\hline
\end{tabular}
\label{Table:Items_Costs}
\end{table}



\begin{table}[h!]
   \small
   \center
      \renewcommand{\arraystretch}{0.4}
   \caption{Mean demand for water, food, and medical kits generated by a minor and major disaster at each potential landfall node. }
   \begin{tabular}{c l l l l l l l l l l l l l l }
     \hline 
      \textbf{Landfall}  & \multicolumn{2}{c}{\textbf{Water}} & &  \multicolumn{2}{c}{\textbf{Food}} &  & \multicolumn{2}{c}{\textbf{Medical kits}}  \\ \cline{2-3} \cline{5-6} \cline{8-9}
      \textbf{node} & minor & major  & & minor & major && minor & major \\ 
       \hline 
        2 & 500 & 2500 & & 1000 & 2000 & & 800 & 2000 \\
        5 & 500& 2000& &  500 & 1500 & & 500 & 1500\\
        11&  1500 & 7500& & 1800& 7500 & & 500& 2000\\
        13&  1000& 1500 && 500 & 9000 & & 1000 & 50000 \\
        14& 1000& 2200 & & 500 & 1500 & & 1000& 10500\\
        15&  1000 & 12000& & 1800& 4000 && 18000 & 4500\\
        21&  600 & 4000 & & 500 & 1800 & & 600 & 12500\\
        22& 1500 & 9000 & & 1500 & 4000 && 2500 & 28000\\
        29& 1500 & 7500 && 1800& 9000 & & 2000& 5000\\
        30& 1000 & 2200 & &500 & 1500 & & 1500& 10500 \\
        Total& 9600 & 47900 && 9400& 39800 &&11400& 169500\\
        \hline 
   \end{tabular}\label{Table:Demand}
\end{table}

\newpage
\clearpage

\section{MCO Convergence Results}\label{Appex:MCO_results}

\noindent For each instance, our process was as follows. For the SP model, we first optimized the sample size. We ran the MCO algorithm (Algorithm \ref{alg:algorithm1}) with $N_{o}=5$, $N^\prime$=10000, $K=20$, and $\epsilon=0.1$. These results (approximate optimality index; confidence intervals; solution times) are presented in Tables~\ref{table:SampleSize}. Based on these results, we find that $N=100$ is an appropriate sample size to obtain near-optimal solutions and tight estimates on the SP model's objective value via its SAA within a reasonable time. Then, we used this value for the Case Study 1 SP experiments.  
\begin{table}[h!]
\center 
 \footnotesize
   \renewcommand{\arraystretch}{0.4}
\caption{The Approximate Optimality Index ($AOI_{N}$) between the statistical lower bound  $\overline v_N$  and upper bound $\overline v_{N^\prime}$ on the objective values of  SP and their 95\% Confidence Interval ($95\%\text{CI}$) for each instance and each sample size, $N$} 
\begin{tabular}{llllllllllllllclllllllllllll}
\hline 
\multicolumn{6}{c}{Medium facility ($S_i$= 408,200, $f_i=\$ 188,400$)}\\ \cline{1-6}
\hline 
(nMinor, nMajor)	&	$N$	&	$95\% CI^{\tiny \bar{v}_{N}}$	&	$95\% CI^{\tiny \bar{v}_{N'}}$	&	$AOI_N$	&	Time 	\\
\hline
(2, 1)	& 5&	[48369558, 	63216812]	&	[66156537, 	71465964]	&	0.2	&	0.12\\
&	10	&	[52332546, 	60597004]	&	[63003935, 	65643255]	&	0.122	&	0.46	\\
	&	20	&	[53514528,	 60874232]	&	[61316102,	 62978588]	&	0.080	&	0.48	\\
	&	30	&	[57846678,	 63731832]	&	[61152383, 	62578997]	&	0.017	&	0.70	\\
	&	40	&	[56095622, 	61365578]	&	[60982111, 	61576539]	&	0.030	&	1.95	\\
	&	50	&	[60338000, 	61455660]	&	[60028447, 63421693]	&	0.013	&	1.00	\\
	&	100	&	[59389165, 	62447305]	&	[60424825, 	61764305]	&	0.003	&	2.30	\\
\\											
(4, 2)	&	10	&	[97907218, 	109527052]	&	[106418968, 	112067632]	&	0.051	&	0.45	\\
	&	20	&	[100012921,	107700430]	&	[106799402, 107313598]	&	0.030	&	1.00	\\
	&	30	&	[101765522, 	107163628]	&	[105463970,  106932931]	&	0.016	&	0.81	\\
	&	40	&	[103218419,	108630911]	&	[106395116, 	107607984]	&	0.010	&	1.74	\\
	&	50	&	[103334694,	 107666076]	&	[106642455,	 107628645]	&	0.015	&	1.97	\\
	&	100	&	[104643866, 108416344]	&	[106037279, 107158322]	&	0.001	&	222	\\
\\											
(6, 3)	&	10	&	[139350470,  151548430]	&	[150611748, 	152423052]	&	0.040	&	0.30	\\
	&	20	&	[142104626, 	149846274]	&	[149599196,  151369804]	&	0.030	&	1.34	\\
	&	30	&	[145259753,	 152343848]	&	[148963381,	 149576219]	&	0.010	&	2.15	\\
	&	40	&	[146396483,	 152199917]	&	[149222838, 	150526362]	&	0.010	&	2.37	\\
	&	50	&	[147858180, 151750721]	&	[148808189,	149405411]	&	0.010	&	2.50	\\
	&	100	&	[146965469, 	149626531]	&	[148128318, 149282882]	&	0.003	&	5.00	\\
\hline
\multicolumn{6}{c}{Large facility ($S_i$= 780,000, $f_i=\$ 300,000$)}\\ \cline{1-6}
\hline 
(nMinor, nMajor)	&	$N$	&	$95\% CI^{\tiny \bar{v}_{N}}$	&	$95\% CI^{\tiny \bar{v}_{N'}}$	&	$AOI_N$	&	Time 	\\
\hline
(2, 1)	&	5	&	[48475405, 63281824]	&	[66550531, 1875919]	&	0.2	&	0.1	\\
	&	10	&	[55213854, 65592216]	&	[63311092, 	66898918]	&	0.08	&	0.2	\\
	&	20	&	[55033183,	 60866807]	&	[61458118, 	62507202]	&	0.07	&	0.3	\\
	&	30	&	[58157970, 	64111460]	&	[61024550, 	62981960]	&	0.03	&	2	\\
	&	40	&	[56483883,  63032347]	&	[61295188,	 62350282]	&	0.04	&	1	\\
	&	50	&	[56723273,  61175247]	&	[60107528,  61222862]	&	0.03	&	2	\\
	&	100	&	[57822125, 	60857705]	&	[59793941, 	61132239]	&	0.01	&	2	\\
\\											
(4, 2)	&	5	&	[98556082, 113502488]	&	[112644835, 122872265]	&	0.1	&	0.12	\\
	&	10	&	[99145976, 109361034]	&	[108198746, 108696354]	&	0.04	&	0.27	\\
	&	20	&	[100432385, 	106943995]	&	[105403174,  106504826]	&	0.03	&	1	\\
	&	30	&	[102573753, 106162307]	&	[104801669, 	105958731]	&	0.02	&	1	\\
	&	40	&	[104457581, 	109726999]	&	[106129241, 	106975560]	&	0.02	&	1.2	\\
	&	50	&	[99486180, 104504499]	&	[105603202, 	106813398]	&	0.04	&	1.2	\\
	&	100	&	[105332261, 106715439]	&	[105707875, 	107434125]	&	0.003	&	4	\\
\\											
(6, 3)	&	5	&	[130143459, 	141601041]	&	[149867389, 152337311]	&	0.11	&	0.1	\\
	&	10	&	[140346961, 149372139]	&	149414392, 	152297508]	&	0.03	&	0.1	\\
	&	20	&	[142323065, 149157035]	&	[148611699,  	149925100]	&	0.03	&	0.4	\\
	&	30	&	[143941057,	 149834843]	&	[147506972, 	149107128]	&	0.02	&	0.8	\\
	&	40	&	[141654411,  148184289]	&	[147870704, 	149179596]	&	0.02	&	1	\\
	&	50	&	[145460556,  149499644]	&	[147500429, 	148996972]	&	0.001	&	1	\\
	&	100	&	[145694430,  149487569]	&	[147597746,	148355253]	&	0.001	&	2.2	\\
\hline 											
\end{tabular}\label{table:SampleSize}
\end{table}

\section{Data related to Earthquake Case Study}\label{Appx:Case2_data}

\begin{figure}[h!]
    \includegraphics[width=\textwidth]{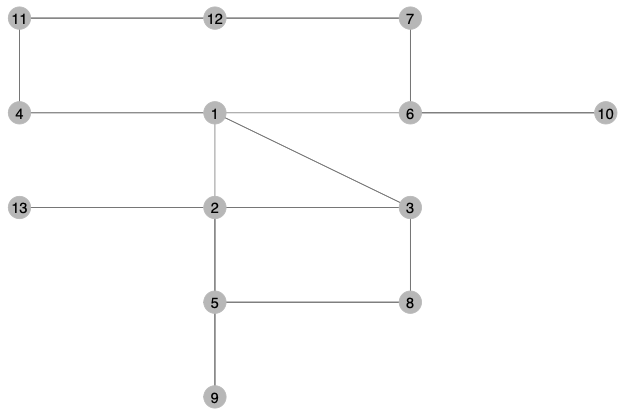}
    \caption{Map of the facility and transportation network.}
    \label{Fig_map2}
\end{figure}
\setcounter{table}{0}

\begin{table}[h!]
 \center 
 \footnotesize
   \renewcommand{\arraystretch}{0.5}
  \caption{Input Parameters for Each Node \citep{ni2018location}.}
\begin{tabular}{lllllllllllllllllllllll}
\hline
Node & 1 & 2& 3 & 4 & 5 & 6& 7 & 8 & 9& 10& 11&12&13 \\
\hline
$f_i$ & 203&	193&	130&	117&	292	&174	&130	&157	&134&	161	&234&	220	1&70\\ 
$c^{\mbox{\tiny a}}_i$ &3.4	&2.33&	2	&2.69&	2.63	&3.44	&3.43&	3.53&	2.33&	2.5&	3.37&	2.84	&3.76 \\
$c^{\mbox{\tiny h }}_i$& 2.81&2.58	 &2.86&	2.42	&3.28&	3.05&	2.77& 2.68	& 2.52&	3.14	& 2.93	&2.85	 &2.87 \\
$c^{\mbox{\tiny u }}_i$ &  11.48	&  14.32&  	12.14&  	16.19	&  12.01&  	14.9	&  9.42	&  11.91	&  10.68&  	11.24&  	13.1	&  11.09&  	10.1 \\
 $\mu^{\tiny \rho}$ & 0.05&	0.05	&0.2	&0.18	&0.18	&0.72&	0.76	&0.7&0.6	&0.7	&0.78	&0.62&	0.68\\
\hline
	\end{tabular} 
\label{table:Case2_input}
\end{table}

\newpage

\section{Computational Performance of DRO--Decomposition Algorithm}\label{Appix:DROPerform}
\setcounter{table}{0}

\begin{table}[h!]
\color{black}
 \center 
 \footnotesize
   \renewcommand{\arraystretch}{0.6}
  \caption{Computational details of solving larger networks. (nMinor, nMajor)=(2,1).}
\begin{tabular}{cclllllllllll}
\hline
$I$	&	Nodes	&	\# of Iter	&	Time	 (s)&	\# of B\&B  	&	\# of MIPiter 	\\
\hline
40	&	10	&	62	&	12	&	865	&	101199	\\
	&	15	&	103	&	29	&	19526	&	809218	\\
	&	20	&	107	&	50	&	91559	&	3410613	\\
											\\
60	&	10	&	87	&	22	&	1515	&	207893	\\
	&	15	&	114	&	43	&	18986	&	1005776	\\
	&	20	&	172	&	103	&	160737	&	7259463	\\
											\\
100	&	10	&	23	&	47	&	4046	&	348607	\\
	&	15	&	42	&	74	&	27077	&	1918668	\\
	&	20	&	53	&	123	&	126004	&	8321668	\\											
\hline
	\end{tabular} 
\label{table:CPU_2_1}
\end{table}
\begin{table}[h!]
\color{black}
 \center 
 \footnotesize
   \renewcommand{\arraystretch}{0.6}
  \caption{Computational details of solving larger networks. (nMinor, nMajor)=(4,2).}
\begin{tabular}{cclllllllllll}
\hline
$I$	&	Nodes	&	\# of Iter		&	Time	 (s) &	\# of B\&B  	&	\# of MIPiter 	\\
\hline
40	&	10	&	81	&	16	&	86	&	111400	\\
	&	15	&	111	&	27	&	9284	&	731560	\\
	&	20	&	162	&	52	&	15783	&	1248777	\\
											\\
60	&	10	&	85	&	25	&	430	&	179054	\\
	&	15	&	118	&	45	&	6661	&	735994	\\
	&	20	&	172	&	133	&	306238	&	20909011	\\
											\\
100	&	10	&	93	&	53	&	309	&	229197	\\
	&	15	&	157	&	500	&	1336	&	552646	\\
	&	20	&	176	&	359	&	108485	&	8245269	\\
\hline
	\end{tabular} 
\label{table:CPU_4_2}
\end{table}

\begin{table}[h!]
\color{black}
 \center 
 \footnotesize
   \renewcommand{\arraystretch}{0.6}
  \caption{Computational details of solving larger networks. (nMinor, nMajor)=(6,3).}
\begin{tabular}{cclllllllllll}
\hline
$I$	&	Nodes	&	\# of Iter		&	Time (s)	&	\# of B\&B  	&	\# of MIPiter 	\\
\hline
40	&	10	&	67	&	12	&	18	&	76817	\\
	&	15	&	102	&	26	&	2006	&	298344	\\
	&	20	&	172	&	58	&	8584	&	845388	\\
											\\
60	&	10	&	64	&	17	&	240	&	135853	\\
	&	15	&	130	&	47	&	1177	&	394397	\\
	&	20	&	206	&	94	&	7020	&	1161068	\\
											\\
100	&	10	&	4	&	41	&	10	&	162069	\\
	&	15	&	176	&	644	&	213	&	532303	\\
	&	20	&	207	&	843	&	8351	&	1575093	\\
\hline
	\end{tabular} 
\label{table:CPU_6_3}
\end{table}

\color{black}

\end{document}